\documentclass[a4paper,twoside]{book}

\newcommand{\heuteIst}{November 29, 2001 }

\usepackage{amsmath,amsthm,amsfonts,latexsym,amscd,amssymb,enumerate}
\newcommand{\href}[2]{#2}

\swapnumbers
\theoremstyle{plain}
\newtheorem{theorem}{Theorem}[section]
\newtheorem{lemma}[theorem]{Lemma}
\newtheorem{corollary}[theorem]{Corollary}
\newtheorem{proposition}[theorem]{Proposition}
\newtheorem{conjecture}[theorem]{Conjecture}

\theoremstyle{definition}
\newtheorem{definition}[theorem]{Definition}
\newtheorem{example}[theorem]{Example}

\theoremstyle{remark}
\newtheorem{remark}[theorem]{Remark}
\newtheorem{question}[theorem]{Question}

\newcommand{\reals}{\mathbb{R}}
\newcommand{\complexs}{\mathbb{C}}
\newcommand{\naturals}{\mathbb{N}}
\newcommand{\integers}{\mathbb{Z}}
\newcommand{\rationals}{\mathbb{Q}}
\newcommand{\quaternions}{\mathbb{H}}
\newcommand{\hyperbolic}{\mathbb{H}}

\DeclareMathOperator{\id}{id}

\newcommand{\boundedops}{\mathcal{B}}

\newcommand{\abs}[1]{\left\lvert#1\right\rvert} 

\newcommand{\tensor}{\otimes}
\newcommand{\into}{\hookrightarrow}
\newcommand{\onto}{\twoheadrightarrow}
\newcommand{\iso}{\cong}

\newcommand{\subgroup}{\leq}
\newcommand{\supergroup}{\geq}

\DeclareMathOperator{\supp}{supp}   
\DeclareMathOperator{\im}{im}      
\DeclareMathOperator{\coker}{coker}      
\DeclareMathOperator{\vol}{vol}    
\DeclareMathOperator{\End}{End}    
\DeclareMathOperator{\Hom}{Hom}    

\DeclareMathOperator{\tr}{tr}
\DeclareMathOperator{\pr}{pr}

\DeclareMathOperator{\ch}{ch}  
\DeclareMathOperator{\Td}{Td}  
\DeclareMathOperator{\ind}{ind}
\DeclareMathOperator{\sign}{sign}
\DeclareMathOperator{\Ric}{Ric}  
\DeclareMathOperator{\scal}{scal}  
\DeclareMathOperator*{\dirlim}{\varinjlim}

\newcommand{\forget}[1]{}

\newcommand{\innerprod}[1]{\langle #1 \rangle}

{\catcode`@=11\global\let\c@equation=\c@theorem}

\allowdisplaybreaks[2]
\newcommand{\RAgroups}{\mathcal G}

\newcommand{\extendedC}{\mathcal{D}}

\newcommand{\NeumannN}{\mathcal{N}}
\newcommand{\universalU}{\mathcal{U}}
\newcommand{\amenableGroups}{{\mathcal{Y}}}

\begin{document}
\date{Last compiled \today; last edited  \heuteIst or later}

\title{Operator algebras and topology}
\author{Thomas Schick\thanks{
\protect\href{mailto:schick@uni-math.gwdg.de}{e-mail: schick@uni-math.gwdg.de}
\protect\\
\protect\href{http://uni-math.gwdg.de/schick}{www:~http://uni-math.gwdg.de/schick}
\protect\\
Fax: ++49 -251/83 38370
}\\
}
        
\maketitle

  These notes, based on three lectures on operator algebras and
  topology at the ``School on High Dimensional Manifold Theory'' at the
  ICTP in Trieste, introduce a new
  set of tools to high dimensional manifold theory, namely techniques
  coming from the theory of operator algebras, in particular
  $C^*$-algebras. These are extensively
  studied in their own right. We will focus on the basic definitions
  and properties, and on their relevance to the geometry and topology
  of manifolds.

  A central pillar of work in the theory of $C^*$-algebras is the
  Baum-Connes conjecture. This is an isomorphism conjecture, as
  discussed in the talks of L\"uck, but with a certain special
  flavor. Nevertheless, it has important direct applications to the
  topology of manifolds, it implies e.g.~the Novikov conjecture. In the
  first chapter, the Baum-Connes conjecture will be explained and put
  into our context.

  Another application of the Baum-Connes conjecture is to the positive
  scalar curvature question. This will be discussed by Stephan
  Stolz. It implies the so-called ``stable Gromov-Lawson-Rosenberg
  conjecture''. The unstable version of this conjecture said that,
  given a closed spin manifold $M$, a
  certain obstruction, living in a certain (topological) $K$-theory
  group, vanishes if and only $M$ admits a Riemannian metric with
  positive scalar curvature. It turns out that this is wrong, and
  counterexamples will be presented in the second chapter.

  The third chapter introduces another set of invariants, also using
  operator algebra techniques, namely $L^2$-cohomology, $L^2$-Betti
  numbers and other $L^2$-invariants. These invariants, their basic
  properties, and the central questions about them, are introduced 
  in the third chapter.

  Several people contributed to these notes by reading preliminary
  parts and suggesting improvements, in particular Marc Johnson, Roman
  Sauer, Marco Varisco und Guido Mislin. I am very indebted to all of
  them.




\chapter[Index theory and Baum-Connes]{Index theory and the Baum-Connes conjecture}

\section{Index theory}
\label{sec:index}

The Atiyah-Singer index theorem is one of the great achievements of
modern mathematics. It gives a formula for the index of a differential
operator (the index is by definition the dimension of the space of its
solutions minus the
dimension of the solution space for its adjoint operator) in terms
only of topological data associated to the operator and the underlying 
space. There are many good treatments of this subject available, apart 
from the original literature (most found
in~\cite{Atiyah-Collected2}). Much more detailed than the present
notes can be, because of constraints of length and time, are
e.g.~\cite{Lawson-Michelsohn(1989),Berline-Getzler-Vergne(1992),Higson-Roe(2001)}.

\subsection{Elliptic operators and their index}

We quickly review what type of operators we are looking at.

\begin{definition}
  Let $M$ be a smooth manifold of dimension $m$; $E,F$ smooth
  (complex) vector bundles
  on $M$. A \emph{differential operator} (of order $d$) from $E$ to
  $F$ is a $\complexs$-linear map from the space of smooth sections
  $C^\infty(E)$ of $E$ to the space of smooth sections of $F$:
  \begin{equation*}
    D\colon C^\infty(E)\to C^{\infty}(F),
  \end{equation*}
  such that in local coordinates and with local trivializations of the 
  bundles it can be written in the form
  \begin{equation*}
    D= \sum_{\abs{\alpha}\le d} A_\alpha(x)
    \frac{\partial^{\abs{\alpha}}}{\partial x^\alpha}.
  \end{equation*}
  Here $A_\alpha(x)$ is a matrix of smooth complex valued functions,
  $\alpha=(\alpha_1,\dots,\alpha_m)$ is an $m$-tuple of non-negative
  integers and
  $\abs{\alpha}=\alpha_1+\dots+\alpha_m$.
  $\partial^{\abs{\alpha}}/\partial x^\alpha$ is an abbreviation for
  $\partial^{\abs{\alpha}}/\partial x_1^{\alpha_1}\cdots\partial
  x_m^{\alpha_m}$. We require that 
  $A_\alpha(x)\ne 0$ for some $\alpha$ with $\abs{\alpha}=d$ (else,
  the operator is of order strictly smaller than $d$).


Let $\pi\colon T^*M\to M$ be the bundle projection of the cotangent
bundle of $M$. We get pull-backs $\pi^*E$ and $\pi^*F$ of the bundles
$E$ and $F$, respectively, to $T^*M$.

The \emph{symbol} $\sigma(D)$ of the differential operator $D$ is the
section of the bundle
$\Hom(\pi^*E,\pi^* F)$ on $T^*M$ defined as follows:

In the above local coordinates, using $\xi=(\xi_1,\dots,\xi_m)$ as coordinate for the
cotangent vectors in $T^*M$,
in the fiber of $(x,\xi)$, the symbol $\sigma(D)$ is given by
multiplication with
\begin{equation*}
  \sum_{\abs{\alpha}=m} A_{\alpha}(x) \xi^\alpha.
\end{equation*}
Here $\xi^\alpha=\xi_1^{\alpha_1}\cdots \xi_m^{\alpha_m}$.

The operator $D$ is called \emph{elliptic}, if
$\sigma(D)_{(x,\xi)}\colon \pi^*E_{(x,\xi)}\to \pi^*F_{(x,\xi)}$ is
invertible outside the zero section of $T^*M$, i.e.~in each fiber over 
$(x,\xi)\in T^*M$ with $\xi\ne 0$. Observe that elliptic operators can 
only exist if the fiber dimensions of $E$ and $F$ coincide.

In other words, the symbol of an elliptic operator gives us two vector 
bundles over $T^*M$,  namely $\pi^*E$ and $\pi^*F$, together with a
choice of an isomorphism of the fibers of these two bundles outside
the zero section. If $M$ is compact, this gives an element of the
relative $K$-theory group $K^0(DT^*M,ST^*M)$, where $DT^*M$ and
$ST^*M$ are the disc bundle and sphere bundle of $T^*M$, respectively
(with respect to some arbitrary Riemannian metric).
\end{definition}

Recall the following definition:
\begin{definition}
  Let $X$ be a compact topological space. We define the $K$-theory of
  $X$, $K^0(X)$, to be the Grothendieck group of (isomorphism classes
  of) complex vector bundles over $X$ (with finite fiber
  dimension). More precisely, $K^0(X)$ consists of equivalence classes 
  of pairs $(E,F)$ of
  (isomorphism classes of) vector bundles over $X$, where $(E,F)\sim
  (E',F')$ if and only if there exists another vector bundle $G$ on
  $X$ such that $E\oplus F'\oplus G\iso E'\oplus F\oplus G$. One often 
  writes $[E]-[F]$ for the element of $K^0(X)$ represented by $(E,F)$.

  Let $Y$ now be a closed subspace of $X$. The \emph{relative
    $K$-theory} $K^0(X,Y)$ is given by equivalence classes of triples
  $(E,F,\phi)$, where $E$ and $F$ are complex vector bundles over $X$, 
  and $\phi\colon E|_Y\to F|_Y$ is a given isomorphism between the
  restrictions of $E$ and $F$ to $Y$. Then $(E,F,\phi)$ is isomorphic to
  $(E',F',\phi')$ if we find isomorphisms $\alpha\colon E\to E'$ and
  $\beta\colon F\to F'$ such that the following diagram commutes.
  \begin{equation*}
    \begin{CD}
      E|_Y @>{\phi}>> F|_Y\\
      @VV{\alpha}V @VV{\beta}V\\
      E'|_Y @>{\phi'}>> F'|_{Y}
    \end{CD}
  \end{equation*}
  Two pairs $(E,F,\phi)$ and $(E',F',\phi')$ are equivalent, if there
  is a bundle $G$ on $X$ such that $(E\oplus G,F\oplus
  G,\phi\oplus\id)$ is isomorphic to $(E'\oplus G,F'\oplus
  G,\phi'\oplus \id)$. 
\end{definition}

\begin{example}
  The element of $K^0(DT^*M,ST^*M)$ given by the symbol of an elliptic 
  differential operator $D$ mentioned above is represented by the
  restriction of the bundles $\pi^*E$ and $\pi^*F$ to the disc bundle
  $DT^*M$, together with the isomorphism $\sigma(D)_{(x,\xi)}\colon
  E_{(x,\xi)}\to F_{(x,\xi)}$ for $(x,\xi)\in ST^*M$.
\end{example}

\begin{example}
  Let $M=\reals^m$ and $D=\sum_{i=1}^m (\partial/\partial_i)^2$ be the 
  Laplace operator on functions. This is an elliptic differential
  operator, with symbol $\sigma(D)=\sum_{i=1}^m \xi_i^2$.

  More generally, a second-order differential operator $D\colon C^\infty(E)\to
  C^\infty(E)$ on a Riemannian manifold $M$ is a \emph{generalized
    Laplacian}, if
  $\sigma(D)_{(x,\xi)} = \abs{\xi}^2\cdot\id_{E_x}$ (the norm of the
  cotangent vector $\abs{\xi}$ is
  given by the Riemannian metric).

  Notice that all generalized Laplacians are elliptic.
\end{example}

\begin{definition}
  \emph{(Adjoint operator)}\\
  Assume that we have a differential operator $D\colon C^\infty(E)\to
  C^\infty(F)$ between two Hermitian bundles $E$ and $F$ on a
  Riemannian manifold $(M,g)$. We define an $L^2$-inner product on
  $C^\infty(E)$ by the formula 
  \begin{equation*}
\innerprod{f,g}_{L^2(E)} := \int_M \innerprod{f(x),g(x)}_{E_x}
\;d\mu(x)\qquad\forall f,g\in C^\infty_0(E),
\end{equation*}
  where $\innerprod{\cdot,\cdot}_{E_x}$ is the fiber-wise inner product 
  given by the Hermitian metric, and $d\mu$ is the measure on $M$
  induced from the Riemannian metric. Here $C^\infty_0$ is the space of
  smooth section with compact support.
  The Hilbert space completion of $C^\infty_0(E)$ with respect to this 
  inner product is called $L^2(E)$.

   The \emph{formal adjoint} $D^*$ of $D$ is then defined by
   \begin{equation*}
     \innerprod{Df,g}_{L^2(F)} =
     \innerprod{f,D^*g}_{L^2(E)}\qquad\forall f\in C^\infty_0(E),\;
       g\in C^\infty_0(F).
     \end{equation*}
   It turns out that exactly one operator with this property exists,
   which is another differential operator, and which is elliptic if
   and only if $D$ is elliptic.
\end{definition}

\begin{remark}
  The class of differential operators is quite restricted. Many
    constructions one would like to carry out with differential
    operators automatically lead out of this class. Therefore, one
    often has to use \emph{pseudodifferential operators}. 
    Pseudodifferential operators are defined as a generalization of
  differential operators. There are many well written sources dealing with the theory of
pseudodifferential operators. Since we will not discuss them in detail 
here, we omit even their precise definition and refer e.g.~to
\cite{Lawson-Michelsohn(1989)} and \cite{Shubin(1987)}.
  What we have done so far with elliptic
  operators can all be extended to pseudodifferential operators. In
  particular, they have a symbol, and the concept of ellipticity is
  defined for them. When studying elliptic differential operators,
  pseudodifferential operators naturally appear and play a very
  important role. An pseudodifferential operator $P$  (which could
    e.g.~be a differential operator) is elliptic if
  and only if a pseudodifferential operator $Q$ exists such that
  $PQ-\id$ and $QP-\id$ are so called \emph{smoothing} operators, a
  particularly nice class of pseudodifferential operators. For many
  purposes, $Q$ can be considered to act like an inverse of $P$, and
  this kind of invertibility is frequently used in the theory of elliptic
  operators. However, if $P$ happens to be an elliptic differential
  operator of positive order, then $Q$ necessarily is not a
  differential operator, but only a pseudodifferential operator.

It should be noted that almost all of the results we present here for
differential operators hold also for pseudodifferential operators, and 
often the proof is best given using them.
\end{remark}

We now want to state several important properties of elliptic
operators.

\begin{theorem}
  Let $M$ be a smooth manifold, $E$ and $F$ smooth finite dimensional
  vector bundles over $M$.
  Let $P\colon C^\infty(E)\to C^\infty(F)$ be an elliptic operator.

  Then the following holds.
  \begin{enumerate}
  \item Elliptic regularity:\\
    If $f\in L^2(E)$ is weakly in the null space of $P$,
    i.e.~$\innerprod{f,P^*g}_{L^2(E)}=0$ for all $g\in C^\infty_0(F)$, 
    then $f\in C^\infty(E)$.
  \item Decomposition into finite dimensional eigenspaces:\\
    Assume $M$ is compact and $P=P^*$ (in particular, $E=F$). Then the 
    set $s(P)$ of eigenvalues of $P$ ($P$ acting on $C^\infty(E)$) is 
    a discrete subset of $\reals$, each eigenspace $e_\lambda$
    ($\lambda\in s(P)$) is finite dimensional, and
    $L^2(E)=\oplus_{\lambda\in s(P)} e_\lambda$ (here we use the
    completed direct sum in the sense of Hilbert spaces, which means
    by definition that the
    algebraic direct sum is dense in $L^2(E)$).
  \item If $M$ is compact, then $\ker(P)$ and $\ker(P^*)$ are finite
    dimensional, and then we define the \emph{index of $P$}
    \begin{equation*}
      \ind(P):=\dim_\complexs \ker(P) - \dim_\complexs \ker(P^*).
    \end{equation*}
  \end{enumerate}
  (Here, we could replace $\ker(P^*)$ by $\coker(P)$, because these
  two vector spaces are isomorphic).
\end{theorem}

\subsection{Statement of the Atiyah-Singer index theorem}

There are different variants of the Atiyah-Singer index theorem. We
start with a cohomological formula for the index.

\begin{theorem}
  Let $M$ be a compact oriented manifold of dimension $m$, and $D\colon C^\infty(E)\to
  C^{\infty}(F)$ an elliptic operator with symbol $\sigma(D)$. There
  is a characteristic (inhomogeneous) cohomology class $\Td(M)\in
  H^*(M;\rationals)$ of the tangent bundle of $M$ (called the complex
  Todd class of the complexified tangent bundle). Moreover, to the 
  symbol is associated a certain (inhomogeneous) cohomology class
  $\pi_!\ch(\sigma(D))\in H^*(M;\rationals)$ such that
  \begin{equation*}
    \ind(D) = (-1)^{m(m+1)/2} \innerprod{\pi_!\ch(\sigma(D)) \cup \Td(M), [M]} .
  \end{equation*}
  The class $[M]\in H_m(M;\rationals)$ is the fundamental class of the 
  oriented manifold $M$, and $\innerprod{\cdot,\cdot}$ is the usual
  pairing between homology and cohomology.
\end{theorem}

If we start with specific operators
given by the geometry, explicit calculation usually give more familiar
terms on the right hand
side.

 For example, for the signature operator we obtain Hirzebruch's
signature formula expressing the signature in terms of the $L$-class,
for the Euler characteristic operator we obtain the
Gauss-Bonnet formula expressing the Euler characteristic in terms of
the Pfaffian, and for the spin or spin$^c$ Dirac operator we obtain an 
$\hat{A}$-formula. For applications, these formulas prove to be
particularly useful.

We give some more details about the signature operator, which we are
going to use later again. To define the signature operator, fix a
Riemannian metric $g$ on $M$. Assume $\dim M=4k$ is divisible by four.

The signature operator maps from a certain subspace $\Omega^+$ of the space of
differential forms to another subspace $\Omega^-$. These subspaces are 
defined as follows. Define, on $p$-forms, the operator $\tau:=
i^{p(p-1)+2k}*$, where $*$ is the Hodge-$*$ operator given by the
Riemannian metric, and $i^2=-1$. Since $\dim M$ is divisible by $4$, an
easy calculation
shows that $\tau^2=\id$. We then define $\Omega^{\pm}$ to be the
$\pm1$ eigenspaces of $\tau$.

The signature operator $D_{sig}$ is now simply defined to by
$D_{sig}:= d+d^*$, where $d$ is the exterior derivative on
differential forms, and $d^*=\pm *d*$ is its formal adjoint. We restrict
this operator to $\Omega^+$, and another easy calculation shows that
$\Omega^+$ is mapped to $\Omega^-$. $D_{sig}$ is elliptic, and a
classical calculation shows that its index is the signature of $M$
given by the intersection form in middle homology.

\begin{definition}
  The \emph{Hirzebruch L-class} as normalized by Atiyah and Singer is
  an inhomogeneous characteristic class, assigning to each complex vector
  bundle $E$ over a space $X$ a cohomology class $L(E)\in
  H^*(X;\rationals)$. It is characterized by the following properties:
  \begin{enumerate}
  \item Naturality: for any map $f\colon Y\to X$ we have
    $L(f^*E)=f^*L(E)$.
  \item Normalization: If $L$ is a complex line bundle with first
    Chern class $x$, then
    \begin{equation*}
      L(E) = \frac{x/2}{\tanh (x/2)} = 1 +\frac{1}{12} x^2 -
      \frac{1}{720} x^4 +\cdots \in H^*(X;\rationals).
    \end{equation*}
  \item Multiplicativity: $L(E\oplus F) = L(E)L(F)$.
  \end{enumerate}
\end{definition}
It turns out that $L$ is a \emph{stable} characteristic class,
i.e.~$L(E)=1$ if $E$ is a trivial bundle. This implies that $L$
defines a map from the K-theory $K^0(X)\to H^*(X;\rationals)$.

The Atiyah-Singer index theorem now specializes to
\begin{equation*}
  \sign(M) = \ind(D_{sig}) = \innerprod{ 2^{2k} L(TM),[M]},
\end{equation*}
with $\dim M=4k$ as above.

\begin{remark}
  One direction to generalize the Atiyah-Singer index theorem is to
  give an index formula for manifolds with boundary. Indeed, this is
  achieved in the Atiyah-Patodi-Singer index theorem. However, these
  results are much less topological than the results for manifolds
  without boundary. They are not discussed in these notes.
\end{remark}

Next, we explain the K-theoretic version of the Atiyah-Singer index
theorem. It starts with the element of $K^0(DT^*M,ST^*M)$ given by the 
symbol of an elliptic operator. Given any compact manifold $M$, there is a
well defined
homomorphism 
\begin{equation*}
K^0(DT^*M,ST^*M)\to K^0(*)=\integers,
\end{equation*}
constructed by
embedding $T^*M$ into high dimensional Euclidean space, then using a
transfer map and Bott periodicity. The image of the symbol element
under this homomorphism is denoted the \emph{topological index}
$\ind_t(D)\in K^0(*)=\integers$. The reason for the terminology is
that it is obtained from the symbol only,
using purely topological constructions. Now, the Atiyah-Singer index
theorem states
\begin{theorem}
  $ \ind_t(D) =\ind(D)$.
\end{theorem}

\subsection{The $G$-index}\label{sec:equiv-index}

Let $G$ be a finite group, or more generally a compact Lie group.
The representation ring $RG$ of $G$ is defined to be the Grothendieck
group of all finite dimensional complex representations of $G$,
i.e.~an element of $RG$ is a formal difference $[V]-[W]$ of two finite 
dimensional $G$-representations $V$ and $W$, and we have
$[V]-[W]=[X]-[Y]$ if and only if $V\oplus Y\iso W\oplus X$ (strictly
speaking, we have to pass to isomorphism classes of representations to 
avoid set theoretical problems). The direct sum of representations
induces the structure 
of an abelian group on $RG$, and the tensor product makes it a
commutative unital ring (the unit given by the trivial one-dimensional 
representation). More about this representation ring can be found
e.g.~in \cite{Broecker-Dieck(1995)}.

Assume now that the manifold $M$ is a compact smooth manifold with a smooth
$G$-action, and let $E,F$ be complex $G$-vector bundles on $M$ (this
means that $G$ acts on $E$ and $F$ by vector bundle automorphisms
(i.e.~carries fibers to fibers linearly), and the bundle projection
maps are $G$-equivariant).

Let $D\colon C^\infty(E)\to C^\infty(F)$ be a $G$-equivariant elliptic 
differential operator.

This implies that $\ker(D)$ and $\coker(D)$ inherit a $G$-action by
restriction, i.e.~are finite dimensional $G$-representations. We
define the (analytic) $G$-index of $D$ to be
\begin{equation*}
  \ind^G(D):= [\ker(D)] - [\coker(D)] \in RG.
\end{equation*}

If $G$ is the
trivial group then $RG \iso \integers$ in a canonical way, and then
$\ind^G(D)$ coincides with the usual index of
$D$.

We can also define a \emph{topological}
equivariant index similar to the non-equi\-va\-riant topological index,
using transfer maps and Bott periodicity. This topological index lives
in the $G$-equivariant $K$-theory of a
point, which is canonically isomorphic to the representation ring
$RG$. Again, the Atiyah-Singer index theorem says
\begin{theorem}\label{theo:equivariant_Atiyah_Singer}
  $\ind^G(D) = \ind^G_t(D) \in K_G^0(*) = RG.$
\end{theorem}

\subsection{Families of operators and their index}

Another generalization is given if we don't look at one operator on
one manifold, but a family of operators on a family of manifolds. More 
precisely, let $X$ be any compact topological space, $Y\to X$ a
locally trivial fiber bundle with fiber $M$ a smooth compact manifold,
and structure group the diffeomorphisms of $M$. Let $E,F$ be families
of smooth vector bundles on $Y$ (i.e.~vector bundles which are
fiber-wise smooth), and $C^\infty(E)$, $C^\infty(F)$ the continuous
sections which are smooth along the fibers. Assume that $D\colon
C^\infty(E)\to C^\infty(F)$ is a family $\{D_x\}$ of elliptic differential
operator along the fiber $Y_x\iso M$ ($x\in X$), i.e., in local
coordinates $D$ becomes 
\begin{equation*}
\sum_{\abs{\alpha}\le m}
A_\alpha(y,x)\frac{\partial^{\abs{\alpha}}}{\partial y^\alpha}
\end{equation*}
with $y\in M$ and $x\in X$ such that $A_\alpha(y,x)$ depends
continuously on $x$, and each $D_x$ is an elliptic differential
operator on $Y_x$.

If $\dim_\complexs{\ker(D_x)}$ is independent of $x\in X$, then all of
these vector spaces patch together to give a vector bundle called
$\ker(D)$ on $X$, and similarly for the (fiber-wise) adjoint
$D^*$. This then gives a $K$-theory element $[\ker(D)]-[\ker(D^*)]\in
K^0(X)$.

Unfortunately, it does sometimes happen that these dimensions
jump. However, using appropriate
perturbations, one can always define the K-theory element 
\begin{equation*}
  \ind(D):= [\ker(D)]-[\ker(D^*)] \in K^0(X),
\end{equation*}
the analytic index of the family of elliptic operators $D$. 

There is also a family version of the construction of the topological
index, giving $\ind_t(D)\in K^0(X)$. The Atiyah-Singer index theorem
for families states:
\begin{theorem}
  $\ind(D)=\ind_t(D)\in K^0(X)$.
\end{theorem}

The upshot of the discussion of this and the last section (for the
details the reader is referred to the literature) is that the
natural receptacle for the index of differential operators in various
situations are appropriate K-theory groups, and much of todays index
theory deals with investigating these K-theory groups.

\section{Survey on $C^*$-algebras and their $K$-theory}

More detailed references for this section are, among others,
\cite{Wegge-Olsen(1993)}, \cite{Higson-Roe(2001)}, and \cite{Blackadar(1998)}.

\subsection{$C^*$-algebras}

\begin{definition}
  A \emph{Banach algebra} $A$ is a complex algebra which is a complete normed
  space, and such that $\abs{ab}\le \abs{a}\abs{b}$ for each $a,b\in A$.

  A \emph{$*$-algebra} $A$ is a complex algebra with an anti-linear
  involution $*\colon A\to
  A$ (i.e.~$(\lambda a)^*=\overline\lambda a^*$, $(ab)^*=b^*a^*$, and
  $(a^*)^*=a$ for
  all $a,b\in A$).

  A \emph{Banach $*$-algebra} $A$ is a Banach algebra which is a
  $*$-algebra such that   $\abs{a^*}=\abs{a}$ for
  all $a\in A$. 

  A \emph{$C^*$-algebra} $A$ is a Banach $*$-algebra
  which satisfies $\abs{a^*a}=\abs{a}^2$ for all $a\in A$.

  Alternatively, a $C^*$-algebra is a Banach $*$-algebra which is
  isometrically $*$-isomorphic to a norm-closed subalgebra of the
  algebra of bounded operators on some Hilbert space $H$ (this is the
  Gelfand-Naimark representation theorem, compare
  e.g.~\cite[1.6.2]{Higson-Roe(2001)}).

  A $C^*$-algebra $A$ is called separable if there exists a countable
  dense subset of $A$.
\end{definition}

\begin{example}\label{ex:loc_comp_Cstar_algebra}
  If $X$ is a compact topological space, then $C(X)$, the algebra of
  complex valued continuous functions on $X$, is a commutative
  $C^*$-algebra (with unit). The adjoint is given by complex
  conjugation: $f^*(x) = \overline{f(x)}$, the norm is the
  supremum-norm.

  Conversely, it is a theorem that every abelian unital $C^*$-algebra
  is isomorphic to $C(X)$ for a suitable compact topological space
  $X$ \cite[Theorem 1.3.12]{Higson-Roe(2001)}.

  Assume $X$ is locally compact, and set
  \begin{equation*}
C_0(X):=\{f\colon X\to\complexs\mid f\text{ continuous},
f(x)\xrightarrow{x\to\infty} 0\}.
\end{equation*}
Here, we say $f(x)\to 0$ for $x\to\infty$, or\emph{$f$ vanishes at
  infinity}, if for all $\epsilon >0$ there is a compact subset $K$ of 
$X$ with $\abs{f(x)}<\epsilon$ whenever $x\in X-K$.
This is again a commutative $C^*$-algebra (we use the supremum norm
on $C_0(X)$), and it is unital if and only if $X$ is compact (in this
case, $C_0(X)=C(X)$).
\end{example}

\subsection{$K_0$ of a ring}

Suppose $R$ is an arbitrary ring with $1$ (not necessarily
commutative). A module $M$ over $R$ is called finitely generated
projective, if there is another $R$-module $N$ and a number $n\ge 0$
such that
\begin{equation*}
  M\oplus N \iso R^n.
\end{equation*}
This is equivalent to the assertion that the matrix ring
$M_n(R)=End_R(R^n)$ contains an idempotent $e$, i.e.~with $e^2=e$,
such that $M$ is isomorphic to the image of $e$, i.e.~$M\iso e R^n$.

\begin{example} Description of projective modules.
  \begin{enumerate}
  \item If $R$ is a field, the finitely generated projective $R$-modules
    are exactly the finite dimensional vector spaces. (In this case,
    every module is projective).
  \item If $R=\integers$, the finitely generated projective modules
    are the free abelian groups of finite rank
  \item Assume $X$ is a compact topological space and $A=C(X)$. Then,
    by the Swan-Serre theorem \cite{Swan(1962)}, $M$ is a finitely
    generated projective $A$-module if and only if $M$ is isomorphic
    to the space $\Gamma(E)$ of continuous sections of some complex vector
    bundle $E$ over $X$.
\end{enumerate}
\end{example}

\begin{definition}\label{def:K0_of_rings}
  Let $R$ be any ring with unit. $K_0(R)$ is defined to be the
  Grothendieck group of finitely generated projective modules over
  $R$, i.e.~the group of equivalence classes $[(M,N)]$ of pairs of (isomorphism 
  classes of) finitely generated projective $R$-modules $M$, $N$,
  where $(M,N)\equiv (M',N')$ if and only if there is an $n\ge 0$ with 
  \begin{equation*}
    M\oplus N'\oplus R^n\iso M'\oplus N\oplus R^n.
  \end{equation*}
  The group composition is given by 
  \begin{equation*}
    [(M,N)] + [(M',N')] := [ (M\oplus M', N\oplus N')].
  \end{equation*}
  We can think of $(M,N)$ as the formal difference of modules $M-N$.

  Any unital ring homomorphism $f\colon R\to S$ induces a map 
  \begin{equation*}
f_*\colon
  K_0(R)\to K_0(S)\colon [M]\mapsto [S\tensor_R M],
\end{equation*}
where $S$ becomes a right $R$-module via $f$. We obtain that $K_0$ is a
covariant functor from the category of unital rings to the category of 
abelian groups.
\end{definition}

\begin{example} Calculation of $K_0$.
  \begin{itemize}
  \item 
    If $R $ is a field, then $K_0(R)\iso\integers$, the isomorphism
    given by the dimension: $\dim_R(M,N):=\dim_R(M)-\dim_R(N)$.
  \item
    $K_0(\integers)\iso\integers$, given by the rank.
  \item
    If $X$ is a compact topological space, then $K_0(C(X))\iso
    K^0(X)$, the topological K-theory given in terms of complex vector
    bundles. To each vector bundle $E$ one associates the
    $C(X)$-module $\Gamma(E)$ of continuous sections of $E$.
  \item Let $G$ be a discrete group. The group algebra $\complexs G$
    is a vector space with basis $G$, and with multiplication coming
    from the group structure, i.e.~given by
    $g\cdot h = (gh)$.

    If $G$ is a finite group, then $K_0(\complexs G)$ is the complex
    representation ring of $G$.
\end{itemize}
\end{example}

\subsection{K-Theory of $C^*$-algebras}

\begin{definition}
  Let $A$ be a unital $C^*$-algebra. Then $K_0(A)$ is defined as in
  Definition \ref{def:K0_of_rings}, i.e.~by forgetting the topology of 
  $A$.
\end{definition}

\subsubsection{K-theory for non-unital $C^*$-algebras}

When studying (the K-theory of) $C^*$-algebras, one has to understand
morphisms $f\colon A\to B$. This necessarily involves studying the
kernel of $f$, which is a closed ideal of $A$, and hence a
\emph{non-unital} $C^*$-algebra. Therefore, we proceed by defining the 
$K$-theory of $C^*$-algebras without unit.

\begin{definition}
  To any $C^*$-algebra $A$, with or without unit, we assign in a
  functorial way a new, unital $C^*$-algebra $A_+$ as follows. As
  $\complexs$-vector space, $ A_+:= A\oplus \complexs$, with product
  \begin{equation*}
    (a,\lambda)(b,\mu) := (ab+\lambda a+\mu b,
    \lambda\mu)\qquad\text{for }(a,\lambda),(b,\mu)\in A\oplus\complexs.
  \end{equation*}
  The unit is given by $(0,1)$. The star-operation is defined as
  $(a,\lambda)^*:=
  (a^*,\overline{\lambda})$, and the new norm is given by
   \begin{equation*}
     \abs{(a,\lambda)} =\sup \{\abs{a x +\lambda x}\mid x\in A \text{
     with }\abs{x}=1\}
   \end{equation*}
\end{definition}

\begin{remark}
  $A$ is a closed ideal of $A_+$, the kernel of the canonical
  projection $A_+\onto \complexs$ onto the second factor. If $A$
  itself is unital, the unit of $A$ is of course different from the
  unit of $A_+$.
\end{remark}

\begin{example}
  Assume $X$ is a locally compact space, and let $X_+:=
  X\cup\{\infty\}$ be the one-point compactification of $X$. Then 
  \begin{equation*}
    C_0(X)_+ \iso C(X_+).
  \end{equation*}
  The ideal $C_0(X)$ of $C_0(X)_+$ is identified with the ideal of
  those functions
  $f\in C(X_+)$ such that $f(\infty)=0$.
\end{example}

\begin{definition}
  For an arbitrary $C^*$-algebra $A$ (not necessarily unital) define
  \begin{equation*}
    K_0(A):= \ker( K_0(A_+)\to K_0(\complexs)).
  \end{equation*}
  Any $C^*$-algebra homomorphisms $f\colon A\to B$ (not necessarily
  unital) induces a unital homomorphism $f_+\colon A_+\to B_+$. The
  induced map 
  \begin{equation*}
(f_+)_*\colon K_0(A_+)\to K_0(B_+)
\end{equation*}
maps the kernel of
  the map $K_0(A_+)\to K_0(\complexs)$ to the kernel of $K_0(B_+)\to
  K_0(\complexs)$. This means it restricts to a map $f_*\colon
  K_0(A)\to K_0(B)$. We obtain a covariant functor from the category
  of (not necessarily unital) $C^*$-algebras to abelian groups.
\end{definition}

Of course, we need the following result.
\begin{proposition}
  If $A$ is a unital $C^*$-algebra, the new and the old definition of
  $K_0(A)$ are canonically isomorphic.
\end{proposition}

\subsubsection{Higher topological K-groups}

We also want to define higher topological K-theory groups. We have an
ad hoc definition using suspensions (this is similar to the
corresponding idea in topological K-theory of spaces). For this we
need the following.

\begin{definition}
  Let $A$ be a $C^*$-algebra. We define the cone $CA$ and the
  suspension $SA$ as follows.
  \begin{equation*}
    \begin{split}
      CA &:= \{f\colon [0,1]\to A\mid f(0)=0\} \\
      SA & := \{f\colon [0,1]\to A\mid f(0)=0=f(1)\}.
  \end{split}
\end{equation*}
These are again $C^*$-algebras, using pointwise operations and the
supremum norm.

Inductively, we define
\begin{equation*}
  S^0A := A\qquad S^n A:= S(S^{n-1}A)\quad\text{for }n\ge 1.
\end{equation*}
\end{definition}

\begin{definition} Assume $A$ is a $C^*$-algebra.
  For $n\ge 0$, define
  \begin{equation*}
    K_n(A):= K_0(S^n A).
  \end{equation*}
    These are the \emph{topological K-theory groups of $A$}. For each
    $n\ge 0$, we obtain a functor from the category of $C^*$-algebras
    to the category of abelian groups.
\end{definition}

For unital $C^*$-algebras, we can also give a more direct definition
of higher K-groups (in
particular useful for $K_1$, which is then defined in terms of
(classes of) invertible matrices). This is done as follows:

\begin{definition}
  Let $A$ be a unital $C^*$-algebra. Then $Gl_n(A)$ becomes a
  topological group, and we have continuous embeddings
  \begin{equation*}
    Gl_n(A)\into Gl_{n+1}(A)\colon X\mapsto 
    \begin{pmatrix}
      X & 0\\ 0 & 1
    \end{pmatrix}.
  \end{equation*}
  We set $Gl_\infty(A):= \lim_{n\to\infty} Gl_n(A)$, and we equip
  $Gl_\infty(A)$ with the direct limit topology.
\end{definition}

\begin{proposition}
  Let $A$ be a unital $C^*$-algebra. If $k\ge 1$, then
  \begin{equation*}
    K_k(A) = \pi_{k-1}(Gl_\infty(A)) (\iso \pi_k(B Gl_\infty(A))).
  \end{equation*}

  Observe that any unital morphism $f\colon A\to B$ of unital
  $C^*$-algebras induces a map $Gl_n(A)\to Gl_n(B)$ and therefore also 
  between $\pi_k(Gl_\infty(A))$ and $\pi_k(Gl_\infty(B))$. This map
  coincides with the previously defined induced map in topological
  $K$-theory.
\end{proposition}

\begin{remark}
  Note that the topology of the $C^*$-algebra enters the definition of 
  the higher topological K-theory of $A$, and in general the
  topological K-theory of $A$ will be vastly different from the
  algebraic K-theory of the algebra underlying $A$. For connections in 
  special cases, compare \cite{Suslin-Wodzicki(1992)}.
\end{remark}

\begin{example}
  It is well known that $Gl_n(\complexs)$ is connected for each
  $n\in\naturals$. Therefore
  \begin{equation*}
    K_1(\complexs) = \pi_0(Gl_\infty(\complexs)) = 0.
  \end{equation*}
\end{example}

A very important result about $K$-theory of $C^*$-algebras is the
following long exact sequence. A proof can be found e.g.~in
\cite[Proposition 4.5.9]{Higson-Roe(2001)}.
\begin{theorem}\label{theo:long_exact}
  Assume $I$ is a closed ideal of a $C^*$-algebra $A$. Then,
  we get a short exact sequence of $C^*$-algebras $0\to I\to A\to
  A/I\to 0$, which induces a long exact sequence in K-theory
  \begin{equation*}
    \to K_n(I)\to K_n(A)\to K_n(A/I)\to K_{n-1}(I) \to \cdots \to K_0(A/I).
  \end{equation*}
\end{theorem}

\subsection{Bott periodicity and the cyclic exact sequence}

One of the most important and remarkable results about the K-theory of 
$C^*$-algebras is Bott periodicity, which can be stated as follows.

\begin{theorem}
  Assume $A$ is a $C^*$-algebra. There is a natural isomorphism,
  called the Bott map
  \begin{equation*}
    K_0(A)\to K_0(S^2A),
  \end{equation*}
  which implies immediately that there are natural isomorphism
  \begin{equation*}
    K_n(A)\iso K_{n+2}(A)\qquad\forall n\ge 0.
  \end{equation*}
\end{theorem}

\begin{remark}
  Bott periodicity allows us to define $K_n(A)$ for each
  $n\in\integers$, or to regard the K-theory of $C^*$-algebras as a
  $\integers/2$-graded theory, i.e.~to talk of $K_n(A)$ with
  $n\in\integers/2$. This way, the long exact sequence of Theorem
  \ref{theo:long_exact} becomes a (six-term) cyclic exact sequence
  \begin{equation*}
    \begin{CD}
      K_0(I) @>>> K_0(A) @>>> K_0(A/I)\\
      @AAA && @VV{\mu_*}V\\
      K_1(A/I) @<<< K_1(A) @<<< K_1(I).
    \end{CD}
  \end{equation*}
  The connecting homomorphism $\mu_*$ is the composition of the Bott
  periodicity isomorphism and the connecting homomorphism of Theorem
  \ref{theo:long_exact}. 
\end{remark}

\subsection{The $C^*$-algebra of a group}

Let $\Gamma$ be a discrete group. Define $l^2(\Gamma)$ to be the 
Hilbert space of square summable complex valued functions on
$\Gamma$. We can write an element $f\in l^2(\Gamma)$ as a sum
$\sum_{g\in\Gamma}\lambda_g g$ with $\lambda_g\in\complexs$ and
$\sum_{g\in\Gamma}\abs{\lambda_g}^2<\infty$.


We defined the \emph{complex group algebra} (often also called the
\emph{complex
group ring}) $\complexs\Gamma$ to be the complex vector space  with
basis the elements of $\Gamma$ (this can also be considered as the
space of
complex valued functions on $\Gamma$ with finite support, and as such
is a subspace of $l^2(\Gamma)$). The product in $\complexs\Gamma$ is
induced by the multiplication in $\Gamma$, namely, if
$f=\sum_{g\in\Gamma}\lambda_g g, u=\sum_{g\in\Gamma}\mu_g
g\in\complexs\Gamma$, then
\begin{equation*}
  (\sum_{g\in\Gamma}\lambda_g g)(\sum_{g\in\Gamma}\mu_g
g) := \sum_{g,h\in\Gamma} \lambda_g\mu_h (gh) =\sum_{g\in \Gamma}
\left(\sum_{h\in\Gamma}\lambda_h\mu_{h^{-1}g} \right) g .
\end{equation*}
This is a convolution product.

We have the \emph{left regular representation} $\lambda_\Gamma$ of $\Gamma$ on
$l^2(\Gamma)$, given by
\begin{equation*}
  \lambda_\Gamma(g)\cdot(\sum_{h\in \Gamma}\lambda_h h) := \sum_{h\in\Gamma}\lambda_h gh
\end{equation*}
for $g\in\Gamma$ and $\sum_{h\in\Gamma}\lambda_h h\in l^2(\Gamma)$.

This unitary representation extends linearly to $\complexs \Gamma$.

The \emph{reduced $C^*$-algebra} $C^*_r\Gamma$ of $\Gamma$ is defined
to be the norm closure of the image $\lambda_\Gamma(\complexs\Gamma)$
in the $C^*$-algebra of bounded operators on $l^2(\Gamma)$.

\begin{remark}
  It's no surprise that there is also a \emph{maximal $C^*$-algebra}
  $C^*_{max}\Gamma$ of a
  group $\Gamma$. It is defined using not only the left regular
  representation of $\Gamma$, but simultaneously all of its
  representations. We will not make use of $C^*_{max}\Gamma$ in these
  notes, and therefore will not define it here.

  Given a topological group $G$, one can define $C^*$-algebras
  $C^*_rG$ and $C^*_{max}G$ which take the topology of $G$ into
  account. They actually play an important role in the study of the
  Baum-Connes conjecture, which can be defined for (almost arbitrary)
  topological groups, but again we will not cover this subject
  here. Instead, we will throughout stick to discrete groups.
\end{remark}

\begin{example}
  If $\Gamma$ is finite, then $C^*_r\Gamma =\complexs \Gamma$ is the
  complex group ring of $\Gamma$.

  In particular, in this case $K_0(C^*_r\Gamma)\iso R(\Gamma)$
  coincides with the (additive group of) the complex representation
  ring of $\Gamma$.
\end{example}

\section{The Baum-Connes conjecture}

The Baum-Connes conjecture relates an object from algebraic topology,
namely the K-homology of the classifying space of a given group
$\Gamma$, to representation theory and the world of $C^*$-algebras,
namely to the K-theory of the reduced $C^*$-algebra of $\Gamma$.

Unfortunately, the material is very technical. Because of lack of
space and time we can not go into the details (even of some of the
definitions). We recommend the sources \cite{Valette(2001)}, \cite{Valette(2000)},
\cite{Higson-Roe(2001)}, \cite{Baum-Connes-Higson(1994)},
\cite{Mislin(2001)} and
\cite{Blackadar(1998)}.

\subsection{The Baum-Connes conjecture for torsion-free groups}

\begin{definition}
  Let $X$ be any CW-complex. $K_*(X)$ is the K-homology of $X$, where
  K-homology is the homology theory dual to topological K-theory. If
  $BU$ is the spectrum of topological K-theory, and $X_+$ is $X$ with
  a disjoint basepoint added, then
  \begin{equation*}
    K_n(X):= \pi_n(X_+\wedge BU).
  \end{equation*}
\end{definition}

\begin{definition}
  Let $\Gamma$ be a discrete group. A classifying space $B\Gamma$
  for $\Gamma$ is a CW-complex with the property that
  $\pi_1(B\Gamma)\iso\Gamma$, and $\pi_k(B\Gamma)=0$ if $k\ne 1$. A
  classifying space always exists, and is unique up to homotopy
  equivalence. Its universal covering $E\Gamma$ is a contractible
  CW-complex with a free cellular $\Gamma$-action, the so called
  \emph{universal space for $\Gamma$-actions}.
\end{definition}

\begin{remark}
  In the literature about the Baum-Connes conjecture, one will often
  find the definition 
  \begin{equation*}
    RK_n(X):= \dirlim K_n(Y),
  \end{equation*}
  where the limit is taken over all finite subcomplexes $Y$ of
  $X$. Note, however, that K-homology (like any homology theory in
  algebraic topology) is compatible with direct limits, which implies
  $RK_n(X)=K_n(X)$ as defined above. The confusion comes from the fact 
  that operator algebraists often use Kasparov's bivariant KK-theory to 
  define $K_*(X)$, and this coincides with the homotopy theoretic 
  definition only if $X$ is compact.
\end{remark}

Recall that a group $\Gamma$ is called torsion-free, if $g^n=1$ for
$g\in \Gamma$ and $n>0$ implies that $g=1$.

We can now formulate the Baum-Connes conjecture for torsion-free
discrete groups.
\begin{conjecture}\label{conj:torsion_free_BC}
  Assume $\Gamma$ is a torsion-free discrete group. It is known that
  there is a particular
  homomorphism, the assembly map
  \begin{equation}\label{eq:BC}
    \overline\mu_*\colon K_*(B\Gamma) \to K_*(C^*_r \Gamma)
  \end{equation}
  (which will be defined later). The \emph{Baum-Connes conjecture} says that
  this map is an isomorphism.
\end{conjecture}

\begin{example}\label{ex:finite_group_assembly_not_iso}
  The map $\overline{\mu}_*$ of Equation \eqref{eq:BC} is also defined
  if $\Gamma$ is
  not torsion-free. However, in this situation it will in general not
  be an isomorphism. This can already be seen if $\Gamma=\integers/2$. Then 
  $C^*_r \Gamma=\complexs \Gamma\iso \complexs\oplus\complexs$ as a
  $\complexs$-algebra. Consequently, 
  \begin{equation}\label{eq:K_of_Z2}
K_0(C^*_r\Gamma)\iso
  K_0(\complexs)\oplus K_0(\complexs) \iso \integers\oplus\integers.
\end{equation}
On the other hand, 
using the homological Chern character, 
\begin{equation}\label{eq:H_of_Z2}
 K_0(B\Gamma)\tensor_\integers\rationals \iso \oplus_{n=0}^\infty
 H_{2n}(B\Gamma;\rationals) \iso\rationals.
\end{equation}
(Here we use the fact that the rational homology of every finite group 
is zero in positive degrees, which follows from the fact that the
transfer homomorphism $H_k(B\Gamma;\rationals)\to H_k(\{1\};\rationals)$ is
(with rational coefficients) up to a factor $\abs{\Gamma}$ a left inverse
to the map induced from
the inclusion, and therefore is injective.)

The calculations \eqref{eq:K_of_Z2} and \eqref{eq:H_of_Z2} prevent
$\mu_0$ of \eqref{eq:BC} from being an isomorphism.
\end{example}

\subsection{The Baum-Connes conjecture in general}

To account for the problem visible in Example
\ref{ex:finite_group_assembly_not_iso} if we are dealing with groups
with torsion, one replaces the left hand side by a more complicated
gadget, the equivariant K-homology of a certain $\Gamma$-space
$E(\Gamma,fin)$, the classifying space for proper actions. We will
define all of this later. Then, the Baum-Connes conjecture says the
following.
\begin{conjecture}\label{conj:general_BC}
    Assume $\Gamma$ is a discrete group. It is known that there is a particular
  homomorphism, the assembly map
  \begin{equation}\label{eq:BC_general}
    \mu_*\colon K^{\Gamma}_*(E(\Gamma,fin)) \to K_*(C^*_r \Gamma)
  \end{equation}
  (we will define it later). The conjecture says that this map is an
  isomorphism.
\end{conjecture}

\begin{remark}
  If $\Gamma$ is torsion-free, then
  $K_*(B\Gamma)=K_*^\Gamma(E(\Gamma,fin))$, and the assembly maps
  $\overline\mu$ of
  Conjectures \ref{conj:torsion_free_BC} and $\mu$ of \ref{conj:general_BC}
  coincide (see Proposition \ref{prop:Davis_Lueck_properties}).
\end{remark}

Last, we want to mention that there is also a \emph{real version} of the 
Baum-Connes conjecture, where on the left hand side the K-homology is
replaced by KO-homology, i.e.~the homology dual to the K-theory of
real vector spaces (or an equivariant version hereof), and on the right
hand side $C^*_r\Gamma$ is replaced by the real reduced
$C^*$-algebra $C^*_{r,\reals}\Gamma$.

\subsection{Consequences of the Baum-Connes conjecture}

\subsubsection{Idempotents in $C^*_r\Gamma$}

The connection between the Baum-Connes conjecture and idempotents is
best shown via Atiyah's $L^2$-index theorem, which we discuss first.

Given a closed manifold $M$ with an elliptic differential operator
$D\colon C^\infty(E)\to C^\infty(F)$ between two bundles on $M$, and a 
normal covering $\tilde M\to M$ (with deck transformation group
$\Gamma$, normal means that $M=\tilde M/\Gamma$), we can lift $E$, $F$ 
and $D$ to $\tilde M$, and get an elliptic $\Gamma$-equivariant
differential operator $\tilde D\colon C^\infty(\tilde E)\to
C^\infty(\tilde F)$. If $\Gamma$ is not finite, we can not use the
equivariant index of Section \ref{sec:equiv-index}. However, because
the action is free, it is possible to define an equivariant analytic index
\begin{equation*}
  \ind_\Gamma(\tilde D)\in K_{\dim M}(C^*_r\Gamma).
\end{equation*}
This is described in Example \ref{ex:Atiyah-operator_on_covering}. 

Atiyah used a certain real valued homomorphism, the $\Gamma$-dimension 
\begin{equation*}
\dim_\Gamma\colon K_0(C^*_r\Gamma)\to\reals,
\end{equation*}
to define the
$L^2$-index of $\tilde D$ (on an even dimensional manifold):
\begin{equation*}
  L^2\text{-}\ind(\tilde D):= \dim_\Gamma(\ind_\Gamma(\tilde D)).
\end{equation*}
The $L^2$-index theorem says
\begin{equation*}
  L^2\text{-}\ind(\tilde D) = \ind(D),
\end{equation*}
in particular, it follows that the $L^2$-index is an integer. For a
different point of view of the $L^2$-index theorem, compare Section
\ref{sec:analytic-l2-betti}.

An alternative description of the left hand side of \eqref{eq:BC} and
\eqref{eq:BC_general} shows that, as long as $\Gamma$ is torsion-free, 
the image of $\mu_0$ coincides with the subset of $K_0(C^*_r\Gamma)$
consisting of $\ind_\Gamma(\tilde D)$, where $\tilde D$ is as above. In
particular, if $\mu_0$ is surjective (and $\Gamma$ is torsion-free),
for each $x\in K_0(C^*_r\Gamma)$ we find a differential operator $D$
such that $x=\ind_\Gamma(\tilde D)$. As a consequence,
$\dim_\Gamma(x)\in\integers$, i.e.~the range of $\dim_\Gamma$ is
contained in $\integers$. This is the statement of the so called
\emph{trace conjecture}.
\begin{conjecture}\label{conj:trace_conjecture}
  Assume $\Gamma$ is a torsion-free discrete group. Then
  \begin{equation*}
    \dim_\Gamma(K_0(C^*_r\Gamma))\subset\integers.
  \end{equation*}
\end{conjecture}

On the other hand, if $x\in
K_0(C^*_r\Gamma)$ is represented by a projection $p=p^2\in
C^*_r\Gamma$, then elementary properties of $\dim_\Gamma$ (monotonicity 
and faithfulness) imply that $0\le \dim_\Gamma(p)\le 1$, and
$\dim_\Gamma(p)\notin \{0,1\}$ if $p\ne 0,1$.

Therefore, we have the following consequence of the Baum-Connes
conjecture. If $\Gamma$ is torsion-free and the Baum-Connes map $\mu_0$ is 
surjective, then $C^*_r\Gamma$ does not contain any projection
different from $0$ or $1$.

This is the assertion of the Kadison-Kaplansky conjecture:
\begin{conjecture}
  Assume $\Gamma$ is torsion-free. Then $C^*_r\Gamma$ does not contain 
  any non-trivial projections.
\end{conjecture}

The following consequence of the Kadison-Kaplansky conjecture deserves 
to be mentioned:
\begin{proposition}
  If the Kadison-Kaplansky conjecture is true for a group $\Gamma$,
then the spectrum $s(x)$ of every self adjoint element $x\in C^*_r\Gamma$ is
connected. Recall that the spectrum is defined in the following way:
\begin{equation*}
s(x):=\{\lambda\in\complexs\mid (x-\lambda\cdot 1) \text{ not
  invertible}\}. 
\end{equation*}
\end{proposition}

If $\Gamma$ is not torsion-free, it is easy to construct non-trivial
projections, and it is clear that the range of $\ind_\Gamma$ is
not contained in $\integers$. Baum and Connes originally conjectured
that it is contained in the abelian subgroup $Fin^{-1}(\Gamma)$ of
$\rationals$ generated by $\{1/\abs{F}\mid F \text{ finite subgroup of 
  }\Gamma\}$. This conjecture is not correct, as is shown by an
example of Roy \cite{Roy(1999)}. In \cite{Lueck(2001a)}, L\"uck proves 
that the Baum-Connes conjecture implies that the range of
$\dim_\Gamma$ is contained in the \emph{subring} of $\rationals$
generated by $\{1/\abs{F}\mid F \text{ finite subgroup of 
  }\Gamma\}$.

\subsubsection{Obstructions to positive scalar curvature}

The Baum-Connes conjecture implies the so called ``stable
Gromov-Lawson-Rosenberg'' conjecture. This implication is a theorem
due to Stephan
Stolz. The details of this will be
discussed in the lectures of Stephan Stolz, therefore we can be very
brief. We just mention the result.

\begin{theorem}\label{theo:stable_GLR}
  Fix a group $\Gamma$.
 Assume that $\mu$ in the real version of \eqref{eq:BC_general}
 discussed in Section \ref{sec:real-c-algebras} is
  injective (which follows e.g.~if $\mu$ in \eqref{eq:BC_general} is
  an isomorphism), and assume that $M$ is a closed spin manifold with
  $\pi_1(M)=\Gamma$. Assume that a certain (index theoretic)
  invariant $\alpha(M)\in K_{\dim M}(C^*_{\reals,r}\Gamma)$
  vanishes. Then there is an $n\ge 0$ such that $M\times B^n$ admits a 
  metric with positive scalar curvature.
\end{theorem}

Here, $B$ is any simply connected $8$-dimensional spin manifold with
$\hat{A}(M)=1$. Such a manifold is called a \emph{Bott manifold}.

The converse of Theorem \ref{theo:stable_GLR}, i.e.~positive scalar
curvature implies vanishing of $\alpha(M)$, is true for arbitrary
groups and without knowing anything about the Baum-Connes conjecture.

\subsubsection{The Novikov conjecture about higher signatures}

\paragraph{Direct approach}

The original form of the Novikov conjecture states that higher
signatures are homotopy invariant.

More precisely, let $M$ be an (even dimensional) closed oriented manifold
with fundamental group $\Gamma$. Let $B\Gamma$ be a classifying space
for $\Gamma$. There is a unique (up to homotopy) \emph{classifying
  map} $u\colon M\to B\Gamma$ which is defined by the property that it 
induces an isomorphism on $\pi_1$. Equivalently, $u$ classifies a
universal covering of $M$.

Let $L(M)\in H^*(M;\rationals)$ be the Hirzebruch L-class (as
normalized by Atiyah and Singer). Given any cohomology class $a\in
H^*(B\Gamma,\rationals)$, we define the higher signature
\begin{equation*}
  \sigma_a(M):= \innerprod{L(M)\cup u^*a, [M]} \in\rationals.
\end{equation*}
Here $[M]\in H_{\dim M}(M;\rationals)$ is the fundamental class of the 
oriented manifold $M$, and $\innerprod{\cdot,\cdot}$ is the usual
pairing between cohomology and homology.

Recall that the Hirzebruch signature theorem states that $\sigma_1(M)$
is the signature of $M$, which evidently is an oriented homotopy invariant.

The Novikov conjecture generalizes this as follows.
\begin{conjecture}
  Assume $f\colon M\to M'$ is an oriented homotopy equivalence between 
  two even dimensional closed oriented manifolds, with (common)
  fundamental group
  $\pi$. ``Oriented'' means that $f_*[M]=[M']$. Then all higher
  signatures of $M$ and
  $M'$ are equal, i.e.
  \begin{equation*}
    \sigma_a(M) = \sigma_a(M')\qquad\forall a\in H^*(B\Gamma,\rationals).
  \end{equation*}
\end{conjecture}

There is an equivalent reformulation of this conjecture in terms of K-ho\-mo\-lo\-gy.
To see this, let $D$ be the signature operator of $M$. (We assume here
that $M$ is smooth, and we choose a Riemannian metric on $M$ to define 
this operator. It is an elliptic differential operator on $M$.)
The operator $D$ defines an element in the K-homology of $M$, $[D]\in
K_{\dim M}(M)$. Using the map $u$, we can push $[D]$ to $K_{\dim
  M}(B\Gamma)$. We define the higher signature $\sigma(M):= u_*[D] \in 
K_{\dim M}(B\Gamma)\tensor\rationals$. It turns out that 
\begin{equation*}
  2^{\dim M/2} \sigma_a(M) =\innerprod{a, ch (\sigma(M))}\qquad\forall 
  a\in H^*(B\Gamma;\rationals),
\end{equation*}
where $ch\colon K_*(B\Gamma)\tensor\rationals \to
H_*(B\Gamma,\rationals)$ is the homological Chern character (an
isomorphism).

Therefore, the Novikov conjecture translates to the statement that
$\sigma(M)=\sigma(M')$ if $M$ and $M'$ are oriented homotopy equivalent.

Now one can show \emph{directly} that
\begin{equation*}
  \overline\mu(\sigma(M)) = \overline\mu(\sigma(M')) \in K_*(C^*_r\Gamma),
\end{equation*}
if $M$ and $M'$ are oriented homotopy equivalent. Consequently,
rational injectivity of the Baum-Connes map $\overline\mu$ immediately implies the
Novikov conjecture. If $\Gamma$ is torsion-free, this is part of the
assertion of the Baum-Connes conjecture. Because of this relation,
injectivity of the Baum-Connes
map $\mu$ is often called the ``analytic Novikov conjecture''.

\paragraph{L-theory approach}

There is a more obvious connection between the Baum-Connes isomorphism 
conjecture and the L-theory isomorphism conjecture (discussed in other 
lectures).

Namely, the L-theory isomorphism conjecture is concerned with a
certain assembly map
\begin{equation*}
  A_\Gamma \colon H_*(B\Gamma, \mathbb{L}_{\bullet}(\integers)) \to
  L_*(\integers[\Gamma]). 
\end{equation*}
Here, the left hand side is the homology of $B\Gamma$ with
coefficients  the algebraic surgery spectrum of $\integers$, and the
right hand side is the free quadratic $L$-group of the ring with
involution $\integers[\Gamma]$.

The Novikov conjecture is equivalent to the statement that
this map is rationally injective, i.e.~that
\begin{equation*}
  A_\Gamma\tensor\id_\rationals \colon H_*(B\Gamma, \mathbb{L}_{\bullet}(\integers))\tensor\rationals \to
  L_*(\integers[\Gamma])\tensor\rationals
\end{equation*}
is an injection. This formulation has the advantage that, tensored
with $\rationals$, all the different flavors of L-theory are
isomorphic (therefore, we don't have to and we won't discuss these
distinctions here).

Now, we get a commutative diagram
\begin{equation}\label{eq:Novikov_diagram}
  \begin{CD}
    H_*(B\Gamma, \mathbb{L}_{\bullet}(\integers))\tensor\rationals @>{A_\Gamma}>>
    L_*(\integers[\Gamma])\tensor\rationals\\
    @VVV  @VVV\\
    H_*(B\Gamma, \mathbb{L}_{\bullet}(\complexs))\tensor\rationals
    @>{A_{\Gamma,\complexs}}>>
    L_*(\complexs[\Gamma])\tensor\rationals\\
    @VVV @VVV\\
    K_*(B\Gamma)\tensor\rationals @>{\overline\mu}>> 
    K_*(C_r^*\Gamma)\tensor\rationals =L_*(C^*_r\Gamma)\tensor\rationals.
 \end{CD}
\end{equation}
The maps on the left hand side are given by natural transformations of 
homology theories with values in rational vector spaces. These
transformations are easily seen to be injective for the coefficients.
Since we deal with \emph{rational} homology theories, they are
injective in general.

The maps on the right hand side are the maps in L-theory induced by
the obvious ring homomorphisms $\integers\Gamma\to \complexs\Gamma\to
C^*_r\Gamma$. Then we use the ``folk theorem'' that, for
$C^*$-algebras, K-theory and L-theory are canonically isomorphic (even 
non-rationally). Of course, it remains to establish commutativity of
the diagram \eqref{eq:Novikov_diagram}. For
more details, we
refer to \cite{Rosenberg(1995)}. Using all these facts and the diagram
\eqref{eq:Novikov_diagram}, we see that
for torsion-free groups, rational injectivity of the
Baum-Connes map $\mu$ implies rational injectivity of the L-theory
assembly $A_\Gamma$, i.e.~the Novikov conjecture.

\paragraph{Groups with torsion}

For an arbitrary group $\Gamma$, we have a factorization of $\overline\mu$ as
follows:
\begin{equation*}
  K_*(B\Gamma)\xrightarrow{f} K_*^\Gamma(E(\Gamma,fin)) \xrightarrow{\mu}
  K_*(C^*_r\Gamma).
\end{equation*}
One can show that $f$ is rationally injective, so that rational
injectivity of the Baum-Connes map $\mu$ implies the Novikov
conjecture also in general.

\subsection{The universal space for proper actions}

\begin{definition}
  Let $\Gamma$ be a discrete group and $X$ a Hausdorff space with an
  action of $\Gamma$. We say that the action is \emph{proper}, if
  for all $x,y\in X$ there are open neighborhood $U_x\ni x$ and
  $U_y\ni y$ such that $gU_x\cap U_y$ is non-empty only for finitely
  many $g\in \Gamma$ (the number depending on $x$ and $y$).

  The action is said to be \emph{cocompact}, if $X/\Gamma$ is compact.
\end{definition}

\begin{lemma}
  If the action of $\Gamma$ on $X$ is proper, then for each $x\in X$
  the \emph{isotropy group} $\Gamma_x:=\{g\in\Gamma\mid gx=x\}$ is
  finite.
\end{lemma}

\begin{definition}
    Let $\Gamma$ be a discrete group. A CW-complex $X$ is a
    \emph{$\Gamma$-CW-complex}, if $X$ is a CW-complex with a
  cellular action of $\Gamma$ with the additional property that,
  whenever $g(D)\subset D$ for a cell $D$ of $X$ and some
  $g\in\Gamma$, then $g|_D=\id_D$, i.e.~$g$ doesn't move $D$ at all.
\end{definition}

\begin{remark}
  There exists also the notion of $G$-CW-complex for topological
  groups $G$ (taking the topology of $G$ into account). These have to
  be defined in a different way, namely by gluing together
  $G$-equivariant cells $D^n\times G/H$. In general, such a
  $G$-CW-complex is not an ordinary CW-complex.
\end{remark}

\begin{lemma}
  The action of a discrete group $\Gamma$ on a $\Gamma$-CW-complex is
  proper if and only if every isotropy group is finite.
\end{lemma}

\begin{definition}\label{def:univesal_space_for_proper_actions}
  A proper $\Gamma$-CW-complex $X$ is called \emph{universal}, or more
  precisely \emph{universal for proper actions}, if for every proper
  $\Gamma$-CW-complex $Y$ there is a $\Gamma$-equivariant map $f\colon 
  Y\to X$ which is unique up to $\Gamma$-equivariant homotopy. Any
  such space is denoted $E(\Gamma,fin)$ or $\underline{E}\Gamma$.
\end{definition}

\begin{proposition}
  A $\Gamma$-CW-complex $X$ is universal for proper actions if and
  only if the fixed point set
  \begin{equation*}
X^H:=\{x\in X\mid hx=x\quad\forall h\in H\}
\end{equation*}
is empty whenever $H$ is an infinite subgroup of $\Gamma$, and is
contractible (and in particular non-empty) if $H$ is a finite subgroup 
of $\Gamma$.
\end{proposition}

\begin{proposition}
  If $\Gamma$ is a discrete group, then $E(\Gamma,fin)$ exists and is
  unique up to $\Gamma$-homotopy equivalence.
\end{proposition}

\begin{remark}
  The general context for this discussion are actions of a group
  $\Gamma$ where the isotropy belongs to a fixed family of subgroups
  of $\Gamma$ (in our case, the family of all finite subgroups). For
  more information, compare \cite{Dieck(1972)}.
\end{remark}

\begin{example}\strut
  \begin{itemize}
  \item If $\Gamma$ is torsion-free, then $E(\Gamma,fin)=E\Gamma$, the 
    universal covering of the classifying space $B\Gamma$. Indeed,
    $\Gamma$ acts freely on $E\Gamma$, and $E\Gamma$ is contractible.
  \item If $\Gamma$ is finite, then $E(\Gamma,fin)=\{*\}$.
  \item If $G$ is a connected Lie group with maximal compact subgroup
    $K$, and $\Gamma$ is a discrete subgroup of $G$, then
    $E(\Gamma,fin)= G/K$ \cite[Section 2]{Baum-Connes-Higson(1994)}.
  \end{itemize}
\end{example}

\begin{remark}
  In the literature (in particular, in
  \cite{Baum-Connes-Higson(1994)}), also a slightly different notion of
  universal
  spaces is discussed. One allows $X$ to be any proper metrizable
  $\Gamma$-space, and requires the universal property for all proper
  metrizable $\Gamma$-spaces $Y$. For discrete groups (which are the
  only groups we are discussing here), a
  universal space in the sense of Definition
  \ref{def:univesal_space_for_proper_actions} is universal in this
  sense.

  However, for some of the proofs
  of the Baum-Connes
  conjecture (for special groups) it is useful to use certain models
  of $E(\Gamma,fin)$ (in the broader sense) coming from the geometry
  of the group, which are not
  $\Gamma$-CW-complexes. 
\end{remark}

\subsection{Equivariant K-homology}

Let $\Gamma$ be a discrete group. We have seen that, if $\Gamma$ is
not torsion-free, the
assembly map \eqref{eq:BC} is not an isomorphism. To account for that, 
we replace $K_*(B\Gamma)$ by the equivariant K-theory of
$E(\Gamma,fin)$. Let $X$ be any proper 
$\Gamma$-CW complex. The original definition of equivariant K-homology is
due to Kasparov, making ideas of Atiyah precise. In this definition,
elements of $K_*^\Gamma(X)$ are equivalence classes of generalized
elliptic operators. In \cite{Davis-Lueck(1998)}, a more homotopy
theoretic definition of $K_*^\Gamma(X)$ is given, which puts the
Baum-Connes conjecture in the context of other isomorphism conjectures.

\subsubsection[Homotopy theoretic equivariant K-homology]{Homotopy
  theoretic definition of equivariant K-homology}
\label{sec:homot-ther-defin}

The details of this definition are quite technical, using spaces and
spectra over the orbit category of the discrete group $\Gamma$. The
objects of the orbit category are the orbits $\Gamma/H$, $H$ any
subgroup of $\Gamma$. The morphisms from $\Gamma/H$ to $\Gamma/K$ are
simply the $\Gamma$-equivariant maps.
In this setting, any spectrum over the orbit category gives rise to an
equivariant homology theory. The decisive step is then the
construction of a (periodic) topological K-theory spectrum
$\mathbf{K}^\Gamma$ over the orbit category of $\Gamma$. This gives
us then a functor from the category of (arbitrary)
$\Gamma$-CW-complexes to the category of (graded) abelian groups, the
\emph{equivariant $K$-homology} $K^\Gamma_*(X)$ ($X$ any
$\Gamma$-CW-complex).

The important property (which justifies the name ``topological
K-theory spectrum) is that
\begin{equation*}
  K^\Gamma_k(\Gamma/H) = \pi_k(\mathbf{K}^\Gamma(\Gamma/H)) \iso K_k(C^*_rH)
\end{equation*}
for every subgroup $H$ of $\Gamma$. In particular,
\begin{equation*}
  K^\Gamma_k(\{*\}) \iso K_k(C^*_r\Gamma).
\end{equation*}

Moreover, we have the following properties:
\begin{proposition}\label{prop:Davis_Lueck_properties}
\begin{enumerate}
\item Assume $\Gamma$ is the trivial group. Then
  \begin{equation*}
    K^\Gamma_*(X) = K_*(X),
  \end{equation*}
  i.e.~we get back the ordinary K-homology introduced above.
\item If $H\subgroup \Gamma$ and $X$ is an $H$-CW-complex, then there
  is a natural isomorphism
      \begin{equation*}
        K_*^H(X) \iso K_*^\Gamma(\Gamma\times_H X).
      \end{equation*}
      Here $\Gamma\times_H X = \Gamma\times H/\sim$, where we divide
      out the equivalence relation generated by $(gh,x)\sim (g,hx)$
      for $g\in\Gamma$, $h\in H$ and $x\in X$. This is in the obvious
      way a left $\Gamma$-space.
    \item Assume $X$ is a free $\Gamma$-CW-complex. Then there is a
      natural isomorphism
  \begin{equation*}
    K_*(\Gamma\backslash X)\to K_*^\Gamma(X).
  \end{equation*}
  In particular, using the canonical $\Gamma$-equivariant map
  $E\Gamma\to E(\Gamma,fin)$, we get a natural homomorphism
  \begin{equation*}
    K_*(B\Gamma) \xrightarrow{\iso} K_*^\Gamma(E\Gamma)\to K_*^\Gamma(E(\Gamma,fin)).
  \end{equation*}
\end{enumerate}
\end{proposition}

\subsubsection[Analytic equivariant K-homology]{Analytic definition
  of equivariant K-homology}

Here we will give the original definition, which embeds into the
powerful framework of equivariant KK-theory, and which is used for
almost all proofs of special cases of the Baum-Connes
conjecture. However, to derive some of the consequences of the
Baum-Connes conjecture, most notably about the positive scalar
curvature question ---this is discussed in one of the lectures of 
Stephan Stolz--- the homotopy theoretic definition is used.

\begin{definition}
  A Hilbert space $H$ is called \emph{($\integers/2$)-graded}, if $H$
  comes with an orthogonal sum
  decomposition $H=H_0\oplus H_1$. Equivalently, a
  unitary operator $\epsilon$ with $\epsilon^2=1$ is given on
  $H$. The subspaces $H_0$ and $H_1$ can be recovered as the 
  $+1$ and $-1$ eigenspaces of $\epsilon$, respectively.

  A bounded operator $T\colon H\to H$ is called \emph{even} (with respect to the
  given grading), if $T$ commutes with $\epsilon$, and \emph{odd}, if
  $\epsilon$ and
  $T$ anti-commute, i.e.~if $T\epsilon=-\epsilon T$. An even operator 
  decomposes as $T=\left(
    \begin{smallmatrix}
      T_0 & 0 \\ 0 & T_1
    \end{smallmatrix}\right)$, an odd one as $T= \left(
    \begin{smallmatrix}
      0 & T_0\\ T_1 & 0
    \end{smallmatrix}\right)$ in the given decomposition $H=H_0\oplus H_1$.
\end{definition}

\begin{definition}
  A \emph{generalized elliptic $\Gamma$-operator on $X$}, or a \emph{cycle
    for $\Gamma$-K-homology of the $\Gamma$-space $X$}, simply a \emph{cycle} for short, is a triple
  $(H,\pi,F)$, where 
  \begin{itemize}
  \item $H=H_0\oplus H_1$ is a $\integers/2$-graded $\Gamma$-Hilbert space
    (i.e.~the direct sum of two Hilbert spaces with unitary
    $\Gamma$-action)
  \item $\pi$ is a $\Gamma$-equivariant $*$-representation of
    $C_0(X)$ on even bounded operators of $H$ (equivariant means that
    $\pi(fg^{-1}) = g\pi(f)g^{-1}$ for all $f\in C_0(X)$ and all
    $g\in\Gamma$.
  \item $F\colon H\to H$ is a bounded, $\Gamma$-equivariant, self adjoint
    operator such that $\pi(f) (F^2-1)$ and
    $[\pi(f),F]:=\pi(f)F-F\pi(f)$ are compact operators for all $f\in
    C_0(X)$.
    Moreover, we require that $F$ is odd, i.e.~$F=
\left(    \begin{smallmatrix}
      0 & D^*\\ D & 0
    \end{smallmatrix}\right)$
    in the decomposition $H=H_0\oplus H_1$.
\end{itemize}
\end{definition}

\begin{remark}
  There are many different definitions of cycles, slightly
  weakening or strengthening some of the conditions. Of course, this
  does not effect the equivariant K-homology groups which are
  eventually defined using them.
\end{remark}

\begin{definition}
  We define the direct sum of two cycles in the obvious way.
\end{definition}

\begin{definition}
  Assume $\alpha=(H,\pi,F)$ and $\alpha'=(H',\pi',F')$ are two cycles.
  \begin{enumerate}
  \item They are called (isometrically) isomorphic, if there is a
    $\Gamma$-equivariant grading preserving isometry $\Psi\colon H\to
    H'$ such that $\Psi\circ \pi(f)= \pi'(f)\circ \Psi$ for all $f\in
    C_0(X)$ and $\Psi\circ F=F'\circ \Psi$.
  \item They are called \emph{homotopic} (or \emph{operator
      homotopic}) if $H=H'$, $\pi=\pi'$, and there is a norm
    continuous path $(F_t)_{t\in[0,1]}$ of operators with $F_0=F$ and
    $F_1=F'$ and such that $(H,\pi,F_t)$ is a cycle for each $t\in
    [0,1]$.
  \item $(H,\pi,F)$ is called \emph{degenerate}, if $[\pi(f),F]=0$ and 
    $\pi(f)(F^2-1)=0$ for each $f\in C_0(X)$.
  \item The two cycles are called equivalent if there are degenerate
    cycles $\beta$ and $\beta'$ such that $\alpha\oplus\beta$ is
    operator homotopic to a cycle isometrically isomorphic to
    $\alpha'\oplus\beta'$.
  \end{enumerate}

  The set of equivalence classes of cycles is denoted
  $KK_0^\Gamma(X)$. (Caution, this is slightly unusual, mostly one
  will find the notation $K^\Gamma(X)$ instead of $KK^\Gamma(X)$).
\end{definition}

\begin{proposition}
  Direct sum induces the structure of an abelian group on $KK_0^\Gamma(X)$.
\end{proposition}

\begin{proposition}
  Any proper $\Gamma$-equivariant map $\phi\colon X\to Y$ between two
  proper $\Gamma$-CW-complexes induces a homomorphism
  \begin{equation*}
    KK^\Gamma_0(X)\to KK^\Gamma_0(Y)
  \end{equation*}
  by $(H,\pi,F)\mapsto (H,\pi\circ \phi^*,F)$, where $\phi^*\colon
  C_0(Y)\to C_0(X)\colon f\mapsto f\circ \phi$ is defined since $\phi$
  is a proper map (else
  $f\circ \phi$ does not necessarily vanish at infinity).
\end{proposition}
Recall that a continuous map $\phi\colon X\to Y$ is called
\emph{proper} if the inverse image of every compact subset of $Y$ is
compact .

  It turns out that the analytic definition of equivariant K-homology
  is quite flexible. It is designed to make it easy to construct
  elements of these groups ---in many geometric situations they
  automatically show up. We give one of the most typical examples of
  such a situation, which we will need later.

\begin{example}\label{ex:Atiyah-operator_on_covering}
  Assume that
  $M$ is a compact even dimensional Riemannian manifold. Let
  $X=\overline M$ be a normal covering of $M$ with deck transformation group
   $\Gamma$ (normal means that $X/\Gamma=M$). Of course, the action is
   free, in particular, proper. Let $E=E_0\oplus E_1$ be a graded
  Hermitian vector bundle on $M$, and 
  \begin{equation*}
D\colon C^\infty(E)\to C^\infty(E)
\end{equation*}
an
  odd elliptic self adjoint differential operator (odd means that $D$
  maps the subspace
  $C^\infty(E_0)$ to $C^\infty(E_1)$, and vice versa). If $M$ is
  oriented, the signature operator on $M$ is such an operator, if $M$
  is a spin-manifold, the same is true for its Dirac operator.

  Now we can pull back $E$ to a bundle $\overline E$ on $\overline M$, and
  lift $D$ to an operator $\overline D$ on $\overline E$. The assumptions
  imply that $\overline D$ extends to an unbounded self adjoint operator
  on $L^2(\overline E)$, the space of square integrable sections of
  $\overline E$. This space is the completion of $C^\infty_c(\overline
  E)$ with
  respect to the canonical inner product (compare Definition
  \ref{def:L2_sections}). (The
  subscript c denotes sections with compact support). Using the
  functional calculus, we can replace
  $\overline D$ by 
  \begin{equation*}
F:= (\overline D^2+1)^{-1/2} \overline D\colon L^2(\overline
  E)\to L^2(\overline E).
\end{equation*}
Observe that
\begin{equation*}
 L^2(\overline E)=L^2(\overline
  E_0)\oplus L^2(\overline E_1)
\end{equation*}
is a $\integers/2$-graded Hilbert
  space with a unitary
  $\Gamma$-action, which admits an (equivariant) action $\pi$ of
  $C_0(\overline M)=C_0(X)$ by fiber-wise multiplication. This action
  preserves
  the grading. Moreover, $\overline D$ as well as $F$ are odd,
  $\Gamma$-equivariant, self adjoint operators on $L^2(\overline E)$ and $F$
  is a bounded operator. From ellipticity it follows that
  \begin{equation*}
  \pi(f)(F^2-1) = -\pi(f) 
  (\overline D^2+1)^{-1} 
\end{equation*}
is compact for each $f\in C_0(\overline M)$
  (observe that this is not true for $(\overline D^2+1)^{-1}$ itself, if
  $\overline M$ is not compact). Consequently, $(L^2(E),\pi,F)$ defines
  an (even) cycle for $\Gamma$-K-homology, i.e.~it represents an element 
  in $KK^\Gamma_0(X)$.

  One can slightly reformulate the construction as follows: $\overline M$ 
  is a principal $\Gamma$-bundle over $M$, and $l^2(\Gamma)$ has a (unitary)
  left $\Gamma$-action. We therefore can construct the associated flat 
  bundle
  \begin{equation*}
 L:=l^2(\Gamma)\times_\Gamma\overline M
\end{equation*}
on $M$ with fiber
  $l^2(\Gamma)$. Now we can twist $D$ with this bundle $L$,
  i.e.~define 
  \begin{equation*}
\overline D:=\nabla_L\tensor \id+\id\tensor D\colon
  C^\infty(L\tensor E)\to C^\infty(L\tensor E),
\end{equation*}
using the given flat
  connection $\nabla_L$ on $L$. Again, we can
  complete to $L^2(L\tensor E)$ and define 
  \begin{equation*}
F:=(\overline
  D^2+1)^{-1/2}\overline D.
\end{equation*}
The left action of $\Gamma$ on $l^2\Gamma$
  induces an action of $\Gamma$ on $L$ and then a unitary action on
  $L^2(L\tensor E)$. Since $\nabla_L$ preserves the $\Gamma$-action,
  $\overline D$ is $\Gamma$-equivariant.
  There is
  a canonical $\Gamma$-isometry
  between $L^2(L\tensor E)$ and $L^2(\overline E)$ which identifies the
  two versions of $\overline D$ and $F$. The action of $C_0(\overline M)$ on 
  $L^2(L\tensor E)$ can be described by identifying $C_0(\overline M)$ with the
  continuous sections of $M$ on the associated bundle
  \begin{equation*}
C_0(\Gamma)\times_\Gamma\overline M,
\end{equation*}
where $C_0(\Gamma)$ is the
  $C^*$-algebra of functions on $\Gamma$ vanishing at infinity, and
  then using the obvious action of $C_0(\Gamma)$ on $l^2(\Gamma)$.

  It is easy to see how this examples generalizes to
  $\Gamma$-equivariant elliptic differential operators on manifolds
  with a proper, but not necessarily free, $\Gamma$-action (with the
  exception of the last part, of course).

  Work in progress of Baum and Schick \cite{Baum-Schick(2001)}
  suggests the (somewhat surprising) fact that, given any proper
  $\Gamma$-CW-complex $Y$, we can,
  for each element $y\in KK^\Gamma_0(Y)$, find such a proper
  $\Gamma$-manifold $X$, together with a $\Gamma$-equivariant map
  $f\colon X\to Y$ and an elliptic differential operator on $X$ giving
   an element $x\in KK^\Gamma_0(X)$ as in the example, such that
  $y=f_*(x)$.
\end{example}

Analytic K-homology is homotopy invariant, a proof can be found in
\cite{Blackadar(1998)}.
\begin{theorem}
  If $\phi_1,\phi_2\colon X\to Y$ are proper $\Gamma$-equivariant maps 
  which are homotopic through proper $\Gamma$-equivariant maps,
  then
  \begin{equation*}
    (\phi_1)_*=(\phi_2)_*\colon KK_*^\Gamma(X)\to KK_*^\Gamma(Y).
  \end{equation*}
\end{theorem}

\begin{theorem}
  If $\Gamma$ acts freely on $X$, then
  \begin{equation*}
    KK_*^\Gamma(X)\iso K_*(\Gamma\backslash X),
  \end{equation*}
  where the right hand side is the ordinary K-homology of
  $\Gamma\backslash X$.
\end{theorem}

\begin{definition}
  Assume $Y$ is an arbitrary proper $\Gamma$-CW-complex. Set
  \begin{equation*}
    RK^\Gamma_*(Y):= \dirlim KK^\Gamma_*(X),
  \end{equation*}
  where we take the direct limit over the direct system of
  $\Gamma$-invariant subcomplexes of $Y$ with compact quotient (by the
  action of $\Gamma$).
\end{definition}

\begin{definition}
  To define higher (analytic) equivariant K-homology, there are two
  ways. The short one only works for complex K-homology. One
  considers cycles and an
  equivalence relation exactly as above --- with the notable exception 
  that one does not require any grading! This way, one defines
  $KK^\Gamma_1(X)$. Because of Bott periodicity (which has period $2$), 
  this is enough to define all K-homology groups ($KK_n^\Gamma(X)=
  KK_{n+2k}^\Gamma(X)$ for any $k\in\integers$).

  A perhaps more conceptual approach is the following. Here, one
  generalizes the notion of a graded Hilbert space by the notion of a
  $p$-multigraded Hilbert space ($p\ge 0$). This means that the graded 
  Hilbert space comes with $p$ unitary operators
  $\epsilon_1,\dots,\epsilon_p$ which are odd with respect to the
  grading, which satisfy $\epsilon_i^2=-1$ and
  $\epsilon_i\epsilon_j+\epsilon_j\epsilon_i =0$ for all $i$ and $j$
  with $i\ne j$.
  An operator $T\colon H\to H$ on a $p$-multigraded Hilbert space is
  called \emph{multigraded} if it commutes with
  $\epsilon_1,\dots,\epsilon_p$. Such operators can (in addition) be
  even or odd.

  This definition can be reformulated as saying that a multigraded
  Hilbert space is a (right) module over 
  the Clifford algebra $Cl_p$, and a multigraded operator is a module
  map.

  We now define $KK_p^\Gamma(X)$ using cycles as above, with the
  additional assumption that the Hilbert space is $p$-graded, that the 
  representation $\pi$ takes values in $\pi$-multigraded even operators,
  and that the operator $F$ is an odd $p$-multigraded
  operator. Isomorphism and equivalence of these multigraded cycles is 
  defined as above, requiring that the multigradings are
  preserved throughout.
\end{definition}

This definition gives an equivariant homology theory if we restrict
to \emph{proper} maps. Moreover, it satisfies Bott
periodicity. The
period is two for the (complex) K-homology we have considered so
far. All results mentioned in this section generalize to higher
equivariant K-homology.

If $X$ is a proper $\Gamma$-CW-complex, the analytically defined representable
equivariant $K$-homology groups $RK_p^\Gamma(X)$ are canonically
isomorphic to the equivariant $K$-homology groups $K^\Gamma_p(X)$
defined by Davis and L\"uck in \cite{Davis-Lueck(1998)} as described
in Section \ref{sec:homot-ther-defin}.

\subsection{The assembly map}
Here, 
  we will use the homotopy theoretic description of equivariant K-homology
  due to Davis and L\"uck \cite{Davis-Lueck(1998)} described in
  Section \ref{sec:homot-ther-defin}.
  The assembly map then becomes
  particularly convenient to describe. From the present point of view,
  the main virtue
  is that they define a functor from \emph{arbitrary},
  not necessarily proper, $\Gamma$-CW-complexes to abelian groups.

  The Baum-Connes assembly map is now simply defined
  using the equivariant collapse $E(\Gamma,fin)\to *$:
  \begin{equation}\label{eq:Davis-Lueck_BC}
   \mu\colon K^\Gamma_k(E(\Gamma,fin)) \to K^\Gamma_k(*) = K_k(C^*_r\Gamma).
 \end{equation}
 If $\Gamma$ is torsion-free, then $E\Gamma=E(\Gamma,fin)$, and the
 assembly map of \eqref{eq:BC} is defined as the composition of
 \eqref{eq:Davis-Lueck_BC} with the appropriate isomorphism in
 Proposition \ref{prop:Davis_Lueck_properties}.

\subsection{Survey of KK-theory}

The analytic definition of $\Gamma$-equivariant K-homology can be extended to a 
bivariant functor on $\Gamma$-$C^*$-algebras. Here, a
$\Gamma$-$C^*$-algebra is a $C^*$-algebra $A$ with an action (by
$C^*$-algebra automorphisms) of $\Gamma$. If $X$ is a proper
$\Gamma$-space, $C_0(X)$ is such a $\Gamma$-$C^*$-algebra.

Given two $\Gamma$-$C^*$-algebras $A$ and $B$, Kasparov defines the
bivariant KK-groups $KK^\Gamma_*(A,B)$.
The most important property of this bivariant KK-theory is that it
comes with a (composition) product, the \emph{Kasparov product}. This can be
stated most conveniently as follows:

Given a discrete group $\Gamma$, we have a category $KK^\Gamma$ whose
objects are $\Gamma$-$C^*$-algebras (we restrict here to separable
$C^*$-algebras). The morphisms in this category between two
$\Gamma$-$C^*$-algebras $A$ and $B$ are called
$KK_*^\Gamma(A,B)$. They are $\integers/2$-graded abelian groups, and
the composition preserves the grading, i.e.~if $\phi\in
KK_i^\Gamma(A,B)$ and $\psi\in KK_j^\Gamma(B,C)$ then $\psi\phi\in
KK^\Gamma_{i+j}(A,C)$.

There is a functor from the category of separable
$\Gamma$-$C^*$-algebras (where morphisms are $\Gamma$-equivariant
$*$-homomorphisms) to the category $KK^\Gamma_*$ which maps an object
$A$ to $A$, and such that the image of a morphism $\phi\colon A\to B$
is contained in $KK^\Gamma_0(A,B)$.

If $X$ is a proper cocompact $\Gamma$-CW-complex then (by definition)
\begin{equation*}
KK^\Gamma_p(C_0(X),\complexs)= KK^\Gamma_{-p}(X).
\end{equation*}
Here, $\complexs$
has the trivial $\Gamma$-action.

On the other hand, for any $C^*$-algebra $A$ without a group action
(i.e.~with trivial action of hte trivial group $\{1\}$),
$KK^{\{1\}}_*(\complexs,A)= K_*(A)$.

There is a functor from $KK^\Gamma$ to $KK^{\{1\}}$, called \emph{descent},
which assigns to every $\Gamma$-$C^*$-algebra $A$ the \emph{reduced
  crossed product} $C^*_r(\Gamma,A)$. The crossed product has the
property that
$C^*_r(\Gamma,\complexs)=C^*_r\Gamma$.

\subsection{KK assembly}

We now want to give an account of the analytic definition of the
assembly map, which was the original definition. The basic idea is
that the assembly map is given by taking an index. To start with,
assume that we have an even generalized elliptic $\Gamma$-operator
$(H,\pi,F)$, representing an element in $K^\Gamma_0(X)$, where $X$ is
a proper $\Gamma$-space such that $\Gamma\backslash X$ is compact. The
index of
this operator should give us an element in $K_0(C^*_r\Gamma)$. Since
the cycle is even, $H$ split as $H= H_0\oplus H_1$, and $F=\left(
  \begin{smallmatrix}
    0 & P\\ P^* & 0
  \end{smallmatrix}\right)$ with respect to this splitting. Indeed,
now, the kernel and cokernel of $P$ are modules over
$\complexs\Gamma$, and should, in most cases, give modules over
$C^*_r\Gamma$.

If $\Gamma$ is finite, the latter is indeed the case (since
$C^*_r\Gamma=\complexs\Gamma$). Moreover, since $\Gamma\backslash X$
is compact and $\Gamma$ is finite, $X$ is compact, which implies that
$C_0(X)$ is unital. We may then assume that $\pi$ is unital (switching 
to an equivalent cycle with Hilbert space $\pi(1)H$, if
necessary). But then the axioms for a cycle imply that $F^2-1$ is
compact, i.e.~that $F$ is invertible modulo compact operators, or that
$F$ is Fredholm, which means that $\ker(P)$ and $\ker(P^*)$
are finite dimensional. Since $\Gamma$ acts on them,
$[\ker(P)]-[\ker(P^*)]$ defines an element of the representation ring
$R\Gamma=K_0(C^*_r\Gamma)$ for the finite group $\Gamma$. It remains
to show that this map respects the equivalence relation defining
$K_0^\Gamma(X)$. 

However, if $\Gamma$ is not finite, the modules $\ker(P)$ and
$\ker(P^*)$, even if they are $C^*_r\Gamma$-modules, are
in general not finitely generated projective.

To grasp the difficulty, consider Example
\ref{ex:Atiyah-operator_on_covering}. Using the description where $F$
acts on a bundle over the base space $M$ with infinite dimensional
fiber $L\tensor E$, we see that loosely speaking, the null space of
$F$ should rather
``contain'' certain copies of $l^2\Gamma$ than copies of $C^*_r\Gamma$ 
(for finite groups, ``accidentally'' these two are the
same!). However, in general $l^2\Gamma$ is not projective over
$C^*_r\Gamma$ (although it is a module over this algebra). To be
specific, assume that $M$ is a point, $E_0=\complexs$ and $E_1=0$, and
$D=0$. Here we obtain, $L^2(E_0)=l^2\Gamma$, $L^2(E_1)=0$, $F=0$, and
indeed, $\ker(P)=l^2\Gamma$.

In the situation of our example, there is a way around this problem:
Instead of twisting the operator $D$ with the flat bundle
$l^2(\Gamma)\times_\Gamma \overline M$, we twist with
$C^*_r(\Gamma)\times_\Gamma\overline M$, to obtain an operator $D'$
acting on a bundle with fiber $C^*_r\Gamma\tensor \complexs^{\dim
  E}$. This way, we replace
$l^2\Gamma$ by $C^*_r\Gamma$ throughout. Still, it is not true in
general that the kernels we get in this way are finitely generated
projective modules over $C^*_r\Gamma$. However, it is a fact that one
can always add to the new
$F'$ an appropriate compact operator such that this is the case. Then
the
obvious definition gives an element
\begin{equation*}
\ind(D')\in K_0(C^*_r\Gamma).
\end{equation*}
This is the Mishchenko-Fomenko index
of $D'$ which does not depend on the chosen compact
perturbation. Mishchenko and Fomenko give a formula for this index
extending the Atiyah-Singer index formula.

One way to get around the difficulty in the general situation (not
necessarily studying a lifted differential operator) is to deform
$(H,\pi,F)$ to an equivalent
$(H,\pi,F')$ which is better behaved (reminiscent to the compact
perturbation above). This allows to proceeds with a rather elaborate
generalization of the Mishchenko-Fomenko
example we just considered, essentially replacing
$l^2(\Gamma)$ by $C^*_r\Gamma$ again. In this way, one defines an
index as an element of
$K_*(C^*_r\Gamma)$.

This gives a homomorphism $\mu^\Gamma\colon KK_*^\Gamma(C_0(X))\to
K_*(C^*_r\Gamma)$ for each proper $\Gamma$-CW-complex $X$ where
$\Gamma\backslash X$ is compact. This passes 
to direct limits and defines, in particular,
\begin{equation*}
  \mu_*\colon RK^\Gamma_*(E(\Gamma,fin))\to K_*(C^*_r\Gamma).
\end{equation*}

Next, we proceed with an alternative definition of the Baum-Connes map
using KK-theory
and the Kasparov product. The basic observation here is that, given
any proper $\Gamma$-CW-space $X$, there is
a specific projection $p\in C^*_r(\Gamma,C_0(X))$ (unique up to an
appropriate
type of homotopy) which gives rise to a canonical element $[L_X]\in
K_0(C^*_r(\Gamma,C_0(X)))=
KK_0(\complexs,C^*_r(\Gamma,C_0(X)))$. This defines by composition the 
homomorphism
\begin{multline*}
 KK^\Gamma_*(X)=  KK^\Gamma_*(C_0(X),\complexs)\xrightarrow{\text{descent}}
 KK_*(C^*_r(\Gamma,C_0(X)),C^*_r\Gamma)\\
   \xrightarrow{[L_X]\circ\cdot} 
  KK_*(\complexs,C^*_r\Gamma) =K_*(C^*_r\Gamma).
\end{multline*}
Again, this passes to direct limits and defines as a special case the
Baum-Connes assembly map
\begin{equation*}
\mu\colon  RK_*^\Gamma(E(\Gamma,fin)) \to K_*(C^*_r\Gamma).
\end{equation*}

\begin{remark}
  It is a non-trivial fact (due to Hambleton and Pedersen
  \cite{Hambleton-Pedersen(2001)}) that this
  assembly map coincides with the map $\mu$ of \eqref{eq:BC_general}.
\end{remark}

Almost all positive results about the Baum-Connes have been obtained
using the powerful methods of KK-theory, in particular the so called
Dirac-dual Dirac method, compare e.g.~\cite{Valette(2001)}.

\subsection{The status of the conjecture}\label{sec:status-BC-conjecture}

The Baum-Connes conjecture is known to be true for the following classes of
groups.
\begin{enumerate}
\item discrete subgroups of $SO(n,1)$ and $SU(n,1)$ \cite{Julg-Kasparov(1995)}
\item Groups with the \emph{Haagerup property}, sometimes called
  \emph{a-T-menable groups}, i.e.~which admit an
  isometric action on some affine Hilbert $H$ space which is proper,
  i.e.~such that $g_nv\xrightarrow{n\to\infty} \infty$ for every $v\in 
  H$ whenever $g_n\xrightarrow{n\to\infty} \infty$ in $G$
  \cite{Higson-Kasparov(1997)}. 
  Examples of groups with the Haagerup property are amenable groups,
  Coxeter groups, groups acting properly on trees, and groups acting
  properly on simply connected CAT(0) cubical complexes
\item One-relator groups, i.e.~groups with a presentation
  $G=\innerprod{g_1,\dots,g_n\mid r}$ with only one defining relation
  $r$ \cite{Beguni-Bettaieb-Valette(1999)}.
\item Cocompact lattices in $Sl_3(\reals)$, $Sl_3(\complexs)$ and
  $Sl_3(\rationals_p)$ ($\rationals_p$ denotes the $p$-adic numbers) \cite{Lafforgue(1999)}
\item Word hyperbolic groups in the sense of Gromov \cite{Mineyev-Yu(2001)}.
\item Artin's full braid groups $B_n$ \cite{Schick(2001a)}.
\end{enumerate}

Since we will encounter amenability later on, we recall
the definition here.
\begin{definition}\label{def:amenable}
  A finitely generated discrete group $\Gamma$ is called amenable, if
  for any given finite set of generators $S$ (where we require $1\in
  S$ and require that $s\in
  S$ implies $s^{-1}\in S$) there exists a sequence of finite subsets $X_k$
  of $\Gamma$ such that
  \begin{equation*}
    \frac{\abs{SX_k:=\{sx\mid s\in S, x\in
        X_k\}}}{\abs{X_k}}\xrightarrow{k\to\infty} 1.
  \end{equation*}
  $\abs{Y}$ denotes the number of elements of the set $Y$.
  
  An arbitrary discrete group is called amenable, if each finitely
  generated subgroup is amenable.

  Examples of amenable groups are all finite groups, all abelian,
  nilpotent and solvable groups. Moreover, the class of amenable
  groups is closed under taking subgroups, quotients,  extensions, and 
  directed unions.

  The free group on two generators is not amenable. ``Most'' examples of
  non-amenable groups do contain a non-abelian free group.
\end{definition}

There is a certain stronger variant of the Baum-Connes conjecture, the
\emph{Baum-Connes conjecture with coefficients}. It has the following
stability properties:
\begin{enumerate}
\item If a group $\Gamma$ acts on a tree such that the stabilizer of
  every edge and every vertex satisfies the Baum-Connes conjecture
  with coefficients,
  the same is true for $\Gamma$ \cite{Oyono(1998)}.
\item If a group $\Gamma$ satisfies the Baum-Connes conjecture with
  coefficients, then so does every subgroup of $\Gamma$ \cite{Oyono(1998)}
\item If we have an extension $1\to\Gamma_1\to \Gamma_2\to\Gamma_3\to
  1$, $\Gamma_3$ is torsion-free and $\Gamma_1$ as well as $\Gamma_3$
  satisfy the Baum-Connes conjecture with coefficients, then so does $\Gamma_2$.
\end{enumerate}

It should be remarked that in the above list, all groups except for
word hyperbolic groups, and cocompact subgroups of $Sl_3$ actually
satisfy the Baum-Connes
conjecture with coefficients.

The Baum-Connes assembly map $\mu$ of \eqref{eq:BC_general} is
known to be rationally injective for considerably larger classes of
groups, in particular the following.
\begin{enumerate}
\item Discrete subgroups of connected Lie groups
  \cite{Kasparov(1995)}
\item Discrete subgroups of $p$-adic groups
  \cite{Kasparov-Skandalis(1991)}
\item Bolic groups (a certain generalization of word hyperbolic
  groups) \cite{Kasparov-Skandalis(1994)}.
\item Groups which admit an amenable action on some compact space
\cite{Higson-Roe(2000)}.
\end{enumerate}

Last, it should be mentioned that recent constructions of Gromov show
that certain variants of the Baum-Connes conjecture, among them the
Baum-Connes conjecture with coefficients, and an extension called the
\emph{Baum-Connes conjecture for groupoids}, are false
\cite{Higson-Lafforgue-Skandalis(2001)}. At the
moment, no counterexample to the Baum-Connes conjecture
\ref{conj:general_BC} seems to be known. However, there are many
experts in the field who think that such a counterexample eventually
will be constructed \cite{Higson-Lafforgue-Skandalis(2001)}.

\section{Real $C^*$-algebras and K-theory}\label{sec:real-c-algebras}

\subsection{Real $C^*$-algebras}

The applications of the theory of $C^*$-algebras to geometry and
topology we present here require at some point that we work with real
$C^*$-algebras. Most of the theory is parallel to the theory of
complex $C^*$-algebras.

\begin{definition}
  A unital real $C^*$-algebra is a Banach-algebra $A$ with unit over
  the real numbers,
  with an isometric involution $*\colon A\to A$, such that
  \begin{equation*}
    \abs{x}^2= \abs{x^*x}\qquad\text{and }1+x^*x\text{ is
      invertible}\quad\forall x\in A.
  \end{equation*}

  It turns out that this is equivalent to the existence of a
  $*$-isometric embedding of $A$ as a closed subalgebra into
  $\boundedops{H_\reals}$, the bounded operators on a suitable real
  Hilbert space (compare \cite{Palmer(1970)}).
\end{definition}

\begin{example}
  If $X$ is a compact topological space, then $C(X;\reals)$, the
  algebra of real valued continuous function on $X$, is a real
  $C^*$-algebra with unit (and with trivial involution).

  More generally, if $X$ comes with an involution $\tau\colon X\to X$
  (i.e. $\tau^2=\id_X$), then $C_\tau(X):=\{f\colon X\to\complexs\mid
  f(\tau x)=\overline{f(x)}\}$ is a real $C^*$-algebra with involution 
  $f^*(x)=\overline{f(\tau x)}$.

  Conversely, every commutative unital real $C^*$-algebra is
  isomorphic to some $C_\tau(X)$.

  If $X$ is only locally compact, we can produce examples of
  non-unital real $C^*$-algebras as in Example \ref{ex:loc_comp_Cstar_algebra}.
\end{example}

Essentially everything we have done for (complex) $C^*$-algebras
carries over to real $C^*$-algebras, substituting $\reals$ for
$\complexs$ throughout. In particular, the definition of the K-theory
of real $C^*$-algebras is literally the same as for complex
$C^*$-algebras (actually, the definitions make sense for even more
general topological algebras), and a short exact sequence of real
$C^*$-algebras gives rise to a long exact K-theory sequence.

The notable exception is Bott periodicity. We don't get the period
$2$, but the period $8$.

\begin{theorem}
  Assume that $A$ is a real $C^*$-algebra. Then we have a Bott
  periodicity isomorphism
  \begin{equation*}
    K_0(A)\iso K_0(S^8A).
  \end{equation*}
  This implies 
  \begin{equation*}
    K_n(A) \iso K_{n+8}(A)\qquad\text{for }n\ge 0.
  \end{equation*}
\end{theorem}

\begin{remark}
  Again, we can use Bott periodicity to define $K_n(A)$ for arbitrary
  $n\in\integers$, or we may view $K_n(A)$ as an $8$-periodic theory,
  i.e.~with $n\in\integers/8$.

  The long exact sequence of Theorem \ref{theo:long_exact} becomes a
  24-term cyclic exact sequence.
\end{remark}

The \emph{real reduced $C^*$-algebra} of a group $\Gamma$, denoted
$C^*_{\reals,r}\Gamma$, is the norm closure of $\reals\Gamma$ in the
bounded operators on $l^2\Gamma$.

\subsection{Real K-homology and Baum-Connes}

A variant of the cohomology theory given by complex vector bundles is
KO-theory, which is given by real vector bundles. The homology theory
dual to this is KO-homology. If $KO$ is the spectrum of topological
KO-theory, then $KO_n(X) = \pi_n(X_+\wedge KO)$.

The homotopy theoretic definition of equivariant K-homology can be
varied easily to define equivariant KO-homology. The analytic
definition can also be adapted easily, replacing $\complexs$ by
$\reals$ throughout, using in particular real Hilbert spaces. However, 
we have to stick to $n$-multigraded cycles to define $KK^\Gamma_n(X)$, 
it is not sufficient to consider only even and odd cycles.

All the constructions and  properties translate appropriately from the 
complex to the real situation, again with the notable
exception that Bott periodicity does not give the period $2$, but
$8$. The upshot of all of this is that we get a real version of the
Baum-Connes conjecture, namely
\begin{conjecture}
  The real Baum-Connes assembly map
  \begin{equation*}
    \mu_n\colon KO^\Gamma_n(E(\Gamma,fin)) \to KO_n(C^*_{\reals,r}\Gamma),
  \end{equation*}
  is an isomorphism.
\end{conjecture}

It should be remarked that all known results about injectivity or
surjectivity of the Baum-Connes
map can be proved for the real version as well as for the
complex version, since each proof translates without
too much difficulty. Moreover, it is known that the complex version of
the
Baum-Connes conjecture for a group $\Gamma$ implies the real version
(for this abstract result, the isomorphism is needed as input, since
this is based on the use of the five-lemma at a certain point).

\chapter[Counterexample to GLR]{A counterexample to the Gromov-Lawson-Rosenberg conjecture}

The Gromov-Lawson-Rosenberg conjecture is discussed in the
notes by Stephan Stolz. As a reminder, we quickly recall the
problem:

\begin{question}
  Given a compact smooth spin-manifold $M$ without boundary, when does
  $M$ admit a Riemannian metric with positive scalar curvature?
\end{question}

Recall that a spin-manifold is a manifold for which the first and
second Stiefel-Whitney class of the tangent bundle vanish. The spin
condition can be compared to the condition that a manifold is
orientable. Indeed, every spin-manifold is orientable. But the spin
condition is considerably
stronger (it is like orientability ``squared'').

 The reason that we concentrate on spin-manifolds is that
 powerful obstructions to the existence of a metric with positive
 scalar curvature have been developed  for them.

\section{Obstructions to positive scalar curvature}

\subsection{Index theoretic obstructions}

We start with a discussion of
 the index obstruction for spin manifolds to admit a metric with $\scal>0$, 
constructed by
Lichnerowicz \cite{Lichnerowicz(1963)}, Hitchin \cite{Hitchin(1974b)} and in
the following refined version due to Rosenberg \cite{Rosenberg(1987)}.
\begin{theorem}\label{theo:index_homomorphism}
One can construct a
homomorphism, called index, from the singular spin bordism group
$\Omega_*^{spin}(B\pi)$ to the (real) $KO$-theory of the reduced real
$C^*$-algebra of $\pi$:
\[  \ind\colon \Omega_*^{spin}(B\pi)\to KO_*(C^*_{\reals,r}\pi)\]
 (this homomorphism is often called $\alpha$
instead of $\ind$). Assume $f\colon N\to B\pi$ represents an element of
$\Omega_m^{spin}(B\pi)$.
If $N$ admits a metric with positive scalar curvature, then
\[ \ind([f\colon N\to B\pi]) =0\in KO_m(C^*_{\reals,r}\pi) \]
\end{theorem}

The converse of this theorem is the content of the following
\emph{Gromov-Lawson-Rosenberg conjecture}.
\begin{conjecture}\label{conj:GLR}
 Let $M$ be a
  compact spin-manifold without boundary, $\pi=\pi_1(M)$, and $u\colon 
  M\to
  B\pi$ be the classifying map for a universal covering
  of $M$. Assume that $m=\dim(M)\ge 5$.

  Then $M$ admits a metric with $\scal>0$ if and only if 
  \begin{equation*}
\ind[u\colon 
  M\to B\pi] =0 \in KO_m(C^*_{\reals,r}\pi).
\end{equation*}
\end{conjecture}
This conjecture was developed in \cite{Gromov-Lawson(1983)} and
\cite{Rosenberg(1983)}. 

The restriction to dimensions $\ge 5$ comes from the
observation that in these dimensions (and not below) the question of
existence of metrics with $\scal>0$ in a certain sense is a bordism
invariant, which of course fits with the structure of the obstruction
described in Theorem \ref{theo:index_homomorphism}.
Failure of this bordism invariance in dimension $4$ is also reflected
by the fact that for $4$-dimensional manifolds, 
the \emph{Seiberg-Witten invariants} provide additional obstructions
to the existence of a metric with positive scalar curvature, which
show in particular that the conjecture is not true if $m=4$.

 The conjecture was proved by Stefan Stolz \cite{Stolz(1992)}
 for $\pi=1$, and subsequently by him and other authors
 also for some other groups 
\cite{Rosenberg(1983),Kwasik-Schultz(1990),%
Botvinnik-Gilkey-Stolz(1995),Rosenberg-Stolz(1995),Joachim-Schick(2000)}.

\subsection{Minimal surface obstructions}

In dimension $\ge 5$ there is only one known
additional obstruction for positive
scalar curvature metrics,
the minimal surface method of  Schoen and Yau, which we will recall now.
(In dimension $4$, the Seiberg-Witten theory yields additional obstructions).

The first 
theorem is the differential geometrical backbone for the application of
minimal surfaces to the positive scalar curvature problem:
\begin{theorem}\label{diffgeo}
Let $(M^m, g)$ be a manifold with $\scal>0$,
 $\dim M=m\ge 3$.
 If $V$ is a smooth $(m-1)$-dimensional
submanifold of $M$ with trivial normal bundle,
 and if $V$ is a local minimum of the
volume functional, then $V$ admits a metric of positive scalar
curvature, too. ``Local minimum'' means that for any
small deformation of the hypersurface, the $(m-1)$-volume of the surface 
increases.

Actually, $V$ can be a ``minimal hypersurfaces'' in the
sense of differential geometry, defined in terms of curvature and second
fundamental form of the hypersurface. Every local minimum for the
$(m-1)$-volume is such a minimal hypersurface; the converse is not true.
\end{theorem}
\begin{proof}
Schoen/Yau: \cite[5.1]{Schoen-Yau(1979a)} for $m=3$, 
\cite[proof of Theorem 1]{Schoen-Yau(1979)} for $m>3$.
 We outline the proof, following closely \cite[Theorem
 1]{Schoen-Yau(1979)}.

 Given $V$, since its normal bundle is trivial, any smooth function $\phi$ on $V$
 gives rise to a variation (if $\nu$ is the unit normal vector field,
 pushing $V$ in normal direction $\phi\nu$). Lt $\Ric\in
 \Gamma(\End(TM))$ be the Ricci curvature, considered as an operator
 on each fiber of $TM$. Let $l$ be the second
 fundamental form of $V$. It is well known that minimality implies
 $\tr(l)=0$. Moreover, the 
 second variation of the area is non-negative. It is given by (see
 \cite{Chern(1968)})
 \begin{equation}\label{eq:variation}
   - \int_V \left(\innerprod{\Ric(\nu),\nu} + \abs{l}^2\right) \phi^2
   +\int_V \abs{\nabla \phi}^2 \ge 0.
 \end{equation}
 We now use the Gauss curvature equation (the ``theorema egregium'')
 to relate this to the scalar curvature of the submanifold. Taking
 appropriate traces of the Gauss equations, we obtain
 \begin{equation}
   \label{eq:theorema_egregium}
   \scal_V = \scal_M - 2 \innerprod{\Ric(\nu),\nu} + (\tr l)^2 - \abs{l}^2,
 \end{equation}
 where $\scal_V$ is the scalar curvature of $V$ with the induced
 Riemannian metric and $\scal_M$ the scalar curvature of $M$.
Putting Equation \eqref{eq:theorema_egregium} into Inequality
\eqref{eq:variation}, we have
\begin{equation}\label{eq:Ric}
  \int_V \scal_M\phi^2 - \int_V \scal_V \phi^2 +\int_V \abs{l}^2\phi^2 
  \le 2\int_V \abs{\nabla\phi}^2
\end{equation}
for all smooth functions $\phi\colon V\to \reals$. We assume that
the scalar curvature of $M$ is everywhere strictly positive. Hence
\eqref{eq:Ric} implies
\begin{equation*}
  -\int_V \scal_V \phi^2 <2\int_V \abs{\nabla\phi}^2,
\end{equation*}
as long as $\phi$ is not identically zero.

Consider the \emph{conformal Laplacian} $\Delta_c:=\Delta +\frac{n-3}{4(n-2)} 
\scal_V$ on $V$ (where $\Delta$ is the positive Laplacian on
functions). Then all eigenvalues of $\Delta_c$ are strictly
positive. Assume, otherwise, that $\phi$ is an eigenfunction to the
eigenvalue $\lambda\le 0$, i.e.
\begin{equation*}
  \Delta\phi = - \frac{m-3}{4(m-2)} \scal_V \phi +\lambda\phi.
\end{equation*}
Taking the $L^2$-inner product of this equation with $\phi$ (and
integration by parts) gives
\begin{equation*}
  \int_V \abs{\nabla\phi}^2 = -\frac{m-3}{4(m-2)} \int_V\scal_V \phi^2 
  +\lambda\int_V \phi^2 < \frac{m-3}{2(m-2)} \int_V\abs{\nabla\phi}^2,
\end{equation*}
which is a contradiction. Now it's a standard fact in conformal
geometry that, if the conformal Laplacian has only positive
eigenvalues, then one can conformally deform the metric to a metric
with positive scalar curvature (compare
\cite{Kazdan-Warner(1975)}). 
This is done in two steps: a generalized maximum principle implies
that we can find an eigenfunction $f$ to the first eigenvalues of
$\Delta_c$ which is strictly positive everywhere.
Then, explicit formulas for the scalar curvature of a conformally
changed metric show that $f^{4/(m-3)}g$ indeed has a metric with $\scal>0$.

Hence, on $V$ there exists a metric
with $\scal>0$ (observe, however, that it is not necessarily the
metric induced from $M$, but only conformally equivalent to this metric).
\end{proof}

The next statement due to Simons and Smale (special cases 
due to Fleming and Almgren) from geometric measure theory implies
applicability of the
previous theorem if $\dim(M)\le 8$.

\begin{theorem}\label{geomeasure}
  Suppose $M$ is an orientable Riemannian manifold of dimension $\dim M=m\le
  8$. Furthermore let $\alpha\in H^1(M,\integers)$. Then
  \[ x:=\alpha\cap [M]\in H_{m-1}(M,\integers) \]
can be represented by an embedded hypersurface $V$ with trivial 
normal bundle which is a local minimum for
$(m-1)$-volume (if $m=8$ with respect to suitable metrics arbitrarily
close in $C^3$ to the metric we started with).
\end{theorem}
\begin{proof}
For $m\le 7$
this is a classical result of geometric measure
theory (cf. \cite[Chapter 8]{Morgan(1988)}) and references therein,
in particular \cite[5.4.18]{Federer(1969)}. 

The case $m=8$ follows 
from the following result of Nathan Smale \cite{Smale(1993)}: the 
set of $C^k$-metrics for which the regularity statement holds is open
and dense in the set
of all $C^k$-metrics ($k\ge3$ and with the usual Banach-space topology). We are
only interested in
$C^\infty$-metrics. But these are dense in the set of $C^k$-metrics,
which concludes the proof.
\end{proof}

Unfortunately, the proofs of the theorems we have cited are very
involved and require a lot of technical work. Therefore, we don't
attempt to indicate the arguments here.

Recall that if we are given a class $\alpha \in H^{1}(M,\integers)$ we may 
represent it by a map $f\colon M \to S^{1}$ being transverse to $1\in S^{1}$.
Then $V=f^{-1}(1) \subset M$  represents $\alpha \cap [M]$ (and
conversely, every hypersurface representing $\alpha \cap [M]$ is
obtained in this way). Furthermore,
if $f'\colon M\to S^{1}$ is a second map as above and $V' = {f'}^{-1}(1)$ then 
$f$ and $f'$ are homotopic, and a homotopy $H\colon f\simeq f'$ being transverse to 
$1\in S^{1}$ provides a bordism $W= H^{-1}(1)\colon V \sim V'$ embedded in
$M\times [0,1]$. Since the normal bundle
of $V \subset M$ and  $W \subset M\times [0,1]$,respectively, is trivial,
the manifolds and bordisms we construct this way belong to the
spin-category, if we start with a spin-manifold $M$.

We want to use these ideas to construct, for an arbitrary space $X$, a map
\begin{equation}
  \label{pairing}
 \cap\colon  H^1(X,\integers)\times\Omega_m^{spin}(X) \to \Omega^{spin}_{m-1}(X) .
\end{equation}

To do this, let $\phi\colon M\to X$ be a singular spin manifold for
$X$, representing an element in $\Omega_m^{spin}(X)$.
If $f\colon X \to S^{1}$ represents an element $\alpha \in
H^{1}(X,\integers)$,
then $f \circ \phi$ is homotopic to a map $\psi\colon M \to S^{1}$
which is transverse to $1\in S^{1}$.

Restricting $\phi$ to $V:=\psi^{-1}(1)$ then gives a singular
spin manifold $\phi|_V\colon V\to X$, which by definition represents
$\alpha\cap [\phi\colon M\to X]\in \Omega^{spin}_{m-1}(X)$.

We have to check that this is well defined. To do this, let
$\Phi\colon W \to X$ be a spin-bordism
between  $\phi\colon M \to X$ and $\phi'\colon M' \to X$. 
Then $f \circ \Phi$ is homotopic to a map $\Psi\colon W \to S^{1}$ with
$\Psi$ and $\Psi|\partial W$ being transverse to $1\in S^{1}$; moreover, the map $\Psi|\partial W$ with the corresponding properties
may be given in advance. 
Then $\Psi^{-1}(1) \subset W$ is a spin-bordism between the 
hypersurfaces $V=(\Psi)^{-1}(1) \cap M \subset M$ and 
$\Psi^{-1}(1) \cap M'\subset M'$, and restricting $\Phi$ to  
$\Psi^{-1}(1)$ now yields a singular spin-bordism between 
singular spin-hypersurfaces into $X$. \\
If $f'\colon X\to S^1$ is homotopic to $f$,  a similar construction
gives a
singular spin-bordism between the resulting  singular
spin-hypersurfaces into $X$. Together, this implies that our map
indeed is well defined.

The two theorems above now imply
\begin{theorem}\label{bord_schoen_yau}
  Let $X$ be a space and let $3\le m\le 8$. Then \eqref{pairing} 
restricts to a homomorphism:
  \begin{equation}
    \label{pairing+}
     \cap\colon H^1(X,\integers)\times  \Omega_m^{spin,+}(X)\to
     \Omega_{m-1}^{spin,+}(X)
  \end{equation}
where $ \Omega_*^{spin,+}(X) \subset  \Omega^{spin}_*(X)$ is the subgroup of 
bordism classes which can be represented by singular manifolds 
which admit a metric with $\scal >0$ (observe, only one representative 
with $\scal>0$ is required).
\end{theorem}
\begin{proof}
  Theorem \ref{geomeasure} implies that, given $f\colon M\to S^1$ (dual to a
  given class in $H_{m-1}(M,\integers)$) we find a
  homotopic map $g\colon M\to S^1$ which is transverse to $1$ and such that
  the hypersurface $V=g^{-1}(1)$ is minimal for the $(m-1)$-volume (in
  dimension $8$ we  replace the given metric
  by one which
  is $C^3$-close). In any case, since the scalar curvature is
  continuous with respect to the $C^3$-topology on the space of all
  Riemannian metrics, $V$ is volume minimizing with respect to a metric with
  positive scalar curvature whenever we start with such a metric. By
  Theorem \ref{diffgeo}, it admits a metric with $\scal>0$.
\end{proof}

To be honest, this does not quite give an obstruction, but rather a
method to produce counterexamples. Namely, if we know
$\Omega_{n-1}^{spin,+}(B\pi)$, and also the cap-product of
\eqref{pairing} well enough, we can get information about
$\Omega_n^{spin,+}(B\pi)$ (with $n\le 8$).

Obviously, one does need some information to start with. This can be
obtained in dimension $2$ using the Gauss-Bonnet theorem.

\subsection{Gauss-Bonnet obstruction in dimension $2$}

\begin{theorem}\label{theo:Gauss-Bonnet}
  Let $G$ be a discrete group. Then
  \begin{equation}
    \label{omega2plus}
    \Omega^{spin,+}_2(BG):=\left\{
       \text{\begin{tabular}{c}
              bordism classes $[M\to BG] \in\Omega^{spin}_2(BG)$,\\
              where $M$ admits a metric with $\scal>0$
             \end{tabular}} 
                                 \right\}= {0} .
  \end{equation}
\end{theorem}
\begin{proof}
 By the Gauss-Bonnet theorem there is only one orientable $2$-manifold with positive
(scalar) curvature, namely $S^2$. On the other hand, $S^2$ is a
spin-manifold with a unique spin-structure, and is spin-bordant to
zero, being
the boundary of $D^3$. Since $\pi_2(BG)$ is trivial, up to homotopy
only the trivial map from $S^2$ to $BG$ exists. Therefore only the trivial
element in $\Omega_2^{spin}(BG)$ can be represented by a manifold with
positive scalar curvature.
\end{proof}

\section{Construction of the counterexample}

\subsection{Application of the minimal hypersurface obstruction}
\label{sec:appl-minim-surf}

Now, we will construct a particular example of a manifold which does
not admit a metric with positive scalar curvature, using the minimal
hypersurface obstruction.

Let $p\colon S^1\to B\integers/3$ be a map so that $\pi_1(p)$ is
surjective and equip $S^1$
with the spin structure induced from $D^2$.
Consider the
singular manifold 
\[ f = \id\times p\colon T^5=S^1\overbrace{\times\cdots\times}^4
S^1\times S^1 \to
 S^1\overbrace{\times\cdots\times}^4 S^1\times B\integers/3 = B\pi, \]
where $\pi=\integers^4\times\integers/3$.
This represents a certain element $x\in \Omega_5^{spin}(B\pi)$.

We have four distinguished maps from $B\pi$ to $S^1$, given by the
projections $p_i\colon B\pi\to S^1$ onto each of the first four factors.
Let $a_1,\dots,a_4\in H^1(B\pi)$ be the corresponding elements in
cohomology. Using the description of the cap-product given before
\eqref{pairing}, in our situation it is easy to find a
representative for
\begin{equation*}
  z :=a_1\cap (a_2\cap (a_3 \cap w))\in \Omega^{spin}_2(B\pi).
\end{equation*}
Namely, taking inverse images of the base point, this $z$ is given by
\begin{equation*}
  g = *\times *\times *\times \id \times p\colon T^2=*\times *\times
  *\times S^1\times S^1 \to
 S^1\times S^1\times S^1\times S^1\times B\integers/3 = B\pi.
\end{equation*}

We want to show that $z\notin \Omega_2^{spin,+}(B\pi)$, because then,
by \eqref{pairing+}, $x\notin \Omega_5^{spin,+}(B\pi)$, i.e.~whenever
we find a representative
$[f\colon M\to B\pi]=x$, then $M$ does not admit a Riemannian metric
with $\scal>0$ (in particular, this follows then for $T^5$, which,
however, is not the manifold we are interested in here).

Because of Theorem \ref{theo:Gauss-Bonnet} we only have to show that
$z$ is a non-trivial element of $\Omega_2^{spin}(B\pi)$. We have the natural
homomorphism $\Omega_*^{spin}(B\pi)\to H_*(B\pi,\integers)$, which
maps $[f\colon M\to B\pi]$ to $f_*[M]$, i.e.~to the image of the
fundamental class of $M$, and the K\"unneth
theorem implies immediately that the image of $z$ under this
homomorphism in $H_2(B\pi)$ is non-trivial, therefore the same is true 
for $z$.

\subsection{Calculation of the index obstruction}

We proceed by proving that the index obstruction
\ref{theo:index_homomorphism} does vanish for the example
constructed in Subsection \ref{sec:appl-minim-surf}.

This index obstruction is an element of $KO_5(C^*_{\reals,r}\pi)$,
where $\pi=\integers^4\times \integers/3$. First, we compute this $K$-theory
group to the extent needed here. By a K\"unneth theorem for the
K-theory of
$C^*$-algebras, the K-theory of $C^*_{\reals,r}(G\times \integers)$ can 
easily be computed from the K-theory of $C^*_{\reals,r}G$.
 Namely, by 
\cite[p.\ 14--15 and 1.5.4]{Schroeder(1993)}
\[ KO_n(C^*_{\reals,r}(\integers^4\times \integers/3))\cong 
\bigoplus_{i=1}^{16} KO_{n-n_i}(C^*_{\reals,r}
(\integers/3));\qquad \text{for suitable $n_i\in\naturals $.} \]
For a finite group $G$, it is well known that  $KO_*(C^*_{\reals,r}(G))$
 is a direct sum
of copies of the (known) $KO$-theories of $\reals$, $\complexs$
and
$\quaternions$. In particular, it is
a direct sum of copies of $\integers$ and $\integers/2$. Therefore,
the same is true for
$\pi$. This implies the following Proposition.
\begin{proposition}\label{prop:no_3_torsion_in_K}
  $KO_*(C^*_{\reals,r}\pi)$ is a direct sum of copies of $\integers$ and
  $\integers/2$. In
particular, its torsion is only $2$-torsion.
\end{proposition}

Let $p\colon S^1\to B\integers/3$ be the map of Subsection
\ref{sec:appl-minim-surf} so that $\pi_1(p)$ is surjective (and $S^1$
is equipped
with the spin structure induced from $D^2$).
This represents a $3$-torsion element $y$ in
$\Omega_1^{spin}(B\integers/3)$ since
\begin{equation*}
\tilde\Omega_1^{spin}(B\integers/3)\cong
 H_1(B\integers/3,\integers)\cong\integers/3
\end{equation*}
(using 
 e.g.~the Atiyah-Hirzebruch spectral sequence). 

It follows that 
\begin{equation*}
x=[\id_{(S^1)^4}\times p\colon T^5\to
B(\integers^4\times\integers/3)]
\end{equation*}
is also $3$-torsion (a zero bordism
for $3x$ is obtained as the product of a zero bordism for $3y$ with
$\id_{(S^1)^4}$).

Since $\ind\colon \Omega_5^{spin}(B\pi)\to KO_5(C^*_{\reals,r}\pi)$ is 
a group homomorphism, 
\begin{equation*}
3\cdot \ind(x)=0\in KO_5(C^*_{\reals,r}\pi).
\end{equation*}
But for 
$\pi=\integers^4\times\integers/3$, by Proposition
\ref{prop:no_3_torsion_in_K} this implies that $\ind(x)=0$, i.e.~the
index obstruction vanishes.

\subsection{Surgery to produce the counterexample}

So far, we have found a bordism class $x\in\Omega_5^{spin}(B\pi)$
($\pi=\integers^4\times\integers/3$) such 
that the index obstruction \ref{theo:index_homomorphism} vanishes for
$x$, but on the other hand no representative $[f\colon M\to B\pi]$ can 
be found such that $M$ has a metric with positive scalar curvature. To 
give a counterexample to Conjecture \ref{conj:GLR}, we have to find a
representative such that $f$ induces an isomorphism on fundamental
groups (i.e.~is the classifying map for the universal covering of
$M$). This is not the case for the tori we have explicitly constructed 
so far (and indeed, for tori one can use the \emph{index method} to show that 
they do not admit a metric with positive scalar curvature).

But adjusting the fundamental group is easy. We only have to perform
surgery on our explicitly given torus $T^5$. That is, we have to
choose an embedded $S^1\to T^5$ which represents the kernel of
$\pi_1(f)\colon \pi_1(T^5)\to \pi_1(B\pi)$ (observe that in this
situation, the kernel actually is cyclic) and which has a trivial
normal bundle. Then a tubular neighborhood of $S^1$ is diffeomorphic
to $S^1\times D^4$, with boundary $S^1\times S^3$. We can now cut away 
this tubular neighborhood and glue in $D^2\times S^3$ (also with
boundary $S^1\times S^3$) instead. The fundamental group of the new
manifold $M'$ is the quotient of the fundamental group of the old
manifold by the (normal) subgroup generated by the loop we started
with, i.e.~is isomorphic to $\pi$. Let $u\colon M'\to B\pi$ be the
classifying map for the universal covering. Using classical ``surgery
below the middle dimension'', we can
arrange all this in such a way that 
\begin{equation*}
[f\colon T^5\to B\pi]=[u\colon
M'\to B\pi] \in \Omega_5^{spin}(B\pi)
\end{equation*}
(compare \cite[Lemma
5.6]{Stolz(1998)}). Consequently, $M'$ is a
counterexample to the Gromov-Lawson-Rosenberg conjecture \ref{conj:GLR}.

\section{Other questions, other examples}

The index map of Theorem \ref{theo:index_homomorphism} admits a
factorization
\begin{equation*}
  \ind\colon \Omega_*^{spin}(B\pi)\xrightarrow{D}
  ko_*(B\pi)\xrightarrow{p} KO_*(B\pi) \xrightarrow{\mu}
  KO_*(C^*_{r,\reals}\pi).
\end{equation*}

Here, $ko_*$ is connective real $K$-homology, $KO_*$ the periodic real
$K$-homology we have considered so far, $D$ is the $ko$-theoretic
orientation, $p$ the canonical
map between the connective and the periodic theory, and $\mu$ the
assembly map in topological $K$-theory. Note that for torsion free
groups, this $\mu$ is the Baum-Connes map, and the Baum-Connes
conjecture
states that this map is an isomorphism.

The original conjecture of Gromov and Lawson asserted that the
vanishing of the image
of $[u\colon M\to B\pi_1(M)]$ in $KO_m(B\pi)$ decides whether $M$
admits a metric with positive scalar curvature. Rosenberg observed
that there are manifolds with $\scal>0$ for which this element does
not vanish, and proposed to modify the conjecture as stated in
Conjecture \ref{conj:GLR}. We adopt the convention that $u\colon M\to
B\pi_1(M)$ denotes the classifying map for the universal covering of $M$.

However, the following question remains.
\begin{question}\label{question:strong_vanish}
  Is the stronger vanishing condition that 
  \begin{equation*}
pD[u\colon M\to
  B\pi_1(M)]=0
\end{equation*}
sufficient for the existence of metrics with positive
  scalar curvature?
\end{question}
If even $D[u\colon M\to
B\pi_1(M)]=0$, then $M$ admits a metric with positive scalar curvature 
by a result of Stephan Stolz
\cite{Stolz(1994)} (as usual, we have to assume that $\dim(M)\ge 5$).

In \cite{Joachim-Schick(2000)}, a counterexample to question
\ref{question:strong_vanish}
is given. The first step to construct the counterexample is to find a
group such that 
\begin{equation*}
p\colon ko_*(B\pi)\to KO_*(B\pi)
\end{equation*}
has a kernel, and since 
we want to use the minimal surface method, this kernel should be given 
for $*=2$. In \cite{Joachim-Schick(2000)}, this is done using explicit 
$K$-homology calculations for finite groups. The remaining proof is
very much along the lines of the proof we have given above.

One of the virtues of the example we have given is that we avoid the
calculation of the index. This is replaced by some (easy)
considerations about torsion. To be able to do this, we used a
fundamental group $\pi$ with torsion. Dwyer and Stolz (unpublished)
have constructed a counterexample to the Gromov-Lawson-Rosenberg 
conjecture with torsion-free fundamental group. In
\cite{Schick-Stolz(2000)} a refinement of this is given where the
classifying space $B\pi$ is a manifold with negative curvature. The
first key idea is the same as in the example in
\cite{Joachim-Schick(2000)} just described, namely to find an element
in the kernel of $ko_2(B\pi)\to KO_2(B\pi)$. To find a $\pi$ such that 
$B\pi$ is particularly nice (e.g.~finite dimensional, which implies
that $\pi$ is torsion-free, or even a manifold of negative curvature)
one uses asphericalization procedures of Baumslag, Dyer and Heller, or
Charney, Davis,
and Januszkiewicz, which produce nice $B\pi$ with certain prescribed
homological properties (starting with (worse) spaces which have these
same homological properties). More constructions of this kind are
described in the lectures of Mike Davis.

The positive scalar curvature question makes sense also for manifolds
which are not spin manifolds. There are ``twisted'' index obstructions 
as long as the universal covering is a spin manifold, and one can
formulate an appropriate ``twisted Gromov-Lawson-Rosenberg''
conjecture. In \cite{Joachim-Schick(2000)}, counterexamples to this
twisted conjecture are given, as well.

\chapter[$L^2$-cohomology]{$L^2$-cohomology and the conjectures of Atiyah, Singer, and
  Hopf}

$L^2$-cohomology and $L^2$-Betti numbers are certain ``higher
invariants'' of manifolds and more general spaces. They were 
introduced 1976 by Michael Atiyah in \cite{Atiyah(1976)}, and since
then have proved to be useful  invariants with connections and
applications in many other mathematical fields, from differential
geometry to group theory and algebra. Apart from the original
literature, there exists L\"ucks informative survey article
\cite{Lueck(1997)}, and also Eckmann lecture notes
\cite{Eckmann(2000)}. Moreover, at the time of writing of this
article,
L\"uck's very comprehensive textbook/research monograph \cite{Lueck(2001)}
is almost finished, and the current version is available from the author's 
homepage. This chapter is a survey style
article which focuses on the main points of the very extensive
subject, leaving out many of the less illuminating details, which can
be found e.g.~in \cite{Lueck(2001)}.

\section{Analytic $L^2$-Betti numbers}\label{sec:analytic-l2-betti}

\begin{definition}\label{def:L2_sections}
  Let $\overline M$ be a (not necessarily compact) Riemannian manifold 
  without boundary, which is complete as a metric space. Define
  \begin{equation*}
    L^2\Omega^p(\overline M) := \{\omega\text{ measurable $p$-form on
      }M\mid \int_{\overline M} \abs{\omega(x)}^2_x \;d\mu(x) <\infty\}.
  \end{equation*}
  Here, $\abs{\omega(x)}_x$ is the pointwise norm (at $x\in\overline
  M$) of $\omega(x)$, which is given by the Riemannian metric, and
  $d\mu(x)$ is the measure induced by the Riemannian metric.

  $L^2\Omega^p(\overline M)$ can be considered as the Hilbert space
  completion of the space of compactly supported $p$-forms on
  $\overline M$. The inner product is given by integrating the
  pointwise inner product, i.e.
  \begin{equation*}
    \innerprod{\omega,\eta}_{L^2} := \int_{\overline M}
    \innerprod{\omega(x),\eta(x)}_x\; d\mu(x).
  \end{equation*}
\end{definition}

\begin{definition}\label{def:L2cohomology}
  Let $M$ be a smooth compact Riemannian manifold without boundary,
  with Riemannian metric $g$. Let $\overline M$ be a normal covering
  of $M$, i.e.~if $\Gamma$ is the deck transformation group, then
  $M=\overline M/\Gamma$. Lift the metric $g$ to $\overline M$. Then
  $\Gamma$ acts isometrically on $\overline M$. Let
  $\overline\Delta_p$ be the Laplacian on $p$-forms on $\overline
  M$. This gives rise to an unbounded operator
  \begin{equation*}
    \overline\Delta_p\colon L^2\Omega^p(\overline M) \to
    L^2\Omega^p(\overline M).
  \end{equation*}
  This operator is an elliptic differential operator (but on the not
  necessarily compact manifold $\overline M$) Let 
  \begin{equation*}
    \pr_p\colon  L^2\Omega^p(\overline M) \to
    L^2\Omega^p(\overline M)
  \end{equation*}
  be the orthogonal projection onto
  $\ker(\overline\Delta_p)$. Ellipticity of $\overline\Delta_p$
  implies that $\pr_p$ has a smooth integral kernel, i.e.~that there
  is a smooth section $\pr_p(x,y)$ over $\overline M\times\overline M$ 
  (of the bundle with fiber $\Hom(\Lambda^p T^*_yM,\Lambda^p T^*_x M)$ 
  over $(x,y)\in\overline M\times \overline M$), such that 
  \begin{equation*}
    \pr_p\omega(x) = \int_{\overline M} \pr_p(x,y)\omega(y) \;d\mu(y)
  \end{equation*}
  for every $L^2$-$p$-form on $\overline M$.

  We define the $L^2$-cohomology of $\overline M$ by
  \begin{equation*}
    H^p_{(2)}(M):=\ker(\overline\Delta^p) = \im(\pr_p).
  \end{equation*}
\end{definition}

It is an easy observation that, given a projection $P\colon V\to V$ on
a finite dimensional vector space $V$, then $  \dim(\im(P)) = \tr(P)$.
On the other hand, the trace of an operator with a smooth integral
kernel can be computed by integration over the diagonal.

We want to use these ideas to define a useful dimension for
$\ker(\Delta_p)$. Note, however, that the Laplacian
$\overline\Delta_p$ is defined using the Riemannian metric on
$\overline M$. It follows that it commutes with the induced action of
$\Gamma$ on $L^2\Omega^p(\overline M)$. Consequently,
$\ker(\overline\Delta_p)$ is $\Gamma$-invariant, and
\begin{equation*}
  \tr_{\overline x}\pr_p(\overline x,\overline x) = \tr_{g\overline
    x}\pr_p(g\overline x,g\overline x)\qquad\forall \overline x\in\overline 
  M, g\in\Gamma.
\end{equation*}
Observe that $\pr_p(\overline x,\overline x)\in\End(\Lambda^p T^*_xM)$ 
is an endomorphism of a finite dimensional vector space for each
$\overline x\in\overline M$, and $\tr_{\overline x}$ is the usual
trace of such endomorphisms. If $\overline M$ is not
compact (i.e.~if $\Gamma$ is infinite), it follows that
\begin{equation*}
  \int_{\overline M}\tr_{\overline x}\pr_p(\overline x,\overline x)\;d\mu(x)
\end{equation*}
does not converge, and indeed, in general $\ker(\overline\Delta_p)$ is 
not a finite dimensional $\complexs$-vector space.

On the other hand, because of the $\Gamma$-invariance, the function
\begin{equation*}
\overline x\mapsto \tr_{\overline x}\pr_p(\overline x,\overline x)
\end{equation*}
``contains the same information many times'', and it doesn't make
sense to try to compute the integral over all of $\overline M$. We
therefore adopt the notion of ``dimension per (unit) volume''.

More concretely, because of $\Gamma$-invariance, the function
$\overline x\mapsto \tr_{\overline x}\pr_p(\overline x,\overline x)$
descents to a smooth function on the quotient $\overline
M/\Gamma=M$. We now define the $L^2$-Betti numbers
\begin{equation*}
  b^p_{(2)}(\overline M,\Gamma) := \dim_\Gamma\ker(\overline\Delta_p)
  :=\int_M \tr_x \pr_p(x,x)\;d\mu(x) \in [0,\infty).
\end{equation*}
This number is a non-negative real number. However, a priori no
further restrictions appear for these values.

Extensions of all these definitions to manifolds with boundary are
possible, compare e.g.~\cite{Schick(1996)}.

  The definition given here is the original definition of $L^2$-Betti numbers as
  given by Atiyah. Of course, the same construction can be applied to
  any elliptic differential operator $D$ on $M$. If $D$ is such an
  elliptic differential operator, and $D^*$ its formal adjoint, Atiyah 
  defined in this way the $\Gamma$-index of the lift $\overline D$ of
  $D$ to $\overline M$ by
  \begin{equation*}
    \ind_\Gamma(\overline D) :=\dim_\Gamma(\ker\overline D)
    -\dim_\Gamma(\ker\overline{D^*}).
  \end{equation*}
  Atiyah's celebrated $L^2$-index theorem now states
  \begin{theorem}\label{theo:L2index}
    \begin{equation*}
    \ind_\Gamma(\overline D)=\ind(D).
  \end{equation*}
\end{theorem}
  Here, recall that ellipticity of $D$ and compactness of $M$ imply
  that $\ker(D)$ and $\ker(D^*)$ are finite dimensional
  $\complexs$-vector spaces, and 
  \begin{equation*}
\ind(D)=\dim_\complexs(\ker
  D)-\dim_\complexs(\ker D^*).
\end{equation*}
 It should be observed that it is  far
  from true in general that $\dim_\Gamma(\ker \overline
  D)=\dim_\complexs(\ker D)$. In particular, the $L^2$-Betti numbers
  and the ordinary Betti numbers usually are quite different from each 
  other. However, Atiyah's $L^2$-index theorem has the following
  consequence for the $L^2$-Betti numbers:
  \begin{equation}\label{eq:Euler_char_L2}
    \chi(M) = \sum_{p=0}^{\dim M} (-1)^p b^p_{(2)}(\overline M,\Gamma).
  \end{equation}

An extension of Atiyah's $L^2$-index theorem to manifolds with
boundary can be found in \cite{Schick(2001)}, which provides one way
to prove a
corresponding result for the Euler characteristic of manifolds with
boundary.

\begin{example}
  Assume $\Gamma$ is finite. Then $\overline M$ itself is a compact
  manifold and, by the above considerations, we get for the ordinary
  Betti numbers of $\overline M$:
  \begin{equation*}
    b^p(\overline M) = \int_{\overline M} \tr_{\overline x}
    \pr_p(\overline x,\overline x)\;d\mu(\overline x).
  \end{equation*}
  Because of $\Gamma$-invariance,
  \begin{equation*}
    \int_{\overline M} \tr_{\overline x}\pr_p(\overline x,\overline x) \;
    d\mu(\overline x) = \abs{\Gamma}\cdot \int_M \tr_x\pr_p(x,x)\;d\mu(x),
  \end{equation*}
  in other words,
  \begin{equation*}
    b^p_{(2)}(\overline M,\Gamma) =\frac{b^p(\overline M)}{\abs{\Gamma}}
  \end{equation*}
\end{example}

\begin{example}
  Let $p=0$, and assume that $\overline M$ is connected. Integration
  by parts shows that $f\in L^2\Omega^0(\overline M)$ belongs to
  $\ker(\overline\Delta_0)$ if and only if $f$ is constant (recall
  that $L^2\Omega^0(\overline M)=L^2(\overline M)$ is the space of
  $L^2$-functions on $\overline M$). If $\vol(\overline M)=\infty$,
  or equivalently $\abs{\Gamma}=\infty$ then an $L^2$-function $f$ is
  constant if and only if it is zero, i.e.~$\ker(\overline
  \Delta_0)=0$. Therefore, $\pr_p(x,y)=0$ for all $x,y\in\overline M$
  and $b^0_{(2)}(\overline M,\Gamma)=0$.

  Note that the $0$-th ordinary Betti number never vanishes. Many of
  the applications of $L^2$-Betti numbers rely on such vanishing
  results, which don't hold for ordinary Betti numbers.
\end{example}

\begin{theorem}\label{theo:Poincare_duality}
  Assume $\overline M$ is orientable. Then, the Hodge-$*$ operator is
  defined and intertwines $p$-forms and $(\dim M-p)$-forms on
  $\overline M$. Since this is an isometry which commutes with the
  Laplace operators, it induces an isometry between
  $H^p_{(2)}(\overline M)$ and $H^{\dim M-p}_{(2)}(\overline
  M)$. Moreover, this isometry is compatible with the action of
  $\Gamma$ and, in particular extends to the integral kernel of
  $\pr_p$. As a consequence, we have Poincar\'e duality for
  $L^2$-Betti numbers:
  \begin{equation*}
    b^p_{(2)}(\overline M,\Gamma) = b^{\dim M-p}_{(2)}(\overline
    M,\Gamma) .
  \end{equation*}
\end{theorem}

\subsection{The conjectures of Hopf and Singer}

\begin{example}\label{ex:L2_Betti_of_symmetric_spaces}
  In general, it will be almost impossible to compute the $L^2$-Betti
  numbers using the Riemannian metric and the integral kernel of
  $\pr_p$. For very nice metrics, however, this is can be done, in
  particular if $(\overline M,\overline g)$ is a symmetric space. One
  obtains e.g.
  \begin{enumerate}
  \item If $M=T^n$ is a flat torus, $\overline M=\reals^n$ is flat
    Euclidean space, then
    \begin{equation*}
      b^p_{(2)}(\reals^n,\integers^n) =0\qquad\forall p\in\integers.
    \end{equation*}
  \item 
    If $(M,g)$ has constant sectional curvature $K=-1$,
    $\Gamma=\pi_1(M)$ and $\overline M=\hyperbolic^m$ is the
    hyperbolic $m$-plane, then
    \begin{equation*}
      b^p_{(2)}(\overline M,\Gamma) =0\qquad\text{if }p \ne m/2,
    \end{equation*}
    and if $m$ is even and $p=m/2$, then
    \begin{equation*}
      b^{m/2}_{(2)}(\overline M,\Gamma)>0.
    \end{equation*}
    In particular, we conclude, using \eqref{eq:Euler_char_L2}, that
    in this situation
    \begin{equation*}
      (-1)^{m/2}\chi(M) >0.
    \end{equation*}
  \item If, more generally, $( M, g)$ is a
    connected, locally symmetric space with strictly negative
    sectional curvature, $\overline M$ is its universal covering (a
    symmetric space) and $\Gamma=\pi_1(M)$, then
     \begin{equation*}
      b^p_{(2)}(\overline M,\Gamma) =0 \qquad\text{if } p\ne \dim M/2,
    \end{equation*}
    and if $\dim M/2$ is an integer, then
    \begin{equation*}
      b^{\dim M/2}_{(2)}(\overline M,\Gamma) >0.
    \end{equation*}
    In particular, if $\dim M$ is even we have again
    \begin{equation*}
      (-1)^{\dim M/2}\chi(M)>0.
    \end{equation*}
  \item 
   If $(M,g)$ is a connected locally symmetric space with
   non-positive, but not strictly negative,
   sectional curvature, then
   \begin{equation*}
     b^{p}_{(2)}(\overline M,\Gamma) =0\qquad\forall p\in\integers,
   \end{equation*}
   with $\overline M$ and $\Gamma$ as above.
   In particular $\chi(M) =0$.
  \end{enumerate}
\end{example}
\begin{proof}
  These calculations are carried out in \cite{Borel(1985)}, another
  account can be found in \cite{Olbricht(2000)}, using the
  ``representation theory of symmetric spaces''. A more
  geometric proof of the hyperbolic case (i.e.~constant curvature
  $-1$) is given in \cite{Dodziuk(1979)}.
\end{proof}

Given any compact Riemannian manifold without boundary, we can compute the Euler
characteristic using the Pfaffian and the Gauss-Bonnet formula in
higher dimensions. This gives rise to another argument for the inequality
\begin{equation*}
  (-1)^{\dim M/2} \chi(M) >0,
\end{equation*}
if $M$ is an even dimensional manifold with constant negative
sectional curvature, and this has been known for a long time. It has lead
to the following conjecture, which is attributed to Hopf.

\begin{conjecture}\label{conj:Hopf}
  Assume $(M,g)$ is a compact Riemannian $2n$-dimensional manifold
  without boundary, and with strictly negative sectional
  curvature. Then
  \begin{equation*}
    (-1)^n \chi(M) >0.
  \end{equation*}
  If the sectional curvature is non-positive, then
  \begin{equation*}
    (-1)^n\chi(M) \ge 0.
  \end{equation*}
\end{conjecture}

\begin{remark}
  The flat torus, or the product of any negatively curved manifold with
  a flat torus, shows that $\chi(M)=0$ is possible if $M$ is a
  manifold with non-positive sectional curvature.
\end{remark}

In view of Atiyah's formula \eqref{eq:Euler_char_L2}, and because of
the calculations of Jozef Dodziuk and the others in Example
\ref{ex:L2_Betti_of_symmetric_spaces}, Singer proposed to use the
$L^2$-Betti numbers to prove the Hopf conjecture \ref{conj:Hopf}. More 
precisely, he made the following stronger conjecture.

\begin{conjecture}\label{conj:Singer}
  If $(M,g)$ is a compact Riemannian manifold without boundary and
  with non-positive sectional curvature, then
  \begin{equation*}
    b^p_{(2)}(\tilde M,\pi_1(M)) = 0,\qquad\text{if }p\ne \dim M/2.
  \end{equation*}
  If $\dim M=2n$ is even and the sectional curvature is strictly
  negative, then
  \begin{equation*}
    b^{\dim M/2}_{(2)}(\tilde M,\pi_1(M)) >0.
  \end{equation*}
\end{conjecture}
 
Because of \eqref{eq:Euler_char_L2}, the Singer conjecture
\ref{conj:Singer} implies the Hopf conjecture \ref{conj:Hopf}. 

Using
estimates for the Laplacians and their spectrum, Ballmann and
Br\"u\-ning \cite{Ballmann-Bruening(2000)} prove (improving earlier
results of Donnelly
and Xavier \cite{Donnelly-Xavier(1984)} and Jost and Xin
\cite{Jost-Xin(2000)}) part of the Singer conjecture
for very negative sectional curvatures:
\begin{theorem}\label{theo:Donnelly-Xavier}
  If $(M,g)$ is a closed Riemannian manifold of even dimension $2n$,
  such that its sectional
  curvature $K$ satisfies $-1\le K\le -(a_n)^2$, with $1\ge a_n> 1-1/n$, then
  \begin{equation*}
    b^p_{(2)}(\tilde M,\pi_1(M)) =0;\qquad\text{if }p\ne \dim M/2.
  \end{equation*}
  Therefore
  \begin{equation*}
   (-1)^{\dim M/2}\chi( M) \ge 0.
  \end{equation*}
\end{theorem}
\begin{proof}
  This result is not explicitly stated in
  \cite{Ballmann-Bruening(2000)}. However, it follows from their
  \cite[Theorem 5.3]{Ballmann-Bruening(2000)} in the same way in which 
  \cite[Theorem 3.2]{Donnelly-Xavier(1984)} of Donnelly and Xavier
  follows from the corresponding \cite[Theorem
  2.2]{Donnelly-Xavier(1984)}. 
\end{proof}

\begin{remark}
  Classical estimates of the Gauss-Bonnet integrand imply that
  $(-1)^{\dim M/2}\chi(M) >0$ if $-1\le K\le b_n$ with $1\ge b_n > 1-
  3/(\dim M+1)$.

  This gives strict positivity, but the curvature bound of Theorem
  \ref{theo:Donnelly-Xavier} is weaker.
\end{remark}

\subsection{Hodge decomposition}

In the classical situation, (de Rham) cohomology is of course not defined as the
kernel of the Laplacian, as is suggested at the beginning of Section
\ref{sec:analytic-l2-betti}, but as the cohomology of the de Rham
cochain complex. Only afterwards, the Hodge de Rham theorem shows that 
these de Rham cohomology groups are canonically isomorphic to the
space of harmonic forms.

The picture is parallel for $L^2$-cohomology. Whenever $(\overline M,g)$ is a
complete Riemannian manifold, we can define the
$L^2$-de Rham complex
\begin{equation*}
  \to L^2\Omega^{p-1}(\overline M)\xrightarrow{d} L^2\Omega^p(\overline M)\xrightarrow{d}
  L^2\Omega^{p+1}(\overline M)\to 
\end{equation*}
where $d$ is the exterior differential considered as an unbounded
operator on the Hilbert space $L^2\Omega^p(\overline M)$. (The fact that $(\overline M,g)$
is complete implies that $d$ has a unique self adjoint extension.
Usually, we work with this self adjoint extension instead of $d$
itself).

We then define the $L^2$-cohomology as 
\begin{equation*}
  H^p_{(2)}(\overline M):=\ker(d)/\overline{\im(d)}.
\end{equation*}
Observe that we divide through the \emph{closure} of the image of
$d$. This way, we stay in the category of Hilbert spaces. Sometimes,
the $L^2$-cohomology groups obtained this way are called the
\emph{reduced} $L^2$-cohomology groups.

We have to check that this definition coincides with the one given in
Definition \ref{def:L2cohomology}. Now, if $(\overline M,g)$ is a
complete Riemannian manifold,
then we have the following Hodge decomposition:
\begin{equation}\label{eq:Hodgedecomp}
  L^2\Omega^p(M) = \ker(\overline\Delta_p) \oplus \overline{\im d}
  \oplus\overline{\im d^*},
\end{equation}
where $d^*$ is the formal adjoint of $d$, and where the sum is an
orthogonal direct sum. This implies that we also have an orthogonal
decomposition
\begin{equation*}
  \ker(d|_{L^2\Omega^p(\overline M)}) =\ker(\overline \Delta_p) \oplus
  \overline{\im d}.
\end{equation*}
Of course this implies immediately that the inclusion of
$\ker(\overline\Delta_p)$ into $L^2\Omega^p(\overline M)$ induces an
isomorphism between $\ker(\overline\Delta_p)$ and
$\ker(d)/\overline{\im d}$.

\subsection{The Singer conjecture and K\"ahler manifolds}

The use of the Singer conjecture to prove the Hopf conjecture in the
examples presented so far is probably not very impressive. In this
section we will discuss a much more striking result, which makes use
of additional structure, namely the presence of a K\"ahler metric. The 
idea to do this is due to Gromov \cite{Gromov(1991)}.

\begin{definition}
  A Riemannian manifold $(M,g)$ is called a K\"ahler manifold, if the (real)
  tangent bundle $TM$ comes with the structure of a complex vector
  bundle (i.e.~$M$ is an almost complex manifold) with a Hermitian
  metric $h\colon TM\times TM\to \complexs$ (for the given complex
  structure on $TM$) such that the following conditions are satisfied:
  \begin{itemize}
  \item $g$ is the real part of the Hermitian metric $h$.
  \item The $2$-form $\omega$ associated to the Hermitian metric by
    \begin{equation*}
      \omega(v,w) = -\frac{1}{2} Im( h(v,w)),
    \end{equation*}
    the so called \emph{K\"ahler form}, is closed, i.e.~$d\omega=0$.
  \end{itemize}
\end{definition}

\begin{remark}
  Under these conditions, the almost complex structure on $M$ is
  integrable, i.e.~$M$ is a complex manifold. Moreover, $\omega$ is
  non-degenerate, i.e.~$\omega^{\dim M/2}$ is a nowhere vanishing
  multiple of the volume form of $(M,g)$.
\end{remark}

\begin{definition}
  A compact K\"ahler manifold $M$ is called \emph{K\"ahler
    hyperbolic}, if we find a $1$-form $\eta$ on the universal
  covering $\tilde M$ such that $\tilde\omega=d\eta$, where
  $\tilde\omega$ is the pullback of the K\"ahler form $\omega$ to
  $\tilde M$, and such that $\abs{\eta}_\infty:=\sup_{x\in\tilde M}
  \abs{\eta(x)}_x<\infty$.

\end{definition}

\begin{example}
  The following manifolds are K\"ahler hyperbolic:
  \begin{enumerate}
  \item closed K\"ahler manifold which are homotopy equivalent 
    to a Riemannian manifold with negative sectional curvature.
  \item closed K\"ahler manifolds with word-hyperbolic fundamental
    group, provided the second homotopy group vanishes
  \item Complex submanifolds of K\"ahler hyperbolic manifolds, or the
    product of two K\"ahler hyperbolic manifolds.
  \end{enumerate}
\end{example}
\begin{proof}
  Compare \cite[Example 0.3]{Gromov(1991)}.
\end{proof}

From the point of view of the Hopf conjecture and the Singer
conjecture, the first example is the most relevant: manifolds with
negative sectional curvature, which also admit a K\"ahler structure,
are covered by the results of this section.

Gromov proved the Singer conjecture for K\"ahler hyperbolic
manifolds. More precisely, he proved in \cite{Gromov(1991)} the following.
\begin{theorem}\label{theo:Kaehler_hyperbolic}
  Let $M$ be a closed K\"ahler hyperbolic manifold of real dimension
  $2n$. Then
  \begin{equation*}
    \begin{split}
      b_{(2)}^p(\tilde M,\pi_1(M)) =0  &\text{ if }p\ne n\\
      b_{(2)}^n(\tilde M,\pi_1(M)) >0.
  \end{split}
\end{equation*}
In particular, because of \eqref{eq:Euler_char_L2}, $(-1)^n\chi(M) >0$.
\end{theorem}
\begin{proof}
  The proof splits into two rather different parts. On one side, we
  have to show vanishing outside the middle dimension. This is done
  using a Lefschetz theorem, which in particular implies that the
  cup-product with the lift of the K\"ahler form $\tilde \omega$
  induces a bounded injective map
  \begin{equation*}
    L\colon \ker(\tilde \Delta_r)\to \ker(\tilde \Delta_{r+2})\qquad
    \text{for }r<n.
  \end{equation*}
  This is a classical fact from complex geometry in the compact case,
  and it extends rather easily to our non-compact situation.

  But, by assumption, $\tilde\omega=d\eta$ where $\eta$ is an
  $L^\infty$-bounded one form.
  The cup product of a closed form $f$ with a
  boundary $d\eta$ is always a boundary, which shows that (on a
  compact manifold) cup product with a boundary $d\eta$ induces the
  zero map on de Rham cohomology.

  The extra boundedness condition on $\eta$ allows us to deduce the
  same for the $L^2$-de Rham cohomology groups, i.e.~cup product with
  $\tilde \omega$ induces the zero map on $H^p_{(2)}(\tilde M)$ if $M$ 
  is K\"ahler hyperbolic. Since the Lefschetz theorem implies that
  this map is also injective, as long as $p<n$, $H^p_{(2)}(\tilde
  M)=0$ for $p<n$. Because of the Poincar\'e duality theorem
  \ref{theo:Poincare_duality}, the same holds for
  $p>n$, thus establishing the vanishing part of the theorem.

  The first step to prove the non-vanishing of the $L^2$-cohomology in 
  the middle dimension is the following result:
  Using similar, but more delicate, arguments as the one above, one shows that
  $\tilde\Delta_r\colon L^2\Omega^r(\tilde M)\to L^2\Omega^r(\tilde 
  M)$  is
  invertible with a bounded inverse when restricted to the orthogonal
  complement of its null space. (If $r\ne n$, this kernel is $0$, but
  we later have to prove that $\tilde \Delta_n$ has a non-trivial null space.)

  The key point now is that one can construct a continuous family of
  twisted $L^2$-de Rham complexes, indexed by $\lambda\in\reals$, such that
  for $\lambda=0$ we obtain the original $L^2$-de Rham complex we are
  interested in. An extension of Atiyah's $L^2$-index theorem
  \ref{theo:L2index} (a proof can be found e.g.~in \cite[Theorem
  3.6]{Mathai(1999)}) implies that the $L^2$-Euler characteristics of these
  twisted complexes are a non-trivial polynomial $p(\lambda)$. Here,
  the twisted $L^2$-Euler characteristic is defined as the alternating 
  sum of the
  $L^2$-dimension of $\ker(\Delta_p(\lambda))$, where the twisted
  Laplacian $\Delta_p(\lambda)$ is obtained from the twisted de Rham
  complex. Of course, $\Delta_p(0)=\tilde \Delta_0$, and $p(0)=\chi(M)$.

  On the other hand, if $\ker(\tilde\Delta_n)$ was zero, then all the
  operators $\overline\Delta_r$ would be invertible, with bounded
  inverse. Because of the properties of the perturbation, the twisted
  Laplacian $\Delta_r(\lambda)$ would also be invertible for $\lambda$ 
  sufficiently close to zero, which would imply that $p(\lambda)=0$
  for $\lambda$ sufficiently close to $0$. This is a contradiction to
  the result that $p(\lambda)$ is a non-trivial polynomial.
\end{proof}

\begin{remark}
  The vanishing part of Theorem \ref{theo:Kaehler_hyperbolic}
  definitely is the easier part. To obtain the corresponding
  statement, the assumptions on the K\"ahler form can be weakened a
  little bit, namely, it suffices that the lift  $\tilde\omega$
  satisfies $\tilde\omega=d\eta$
  for a one form $\eta$ which has at most linear growth (such
  manifolds are called \emph{K\"ahler non-elliptic}). This is carried
  out independently in \cite{Jost-Zuo(2000)} and
  \cite{Cao-Xavier(2001)}. The most important example of K\"ahler
  non-elliptic manifolds are closed Riemannian manifolds with
  non-positive sectional curvature which admit a K\"ahler
  structure. In particular, for such a manifold the assertions of the
  Singer conjecture and of the Hopf conjecture are true.
\end{remark}

\section{Combinatorial $L^2$-Betti numbers}

So far, we have defined $L^2$-cohomology and $L^2$-Betti numbers only
for coverings of Riemannian manifolds, and we have not even shown that 
they do not depend on the chosen Riemannian metric. Here, we will
extend the definition to coverings of arbitrary CW-complexes. Moreover, 
we will show that $L^2$-cohomology is an equivariant homotopy
invariant.

From now on, therefore, assume that $X$ is a compact CW-complex, and
that $\overline X$ is a regular covering of $X$, with deck
transformation group $\Gamma$ (in particular, $\Gamma$ acts freely on
$\overline X$, and $X=\overline X/\Gamma$).
We use the induced structure of a CW-complex on $\overline X$, and
$\Gamma$ acts by cellular homeomorphisms.

Then, the cellular chain complex of $\overline X$ is a chain complex
of finitely generated free $\integers\Gamma$-modules, with
$\integers\Gamma$-basis
given by lifts of cells of $X$ . This way, the basis is unique up to
permutation and multiplication
with $\pm g\in\integers \Gamma$ for $g\in\Gamma$.

The cellular $L^2$-cochain complex is defined by
\begin{equation*}
  C^*_{(2)}(\overline
  X,\Gamma):=\Hom_{\integers\Gamma}(C^{cell}_*(\overline
  X),l^2(\Gamma)).
\end{equation*}

We assume that $\Gamma$ acts on $X$ from the right, and therefore
consider homomorphisms of right $\integers\Gamma$-modules. $\Gamma$
still acts on $C^*_{(2)}(\overline \Gamma)$ by isometries, with
\begin{equation*}
(fg)(x):=g^{-1}\cdot(f(x))\text{ for }\in\Gamma,\;f\in C^*_{(2)}(\overline
X)\text{ and }x\in C^{cell}_*(\overline X).
\end{equation*}
The choice of a cellular
$\integers\Gamma$-basis for $C^{cell}_*(\overline X)$ identifies
$C^p_{(2)}(\overline X)$ with $(l^2\Gamma)^n$ ($n$ being the number of 
$p$-cells in $X$), which induces in particular the structure of a
Hilbert space on the cochain groups (which does not depend on the
particular cellular basis).

The $L^2$-cochain maps 
\begin{equation*}
d_p\colon C^p_{(2)}(\overline X,\Gamma)\to
C^{p+1}_{(2)}(\overline X,\Gamma)
\end{equation*}
are bounded equivariant linear
maps.
Let $d_p^*$ be the adjoint operator, and set 
\begin{equation*}
\Delta_p:=d_p^*d_p+ d_{p-1}d_{p-1}^*.
\end{equation*}
This is the \emph{cellular Laplacian}, a bounded
self adjoint
equivariant operator on $C^p_{(2)}(\overline X)$.

\begin{remark}
  We can also
  define the cellular $L^2$-chain complex by
  \begin{equation*}
    C_*^{(2)}(\overline X):= C^{cell}_*(\overline
    X)\tensor_{\integers\Gamma} l^2(\Gamma).
  \end{equation*}
  The self duality of the Hilbert space $l^2(\Gamma)$ induces a
  duality between $C_*^{(2)}$ and $C^*_{(2)}$ which gives a
  chain isometry between the $L^2$-cochain and the $L^2$-chain
  complex. In particular, in all the definitions
  we are going to make, there will be a canonical isomorphism between
  the cohomological and homological version. Because for manifolds de
  Rham cohomology seems to be more natural to use, we will stick to
  the cohomological version throughout.
\end{remark}

\subsection{Hilbert modules}

\begin{definition}
  A (finitely generated) Hilbert $\NeumannN\Gamma$-module is a Hilbert
  space $V$ with a right $\Gamma$-action which admits an equivariant
  isometric embedding into $l^2(\Gamma)^n$ for some $n$.
\end{definition}

\begin{remark}
  To explain this notation, we remark that an isometric action of
  $\Gamma$ on a (complex) Hilbert space by linearity extends to an
  ``action'' of the integral group ring $\integers[\Gamma]$ and the
  complex group ring $\complexs[\Gamma]$ (recall that, given a ring
  $R$, the group ring $R\Gamma$ is defined to consist of finite formal
  linear combination $\sum_{g\in G} r_g g$ with $r_g\in R$, with
  component-wise addition, and multiplication defined by
  $(r_gg)(r_hh)=(r_gr_h)(gh)$). We can embed $\complexs \Gamma$ into a
  certain completion, the reduced $C^*$-algebra $C^*_r\Gamma$.
\begin{definition}
  The \emph{group von Neumann algebra} $\NeumannN\Gamma$ is an even
  bigger completion of $\complexs\Gamma$. It can be defined to
  consists of those bounded operators on $l^2\Gamma$ which commute
  with the \emph{right} action of $\Gamma$ on $l^2\Gamma$,
  i.e.
  \begin{equation*}
  \NeumannN\Gamma:=\boundedops(l^2\Gamma)^\Gamma.
\end{equation*}
The right
  action is given by 
  \begin{equation*}
(\sum_{g\in\Gamma}\lambda_g g)\cdot v :=
  \sum_{g\in\Gamma}\lambda_g(gv)\text{ for }v\in G\text{ and }
  \sum_{g\in\Gamma}\lambda_g g\in l^2\Gamma.
\end{equation*}
Then $\NeumannN\Gamma$
  is a ring which acts on the left on $l^2\Gamma$, and
  $\complexs\Gamma$ (actually $C^*_r\Gamma$) is contained in
  $\NeumannN\Gamma$.
  
  Equivalently, one can define $\NeumannN\Gamma$ to be the closure of
  $\complexs\Gamma$ (with the left action) in $\boundedops\Gamma$ with
  respect to the weak topology. This is a consequence of von Neumann's
  bicommutant theorem.
\end{definition}

The action of $\Gamma$ on a Hilbert $\NeumannN\Gamma$-module $V$
extends to $\complexs\Gamma$ and then to $\NeumannN\Gamma$, making $V$
indeed a module over $\NeumannN\Gamma$. Observe, however, that we
don't get arbitrary algebraic modules, but modules with a topology and
of a rather special kind. This additional $\NeumannN\Gamma$-module
structure is underlying many of the definitions and proofs in
$L^2$-cohomology. However, in the sequel we will omit explicitly using
this, and instead work with the (for our purposes equivalent) unitary
action of $\Gamma$ and existence of the embedding into
$(l^2\Gamma)^n$.
\end{remark}

Given a Hilbert $\NeumannN\Gamma$-module $V$, let $\pr\colon l^2(\Gamma)^n\to
l^2(\Gamma)^n$ be the orthogonal projection onto the image of any such 
embedding. We define the $\Gamma$-dimension of $V$ by
\begin{equation}\label{eq:gamma_trace}
  \dim_\Gamma(V):=\tr_\Gamma(\pr):= \sum_{i=1}^n
  \innerprod{\pr(e_i),e_i}_{l^2(\Gamma)^n} .
\end{equation}
Here, $e_i=(0,\dots,\delta_1,\dots,0)$ is the standard basis vector of 
$l^2(\Gamma)^n$ with $i$-th entry being the characteristic function of 
the unit of $\Gamma$, and all other entries being zero.

Observe that $\Gamma$-invariance of $V$ implies that $\pr$ is
$\Gamma$-equivariant, i.e.
\begin{equation*}
\innerprod{\pr(e_ig),e_i
  g}=\innerprod{\pr(e_i),e_i}
\end{equation*}
for all $g\in\Gamma$. If we would like
to compute the $\complexs$-dimension of $V$, and therefore take the
ordinary trace of $\pr$ (as endomorphism of
$\complexs$-vector spaces) we would have to sum over
$\innerprod{\pr(e_i g),e_i g}$ for all $g\in\Gamma$. Of course, in
general this
doesn't make sense since $\pr$ is not of trace class. As in the
Definition \ref{def:L2cohomology}, we pick the relevant part of this
trace in \eqref{eq:gamma_trace}, summing over a ``fundamental domain'' 
for the $\Gamma$-action on $l^2(\Gamma)^n$.

It is not hard to check that the above definition is independent of
the choice of the embedding of $V$ into $l^2(\Gamma)^n$.

\begin{example}
  If $\Gamma$ is finite, then every finitely generated Hilbert
  $\NeumannN\Gamma$-module $V$ is a finite dimensional vector space over
  $\complexs$, and 
  \begin{equation*}
    \dim_\Gamma(V) = \frac{1}{\abs{\Gamma}}\dim_\complexs(V).
  \end{equation*}
\end{example}

\begin{example}\label{ex:Hilbert_Z_modules}
  A more interesting example is given by free abelian groups. Assume
  that $\Gamma=\integers$. Then Fourier transform provides an
  isometric isomorphism between $l^2(\Gamma)$ and $L^2(S^1)$. Under
  this isomorphism, the subspace $\complexs[\Gamma]$ corresponds to
  the space of trigonometric polynomials in $L^2(S^1)$, which act by
  pointwise multiplication. The reduced $C^*$-algebra $C^*_r\integers$ 
  becomes $C(S^1)$, and the von Neumann algebra $\NeumannN\integers$
  becomes $L^\infty(S^1)$, also acting by pointwise multiplication. A projection $P$
  of $L^2(\Gamma)$ which commutes with all these trigonometric
  polynomials is itself given by multiplication with a measurable
  function $f$, and being a projection translates to the fact that $f$
  only takes the values $0$ and $1$ (up to a set of measure zero).

  The image of $P$ is the set of functions in $L^2(S^1)$ which vanish
  on the zero set of $f$. The
  $\Gamma$-trace of $P$ is the constant term in the Fourier expansion
  of $f$, which can be computed by integration over $S^1$, i.e.
  \begin{equation*}
    \tr_\Gamma(P) = \int_{S^1} f = \vol(\supp(f)),
  \end{equation*}
  which here is of course just the volume of the support of $f$,
  i.e.~the set of all $x\in S^1$ with $f(x)=1$
  We use the standard measure on $S^1$, normalized in such a way 
  that $\vol(S^1)=1$.

  It should be observed that one can obtain any real number between
  $0$ and $1$ in this way.
\end{example}

The $\Gamma$-dimension has the following useful properties, which in
particular justify the term ``dimension''.
\begin{proposition}\label{prop:L2_dim_properties}
  Let $U,V,W$ be finitely generated Hilbert $\NeumannN\Gamma$-modules.
  \begin{enumerate}
  \item Faithfulness: $\dim_\Gamma(U)=0$ if and only if $U=0$.
  \item Additivity: If we have a weakly exact sequence of Hilbert
    $\NeumannN\Gamma$-modules 
    \begin{equation*}
0\to 
    U\to W\to V\to 0
  \end{equation*}
then
\begin{equation*}
\dim_\Gamma(W)=\dim_\Gamma(U)+\dim_\Gamma(V).
\end{equation*}
    Weakly exact means that the kernel of the outgoing map coincides
    with the \emph{closure} of the image of the incoming map,
    i.e.
    \begin{equation*}
\cdots\xrightarrow{\phi_1} X\xrightarrow{\phi_2}\cdots
\end{equation*}
is
    weakly exact at $X$ if and only if $\ker(\phi_2)=\overline{\im(\phi_1)}$.
  \item Monotonicity: If $U\subset V$ then $\dim_\Gamma(U)\le
    \dim_\Gamma(V)$, and $\dim_\Gamma(U)=\dim_\Gamma(V)$ if and
    only if $U=V$.
  \item Normalization: $\dim_\Gamma(l^2(\Gamma)) =1$.
  \item If $H$ is a subgroup of finite index $d$
  in $\Gamma$, then every finitely generated Hilbert $\NeumannN\Gamma$-module $V$
  becomes by restriction of the action a finitely generated Hilbert
  $\NeumannN H$-module. Then
  \begin{equation*}
    \dim_H(V) = d\cdot \dim_\Gamma(V).
  \end{equation*}
  (Note that $\Gamma$ finite and $H$ trivial is a special case of this
  situation.)
  \end{enumerate}
\end{proposition}

\subsection{Cellular $L^2$-cohomology}

\begin{definition}
  We define the cellular $L^2$-cohomology by
  \begin{equation*}
    H^p_{(2)}(\overline X,\Gamma):= \ker(d^p)/\overline
    \im(d^{p-1}).
  \end{equation*}

  We have a Hodge decomposition
  \begin{equation*}
    C^p_{(2)}(\overline M,\Gamma) = \ker(\Delta_p)\oplus \overline{\im d}
    \oplus \overline{\im d^*}.
  \end{equation*}
  This is similar to Hodge decomposition for differential forms on
  Riemannian manifolds, but, since all operators involved here are
  bounded, is a much more elementary result. From this it follows that 
  we have an isometric $\Gamma$-isomorphism
  \begin{equation*}
    H^p_{(2)}(\overline X,\Gamma)  \iso\ker(\Delta_p)
  \end{equation*}

  We define the $L^2$-Betti numbers
  \begin{equation*}
    b^p_{(2)}(\overline X,\Gamma):= \dim_\Gamma(\ker(\Delta_p)).
  \end{equation*}
\end{definition}

Observe again the important fact that we divide by the closure of the
image of the differential, such as to remain in the category of
Hilbert spaces. This is the decisive difference to the equivariant
cohomology with values in the $\integers\Gamma$-module
$l^2\Gamma$. However, in \cite{Lueck(1998a),Lueck(1998b)}, L\"uck
generalized the concept of $L^2$-Betti numbers from normal coverings
of finite CW-complexes to 
arbitrary spaces with group action, using the usual twisted cohomology 
(with coefficients the group von Neumann algebra
$\NeumannN\Gamma$). The starting point, however, is also in this
treatment the theory
of Hilbert $\NeumannN\Gamma$-modules.

\subsubsection{Matrices over the group ring}
\label{sec:matrices-over-group}

Given a compact CW-complex, we can explicitly compute it's cohomology
using a cellular basis and solving certain systems of linear
equations. A similar approach is possible here (leading to more
complicated, ``non-commutative'', equations in this situation).
 
In the compact case, the choice of an orientation for each cell
identifies $C_p^{cell}(X)$ with $\integers^{c_p}$, where $c_p$ is the
number of $p$-cells of $X$, and this identification is well defined up 
to permutation of the basis, and multiplication of basis elements with 
$\pm 1$. The boundary map, in this representation, is given by
multiplication with an appropriate matrix with integral entries.

To proceed in the $L^2$-case, observe that each cell of the finitely
many cells of $X$
has as inverse image a free $\Gamma$-orbit of cells in $\tilde X$. We
can choose one cell in each orbit. Together with the choice of an
orientation, this identifies the $\integers\Gamma$-module
$C^{cell}_p(\overline X)$ with $(\integers\Gamma)^{c_p}$, and this
identification is unique up to multiplication of the basis elements
with $\pm g$ ($g\in\Gamma$) and permutation. In this realization, the
boundary maps of $C_p^{cell}(\overline X)$ are given by
multiplication with matrices $A_p$ over the integral group ring. 

The chosen $\integers\Gamma$-module isomorphism of
$C^{cell}_p(\overline X)$ with $(\integers\Gamma)^{c_p}$ induces an
isomorphism of $C^*_{(2)}(\overline X,\Gamma)$ with
$(l^2\Gamma)^{c_p}$. An easy calculation shows that the coboundary
maps are given by multiplication with the adjoint matrices
$A_p^*$ (extending the multiplication of elements of $\complexs\Gamma$ 
with elements of $l^2\Gamma$ used before). Here, if
\begin{equation*}
u=\sum_{g\in\Gamma}\lambda_g g\in \complexs\Gamma
\text{ then }u^*:=\sum_{g\in\Gamma} \overline{\lambda_g} g^{-1},
\end{equation*}
and
obviously $u^*\in\integers\Gamma$ if $u\in\integers\Gamma$.
If $A=(A_{ij})\in M(d_1\times d_2,\complexs\Gamma)$, then 
\begin{equation*}
A^*:=
(A^*_{ji})\in M(d_2\times d_1,\complexs\Gamma),
\end{equation*}
and again this
restricts to an operation on matrices over $\integers\Gamma$.

To finish the picture, the combinatorial Laplacian 
$\Delta_p=d_p^*d_p+d_{p-1}d_{p-1}^*$ 
is given by the matrix 
\begin{equation*}
\Delta:= A_p^*A_p+A_{p-1}^*A_{p-1}\in
M(c_p\times c_p,\integers\Gamma).
\end{equation*}
To understand the $L^2$-cohomology
we therefore have to understand the kernel of such matrices, acting
on $(l^2\Gamma)^{c_p}$.

This gives an algebraic way of studying questions about
$L^2$-cohomology ---they translate to questions about matrices over
$\integers\Gamma$.

Actually, if $\Gamma$ is finitely presented (i.e.~has a presentation
with finitely many generators and finitely many relations), given any matrix $A\in
M(d\times d,\integers\Gamma)$, a
standard construction provides us with a compact CW-complex $X$ with
$\pi_1(X)=\Gamma$ such that the kernel of the combinatorial Laplacian
$\Delta_3$ becomes the kernel of $A$, acting on
$l^2(\Gamma)^d$. Therefore, we can also translate questions about
matrices over $\integers\Gamma$ to questions in $L^2$-cohomology.

\subsubsection{Properties of $L^2$-Betti numbers}

\begin{theorem}\label{theo:properties_of_L2Bettinumbers}
$L^2$-cohomology and in particular $L^2$-Betti numbers have the
following basic properties.
Here, let $\overline X$ be a normal covering of a finite CW-complex
    $X$ with covering group $\Gamma$.
  \begin{enumerate}
  \item\label{item:Hodge_de_Rham} Let $(M,g)$ be a closed Riemannian
    manifold, and equip it with the
    CW-structure coming from a smooth triangulation. Let $(\overline
    M,\overline g)$ 
    be a normal covering of $M$ with covering group $\Gamma$. Then
    integration of forms over simplices (the de Rham map) defines a
    $\Gamma$-isomorphism
    \begin{equation*}
      \ker(\overline \Delta_p(\overline g)) \to
      H^p_{(2),cell}(\overline M,\Gamma).
    \end{equation*}
    In particular, the $L^2$-Betti numbers defined using the Riemannian 
    metric and using the $\Gamma$-CW-structure coincide.
  \item\label{item:homotopy_invariance} Let $Y$ be another finite CW-complexes and $f\colon Y\to X$ a
    homotopy equivalence. Let $\overline Y$ be the pullback 
    of $\overline X$ along $f$ (this means that the deck
    transformation group for $\overline Y$ is also $\Gamma$). Then
    \begin{equation*}
      b_{(2)}^p(\overline X,\Gamma)=b_{(2)}^p(\overline
      Y,\Gamma)\qquad\forall p\ge 0.
    \end{equation*}
  \item \label{item:Eulerchar} For the Euler characteristic of $X$, we get
    \begin{equation}\label{eq:Euler_char_for_cW}
      \chi(X) = \sum_{p=0}^\infty (-1)^p b^p_{(2)}(\overline X,\Gamma).
    \end{equation}
  \item \label{item:finite_group}Assume $\Gamma$ is finite. Then
    $\overline X$ is itself a
    finite CW-complex, and its (ordinary) Betti number $b^p(\overline
    X)$ are defined. They satisfy
    \begin{equation*}
      b^p_{(2)}(\overline X,\Gamma) = \frac{1}{\abs{\Gamma}}
      b^p(\overline X).
    \end{equation*}
  \item\label{item:coverings} If $\Gamma$ is infinite, then for the
    zeroth $L^2$-Betti
    number we get
    \begin{equation*}
      b^0_{(2)}(\overline X,\Gamma) =0.
    \end{equation*}
  \item \label{item:multipicative}
  Let $H\subgroup\Gamma$ be a subgroup of finite index $d$, and
    set $X_1:=\overline X/H$. This is a finite $d$-sheeted covering of 
    $X$, and $\overline X$ can be considered to be a normal covering
    of $X_1$ with covering group $H$. Then
    \begin{equation*}
      b^p_{(2)}(\overline X,H) = d\cdot b^p_{(2)}(\overline
      X,\Gamma)\qquad \forall p\ge 0.
    \end{equation*}
  \end{enumerate}
\end{theorem}
Note that the Euler characteristic is multiplicative under finite
covering, in the situation of \ref{item:multipicative} this means that 
$\chi(X_1)=d\cdot \chi(X)$. In view of \ref{item:multipicative} and the
Euler characteristic formula \ref{item:Eulerchar}, the $L^2$-Betti
numbers are the appropriate refinement of the Euler characteristic
which (unlike the ordinary Betti numbers) remain multiplicative under 
finite coverings.

\begin{proof}[Proof of Theorem
  \ref{theo:properties_of_L2Bettinumbers}]
  We only indicate reference and the main points.
  \begin{itemize}
  \item[\ref{item:Hodge_de_Rham}] This was a classical question of Atiyah, proved by Dodziuk in
    \cite{Dodziuk(1977)}.
  \item[\ref{item:homotopy_invariance}] $f$ is covered by a $\Gamma$-homotopy equivalence $\overline
    f\colon \overline Y\to\overline X$. One easily checks that such a
    map induces a map on $L^2$-cohomology, and that two
    $\Gamma$-homotopic maps induce the same map. The claim follows.
  \item[\ref{item:Eulerchar}] The Euler characteristic formula follows
    exactly as in the
    classical situation, using additivity of the $\Gamma$-dimension
    and the normalization
    \begin{equation*}
\dim_\Gamma(l^2\Gamma)=1.
\end{equation*}
  \item[\ref{item:finite_group}] If $\Gamma$ is finite, all the
    Hilbert spaces in question are finite dimensional. Consequently,
    $\im(d)$ is automatically closed, and there is no difference
    between $L^2$-cohomology and ordinary cohomology with complex
    coefficients of $\overline X$. In particular, 
    \begin{equation*}
b^p(\overline
    X)=\dim_\complexs H^p_{(2)}(\overline X,\Gamma)=
    \abs{\Gamma}\cdot \dim_\Gamma(H^p_{(2)}(\overline X,\Gamma))
    =\abs{\Gamma}\cdot b^p_{(2)}(\overline X,\Gamma).
  \end{equation*}
\item[\ref{item:multipicative}] This follows immediately from the
  corresponding formula in Proposition \ref{prop:L2_dim_properties}.
  \end{itemize}
\end{proof}

\section{Approximating $L^2$-Betti numbers}\label{sec:appr-l2-betti}

As mentioned in Section \ref{sec:analytic-l2-betti}, there are almost
no relations between the ordinary Betti numbers of a space $X$ and
the $L^2$-Betti numbers of a covering $\tilde X$ of $X$. However, if
we have a whole sequence of nested coverings $X\leftarrow
X_1\leftarrow X_2\leftarrow\cdots$, ``converging'' to $\tilde X$, in
many cases we can approximate the $L^2$-Betti numbers of $\tilde X$ in 
terms of this sequence. More precisely, let $\tilde X$ be a
$\Gamma$-covering of $X$. Assume that there is a nested sequence of
normal subgroups $\Gamma\supergroup
\Gamma_1\supergroup\Gamma_2\supergroup\cdots$ (each $\Gamma_k$ normal
in $\Gamma$) such that $\bigcap_{k\ge 1} \Gamma_k =\{1\}$. Then
$X_k:=\Gamma_k\backslash \tilde X$ is a normal covering of $X$, with
covering group $\Gamma/\Gamma_k$.

\begin{conjecture}\label{conj:convergence}
  In this situation, 
  \begin{equation}\label{eq:convergence}
    b^p_{(2)}(\tilde X,\Gamma) = \lim_{k\to\infty} b^p_{(2)}(X_k,\Gamma/\Gamma_k).
  \end{equation}
  Observe that convergence is not clear, but part of the statement.
\end{conjecture}

This question was first asked by Gromov if all the groups
$\Gamma/\Gamma_k$ are finite. In this case, the conjecture is true, as 
proved by L\"uck:
\begin{theorem}\label{theo:Lueck_approxi}
  Equation \eqref{eq:convergence} is correct if $\Gamma/\Gamma_k$ is
  finite for each $k\in\naturals$.

  Observe that, in this setting, $X_k$ is a finite covering of
  $X$. Consequently, we can express the $L^2$-Betti numbers in terms
  of ordinary Betti numbers and obtain, in the setting of Conjecture
  \ref{conj:convergence}, 
  \begin{equation*}
    b^p_{(2)}(\tilde X,\Gamma) =\lim_{k\to\infty}
    \frac{b^p(X_k)}{\abs{\Gamma/\Gamma_k}}. 
  \end{equation*}
\end{theorem}

In general, Conjecture \ref{conj:convergence} is still open. However,
in \cite{Schick(1998a)} (with an improvement in
\cite{Dodziuk-Linnell-Mathai-Schick-Yates(2001)}) a quite large class
$\RAgroups$ of groups
is constructed for which the conjecture is true. $\RAgroups$ contains
all amenable and all free 
groups, and is closed under taking subgroups, extensions with
amenable quotients, directed unions, and inverse limits (therefore, it 
contains e.g.~all residually finite groups). 

More precisely
\begin{definition}
  Let $\RAgroups$ be the smallest class of groups which contains the trivial
  group and is closed under the following processes:
  \begin{itemize}
   \item If $H\in\RAgroups$ and $G$ is a generalized amenable
     extension of $H$, then $G\in\RAgroups$.
  \item If $G=\lim_{i\in I}G_i$ is the direct or inverse limit of a
    directed system of groups $G_i\in \RAgroups$,
    then  $G\in\RAgroups$.
  \item If $H\in \RAgroups$ and $U\subgroup H$, then $U\in\RAgroups$.
  \end{itemize}
\end{definition}

Here, the notion of \emph{generalized amenable extension} is defined
as follows:\begin{definition}\label{def:generalized_amenable_extension}
  Assume that $G$ is a finitely generated discrete group with a finite
  symmetric set of generators $S$ (i.e.~$s\in S$ implies $s^{-1}\in
  S$), and let $H$ be an
  arbitrary discrete group. We say that \emph{$G$ is a generalized
    amenable extension of $H$}, if there is a set $X$ with a free $G$-action
  (from the left) and a commuting free $H$-action (from the right),
  such that a sequence of $H$-subsets $X_1\subset X_2\subset
  X_3\subset \cdots \subset X$ exists with $\bigcup_{k\in\naturals} X_k =X$,
  and with $\abs{X_k/H}<\infty$ for every $k\in\naturals$, and such
  that
  \begin{equation*}
   \frac{ \abs{(S\cdot X_k - X_k)/H}}{\abs{X_k/H}}
   \xrightarrow{k\to\infty} 0.
 \end{equation*}
\end{definition}

In \cite{Schick(1998a)} and
\cite{Dodziuk-Linnell-Mathai-Schick-Yates(2001)},
the following theorem is proved.
\begin{theorem}\label{theo:general_approximation}
  Equation \eqref{eq:convergence} is correct if $\Gamma/\Gamma_k$
  belongs to $\RAgroups$ for each $k\in\naturals$.  
\end{theorem}
 In particular, we
obtain the following corollary which we will use later.
\begin{corollary}\label{corol:amenable_approximation}
  Equation \eqref{eq:convergence} is true if $\Gamma/\Gamma_k$ is
  amenable, e.g.~solvable or nilpotent, or virtually solvable, for
  each $k\in\naturals$. Recall that a group $G$ has \emph{virtually}
  a certain property
  $P$, if it contains a subgroup of finite index which has property $P$.
\end{corollary}

\begin{remark}
  There are generalizations of the above approximation results to
  other $L^2$-invariants, in particular to the $L^2$-signature,
  compare \cite{Lueck-Schick(2001)}.
\end{remark}

\section{The Atiyah conjecture}

Fix a discrete group $\Gamma$. The $L^2$-Betti numbers
$b^p_{(2)}(\overline X,\Gamma)$ of a $\Gamma$-covering of a finite
CW-complex $X$ are the $\Gamma$-dimensions of certain Hilbert
$\NeumannN\Gamma$-modules. In Example \ref{ex:Hilbert_Z_modules} we
have seen that a priori arbitrary non-negative real numbers could
occur, even for groups as nice as $\integers$. However, the Euler
characteristic formula \eqref{eq:Euler_char_for_cW} shows that certain
combinations of$L^2$-Betti numbers are always integers. 

The Atiyah conjecture predicts a certain amount of integrality for the 
individual $L^2$-Betti numbers.
\begin{conjecture}\label{conj:Atiyah}
  Fix a discrete group $\Gamma$. Let $Fin^{-1}(\Gamma)$ be the additive
  subgroup of $\rationals$ generated by 
  \begin{equation*}
\{\frac{1}{\abs{F}}\mid
    F\text{ finite subgroup of }\Gamma\}.
  \end{equation*}
  Let $X$ be a finite CW-complex or a compact manifold, $\overline X$ a
  $\Gamma$-covering of $X$.
  \begin{enumerate}
  \item If $\Gamma$ is torsion-free, then
    \begin{equation*}
      b^p_{(2)}(\overline X,\Gamma)\in\integers.
    \end{equation*}
  \item Assume there is a bound on the finite subgroups of $\Gamma$
    (observe that this is equivalent to $Fin^{-1}(\Gamma)$ being a
    discrete subset of $\reals$). Then
    \begin{equation*}
      b^p_{(2)}(\overline X,\Gamma)\in Fin^{-1}(\Gamma).
    \end{equation*}
  \item\label{item:rational_Atiyah} Without any assumption on $\Gamma$,
    \begin{equation*}
      b^p_{(2)}(\overline X,\Gamma)\in \rationals.
    \end{equation*}
  \end{enumerate}
\end{conjecture}

For a while, also the following conjecture was around:
\begin{conjecture}
  Without any assumption on $\Gamma$, $b^p_{(2)}(\overline X,\Gamma)\in
  Fin^{-1}(\Gamma)$. 
\end{conjecture}
This last conjecture is singled out here because it is wrong. In
\cite{Grigorchuk-Linnell-Schick-Zuk(2000)}, a smooth $7$-dimensional
Riemannian manifold $M$ is constructed such that every finite subgroup
of $\pi_1(M)$ is an elementary abelian $2$-group, but
$b^3_{(2)}(\tilde M,\pi_1(M))=\frac{1}{3}$. This example is based on
the explicit calculation of the eigenspaces and their $L^2$-dimensions
of a certain operator in
\cite{Grigurchuk-Zuk(1999)}, using in particular the methods of the
proof of Theorem \ref{theo:Lueck_approxi}. A more direct and slightly
more general computation for such eigenspaces is carried out in
\cite{Dicks-Schick(2001)}.

It should be remarked that none of the above conjectures were
formulated by Atiyah as stated here, although he makes some remarks which
show that he was interested in the question of the possible values of
$L^2$-Betti numbers.

Statement \ref{item:rational_Atiyah} of Conjecture \ref{conj:Atiyah},
which is the oldest version of the Atiyah conjecture, is also quite
unlikely to hold in general. In \cite{Dicks-Schick(2001)}, for each
$r,s\in\naturals$ with $r,s\ge 2$, a manifold $M_{r,s}$ is constructed 
such that
\begin{equation*}
  b^3_{(2)}(\tilde M,\pi_1(M)) = \alpha_{r,s}:=(r-1)^2(s-1)^2\cdot
  \sum_{n=2}^\infty \frac{\phi(n)}{(r^n-1)(s^n-1)},
\end{equation*}
where $\phi(n)$ is Euler's phi-function, i.e.~the number of primitive
$n$-th roots of unity. At the moment, it is unknown whether any of the 
numbers $\alpha_{r,s}$ is irrational. However, the use of computer
algebra shows e.g.~that, if $\alpha_{2,2}$ is rational, both the numerator 
and denominator exceed $10^{100}$. It seems reasonable to assert that
$\alpha_{2,2}$ is not obviously rational.

\subsection{Combinatorial reformulation of the Atiyah Conjecture}

The following assertion is equivalent to Conjecture \ref{conj:Atiyah}.
\begin{conjecture}\label{conj:comb_Atiyah}
  Let $\Gamma$ be a discrete group, and assume that $A\in M(d\times
  d,\integers\Gamma)$. Consider $A$ to be a bounded operator on
  $l^2(\Gamma)^d$, as in Section \ref{sec:matrices-over-group}.
  \begin{enumerate}
  \item If $\Gamma$ is torsion-free, then
    $\dim_\Gamma(\ker(A))\in\integers$.
  \item If $Fin^{-1}(\Gamma)$ is a discrete subset of $\reals$, then
    \begin{equation*}
      \dim_\Gamma(\ker(A))\in Fin^{-1}(\Gamma).
    \end{equation*}
    Without any assumption on $\Gamma$, $\dim_\Gamma(\ker(A))\in\rationals$.
  \end{enumerate}
\end{conjecture}

\begin{remark}
  It is equivalent to require the assertions of Conjecture
  \ref{conj:comb_Atiyah} for all matrices over $\integers\Gamma$, or
  for all square matrices, or for all self-adjoint matrices (these are
  automatically square matrices), or for all matrices of the form
  $A=B^*B$ (these are automatically self-adjoint). This is the case
  since $\ker(A)=\ker(A^*A)$. Moreover, we have the weakly exact
  sequence
  \begin{equation*}
    0\to \ker(A)\into l^2(\Gamma)^d\xrightarrow{A} \overline{\im A}
    \to 0. 
  \end{equation*}
  Because of additivity and normalization of the $\Gamma$-dimension
  (Proposition \ref{prop:L2_dim_properties}), we could replace the
  kernel of $A$ by the closure of the image, throughout.
\end{remark}

The equivalence of Conjecture \ref{conj:Atiyah} and Conjecture
\ref{conj:comb_Atiyah} follows immediately from the principle
described at the end of Section \ref{sec:matrices-over-group}.

\subsection{Atiyah conjecture and non-commutative algebraic
  geometry --- Generalizations}

As an illustration, we now want to study the Atiyah conjecture for the 
group $\Gamma=\integers$, which we understand particularly well
because of Example \ref{ex:Hilbert_Z_modules}.We look at the algebraic 
reformulation. For simplicity, assume
first that $d=1$. Then $A\in\integers[\integers]=\integers[z,z^{-1}]$
is a Laurent polynomial with $\integers$-coefficients. Under Fourier
transform we get the commutative diagram
\begin{equation*}
  \begin{CD}
    l^2(\integers) @>{A}>> l^2(\integers)\\
    @VV{\iso}V @VV{\iso}V\\
    L^2(S^1) @>{A(z)}>> L^2(S^1),
  \end{CD}
\end{equation*}
i.e.~the action of $A$ translates to multiplication with the function
$A(z)$, ($z\in S^1\subset \complexs$). Now, 
\begin{equation*}
\dim_\integers(\ker(A))=
\mu(\{z\in S^1\mid A(z)=0\})
\end{equation*}
is the volume of the set of zeros of
the Laurant polynomial $A(z)$ on $S^1$. But the Laurant polynomial
$A(z)$ has, if it
is not identically zero, only finitely many zero. Therefore
\begin{equation*}
\dim_\integers(\ker(A))=0\text{ if }A\ne 0,
\end{equation*}
and of course
\begin{equation*}
\dim_\integers(\ker(A))=1\text{ if }A=0.
\end{equation*}

Since $\integers$ is commutative, every matrix $A$ can be replaced by
a diagonal matrix without changing the dimension of the kernel, and
this way the above calculation proves the Atiyah conjecture for
$\integers$.

Similar considerations for $\Gamma=\integers^n$ show that the Atiyah
conjecture here amounts to understanding the zeros of polynomials in
several variables. This created the slogan that the Atiyah conjecture
(and more generally $L^2$-cohomology) in a certain sense is
non-commutative algebraic geometry.

Exactly the same proof works if we replace the coefficient ring
$\integers$ by $\complexs$, or by any subring of $\complexs$.

This leads to the following algebraic generalization of the Atiyah
conjecture.
\begin{conjecture}\label{conj:general_algebraic_Atiyah} Fix a discrete 
  group $\Gamma$.
  Let $K$ be any subring of $\complexs$ which is closed under complex 
  conjugation. Let $\Gamma$ be a discrete group, and assume that $a\in 
  M(d\times d, K\Gamma)$.

  If $Fin^{-1}(\Gamma)$ is discrete, then $\dim_\Gamma(\ker(A))\in
  Fin^{-1}(\Gamma)$. 
\end{conjecture}

Since we can multiply any matrix with a non-zero constant (e.g.~a
common denominator of the finitely many non-zero coefficients), the
assertion of Conjecture \ref{conj:general_algebraic_Atiyah} for a ring 
$K\subset\complexs$ and its field of fractions is equivalent. In the
sequel, we will therefore usually assume that the subring $K$ is a field.

It is not clear, however, whether Conjecture
\ref{conj:general_algebraic_Atiyah} is equivalent to the original
geometric Atiyah conjecture \ref{conj:Atiyah} if $K\not\subset
\rationals$.

Observe that the proof of the Atiyah conjecture for $\integers$
extends from the ring of Laurent polynomials (with complex
coefficients) to the ring of meromorphic functions $\complexs$ without 
poles on $S^1$, which is also a subring of $L^\infty(S^1)\iso
\NeumannN\integers$. The question arises whether there are reasonable
generalizations similar to this ring for other groups. One possibility 
would be to look at infinite sums $a=\sum_{g\in\Gamma} \lambda_g g$,
where the coefficients $\lambda_g$ very rapidly tend to zero as
$g\to\infty$ (with respect to a suitable word length metric). (Observe 
that this is the case for the coefficients of the Laurent expansion of 
a meromorphic function on $\complexs$). Under suitable circumstances,
(convolution) multiplication with such an $a$  will indeed give rise
to a (very special) $\Gamma$-equivariant operator on $l^2\Gamma$, and
the question arises whether for the dimension of its kernel the
statement of the Atiyah conjecture holds. This was suggested by Nigel
Higson. It is quite distinct from Conjecture
\ref{conj:general_algebraic_Atiyah} in that it is analytic in flavor,
and no longer algebraic.

\subsection{Atiyah conjecture and zero divisors}

Among the most interesting observations about the Atiyah conjecture
are its strong connections to questions in algebra, in particular to
group rings.

Here, we address the following conjecture, the so called \emph{zero
  divisor conjecture}.
\begin{conjecture}\label{conj:zero_divisor}
  Let $\Gamma$ be a torsion-free discrete group and $K$ a subring of
  the complex numbers. Then, there are no non-trivial zero divisors in 
  the group ring $K\Gamma$, i.e.~if $a,b\in K\Gamma$ with $ab=0$ then
  either $a=0$ or $b=0$.
\end{conjecture}

This is one of the longstanding questions in the theory of group rings 
(which of course makes also sense for other coefficients rings $K$ and 
is studied also in this broader generality by ring theorists).

Observe that, if $g\in\Gamma$ is a torsion element, i.e.~$g\ne 1$ but
$g^n=1$ for some $n>0$, then $a=(1-g)$ and $b=1+g+\cdots+g^{n-1}$ are
two non-zero elements of $\integers\Gamma$ with $ab=0$.

It now turns out that the Atiyah conjecture implies the zero divisor
conjecture. More precisely:
\begin{theorem}\label{theo:Atiyah_implies_zero_divisor}
 Assume $\Gamma$ is a torsion-free discrete group, and
 $K\subset\complexs$ is a ring (closed under complex conjugation).

 If the statement of the algebraic Atiyah conjecture
 \ref{conj:general_algebraic_Atiyah} is true for every $A\in M(1\times 
 1,K\Gamma)$, then there are no non-trivial zero divisors in $K\Gamma$.
\end{theorem}
\begin{proof}
  Fix $a,b\in K\Gamma$ with $ab=0$. We have to show that either $a=0$
  or $b=0$. Now observe that 
  \begin{equation*}
a\in K\Gamma=M(1\times 1,K\Gamma)
\end{equation*}
is a
  $1$-by-$1$ matrix over $K\Gamma$. On the other hand, 
  \begin{equation*}
b\in
  K\Gamma\subset l^2\Gamma
\end{equation*}
can be considered to be an element of
  $l^2\Gamma$. And $0=ab$ is just the result of the action of the
  matrix $a$ on the $l^2$-function $b$. Therefore, $b\in \ker(a)$. Now 
  we know that $\dim_\Gamma(\ker(a))\in\integers$ because the
  Atiyah conjecture is true for the torsion-free group
  $\Gamma$. Evidently, 
  \begin{equation*}
\{0\}\subset \ker(a)\subset l^2\Gamma,
\end{equation*}
with
\begin{equation*}
0=\dim_\Gamma(\{0\})\quad \text{and}\quad1=\dim_\Gamma(l^2\Gamma).
\end{equation*}
Because of
  monotonicity, either $\dim_\Gamma(\ker(a))=0$ or
  $\dim_\Gamma(\ker(a))=1$. In the first case, because of faithfulness 
  $\ker(a)=\{0\}$ which implies $b=0$. In the second case,
  $\ker(a)=l^2\Gamma$ (again because of faithfulness) which means that 
  $a=0$. This proves the statement.

  Here, we used some of the properties of $\dim_\Gamma$ developed in
  Proposition \ref{prop:L2_dim_properties}.
\end{proof}

In Theorem \ref{theo:skew}, we will see that there is actually an even
stronger relation
between the Atiyah conjecture and the zero divisor conjecture for a
torsion free group $\Gamma$.

\subsection{Atiyah conjecture and calculations}

A second possible application of the Atiyah conjecture  could be the
explicit calculation
of $L^2$-Betti
numbers.

Here one would use that by now there are several approximation formulas
for $L^2$-Betti numbers. In particular, we have discussed one of these 
in Section \ref{sec:appr-l2-betti}. Obviously, if we know in advance
that the limit has to be an integer, this can make it much easier to
exactly compute the limit, in particular if (as is the case for some
of the approximation results) error bounds are available.
Up to now, however, to the authors knowledge this idea has not been used 
anywhere.

\subsection{The status of the Atiyah conjecture}

By now, the Atiyah conjecture is known for a reasonably large class of 
groups. We use the following definitions.

\begin{definition}\label{Def:Groups}\strut
  The class of \emph{elementary amenable groups} is the smallest class of groups which contains
  all abelian and all finite groups and is closed under extensions and
  directed unions. It is denoted $\amenableGroups$. Obviously, every
  elementary amenable group is amenable. Moreover, every nilpotent and 
  every solvable group is elementary amenable, as well as every group
  which is virtually solvable. (A group has \emph{virtually} a
  property P, if it contains a subgroup of finite index which actually 
  has property P).
\end{definition}

\begin{definition}\label{defextC}
  Let $\extendedC$ be the smallest non-empty class of groups such that:
  \begin{enumerate}
  \item \label{aex} If $G$ is torsion-free and $A$ is elementary
    amenable, and we have a projection
    $p\colon G\to A$ such that $p^{-1}(E)\in\extendedC$ for every
    finite subgroup $E$ of $A$, then $G\in\extendedC$.
   \item \label{sgr} $\extendedC$ is subgroup closed.
   \item\label{lim} Let $G_i\in \extendedC$ be a
  directed system of groups
    and $G$ its (direct or inverse) limit. Then $G\in\extendedC$.
\end{enumerate}
\end{definition}

\begin{definition}
  A \emph{directed system} of groups is a system of groups
  $G_i$, indexed by an index set $I$ with a partial ordering $<$, and
  either
  with a homomorphism $\phi_{ij}\colon G_i\to G_j$ if $i<j$ such that
  $\phi_{jk}\circ \phi_{ij}=\phi_{ik}$ if $i<j<k$ (then we take the
  direct limit), or with homomorphisms the other way around,
  i.e.~$\phi_{ij}\colon G_j\to G_i$ if $i<j$, with the corresponding
  compatibility condition (then we take the inverse limit).

  Directed means that to each $i,j\in I$ exists $k\in I$ with
  $i<k$ and $j<k$. The most obvious examples are systems of groups
  indexed by $\naturals$ with its usual ordering.
\end{definition}

Observe that $\extendedC$ contains
only torsion-free groups.

\begin{example}\label{ex:groups_in_extendedC}
  The class $\extendedC$ contains all torsion-free elementary amenable
  groups. It
  also contains all free groups and all  braid groups (compare Section 
  \ref{sec:atiy-conj-braid}). Moreover, $\extendedC$ is
  closed under direct sums, direct products, and free products.

  To see this, observe that clearly $\extendedC$ contains all
  elementary amenable groups, as long
  as they are torsion-free. Moreover, if $\Gamma$ contains a sequence
  of normal subgroups $\Gamma\supergroup \Gamma_1\supergroup
  \Gamma_2\supergroup\cdots$ with
  $\bigcap_{k\in\naturals}\Gamma_k=\{1\}$ and such that
  $\Gamma/\Gamma_k$ is torsion-free elementary amenable, then $\Gamma$ 
  is a subgroup of the inverse limit of the sequence of quotients, and 
  consequently belongs to $\extendedC$.

  In particular, every free group admits such a sequence in such a way 
  that the quotients are torsion-free nilpotent, and every braid group 
  admits such a sequence where the quotients are torsion-free and
  virtually nilpotent.
\end{example}

The following theorem is proved in \cite{Schick(2000)} and
\cite{Dodziuk-Linnell-Mathai-Schick-Yates(2001)}. 
\begin{theorem}\label{theo:best_result_about_AC}
  Set $K:=\overline{\rationals}$, the field of algebraic numbers (over 
  $\rationals$) in $\complexs$.

  If $\Gamma\in\extendedC$, then the Atiyah conjecture
  \ref{conj:general_algebraic_Atiyah} is true for $K\Gamma$.
\end{theorem}

To prove this, one only has to prove that the Atiyah conjecture is
preserved when passing to subgroups, doing extensions with
torsion-free elementary
amenable quotients, and under direct and inverse limits of directed
systems of groups in $\extendedC$. The former (almost) translates to
directed unions
of groups, and the latter to the case where a group $\Gamma$ has the
nested sequence of normal subgroups $\Gamma_k$ with trivial
intersection we have discussed previously.

We want to start with two rather elementary observations, which are
not very useful without examples of groups for which the Atiyah
conjecture is true, but which give rise to part of the statements of
Theorem \ref{theo:best_result_about_AC}.

\begin{proposition}\label{subgroups}
  If $G$ fulfills the  Atiyah conjecture, and if $U$ is a
  subgroup of $G$ with $Fin^{-1}(U)=Fin^{-1}(G)$, then
$U$ fulfills the Atiyah conjecture.
\end{proposition}
\begin{proof}
  This follows from the fact that the $U$-dimension of the kernel of
  a matrix over $K U$ acting on $l^2(U)^n$ coincides with the
  $G$-dimension of the
  same matrix, considered as an operator on $l^2(G)^n$
  (compare e.g.~\cite[3.1]{Schick(1998a)}), which follows from a
  simply diagonal
  decomposition argument.
\end{proof}

\begin{proposition}  \label{unions}
  Let $G$ be the directed union of groups $\{G_i\}_{i\in I}$ and
  assume that each
  $G_i$ fulfills the Atiyah conjecture. Then $G$ fulfills the Atiyah
  conjecture.
\end{proposition}
\begin{proof}
  A matrix over $KG$, having only finitely many non-trivial
  coefficients, already is a matrix over $KG_i$ for some $i$. The
  $G_i$-dimension and the $G$-dimension of the kernel of the
  matrix coincide, as in Proposition \ref{subgroups}, compare
  e.g.~\cite[3.1]{Schick(1998a)}. Note that $Fin^{-1}(G_i)$ is
  contained in $Fin^{-1}(G)$ for each $i\in I$.
\end{proof}

Besides of these two results, at the moment three rather different
methods are known to prove the 
Atiyah conjecture for certain groups.

\subsubsection{Fredholm modules and the Atiyah conjecture}

One of these methods, which one could call the method of the
``finite rank Fredholm
module'', was developed by Peter Linnell in \cite{Linnell(1993)} to
prove the Atiyah conjecture for free groups. It extends in an
ingenious way the method used by Connes to prove the trace conjecture
\ref{conj:trace_conjecture} 
(compare e.g.~\cite{Effros(1989)}) for the free group, and therefore
is related to the Baum-Connes conjecture. However, whereas the
KK-theory methods for Baum-Connes turn out to be quite flexible and
generalize to many other groups, compare e.g.~\cite{Schick(2001a)},
nobody so far was able to find a corresponding generalized approach to 
the Atiyah conjecture. Since the Atiyah conjecture for free groups
follows also from one of the other methods to be described later, we
don't discuss Linnell's original approach but refer instead to the
original article \cite{Linnell(1993)} and the review article
\cite{Linnell(1998)}. (Note, however, that the generalization to
products of free groups does not work as described in
\cite{Linnell(1998)}, since the proof of the basic Lemma \cite[Lemma
11.7]{Linnell(1998)} has a gap. One has to rely on a different method
for the proof, instead.)

\subsubsection{Atiyah conjecture and algebra}

We have shown in Theorem \ref{theo:Atiyah_implies_zero_divisor} that
the Atiyah conjecture for a torsion-free group implies the zero
divisor conjecture. The ring $K\Gamma$ evidently has no zero divisors, 
if it can be embedded into a skew field, i.e.~a (not necessary
commutative) ring where every non-zero element has a multiplicative
inverse. The optimal solution to the zero divisor question therefore
is to construct exactly this. It turns out that the Atiyah
conjecture provides us with such an embedding. We need the following definitions.

\begin{definition}\label{defunivU}
  Given the group von Neumann algebra $\NeumannN \Gamma$, we define
  $\universalU \Gamma$ to be the algebra of all unbounded operators
  affiliated to $\NeumannN \Gamma$ (compare e.g.~\cite[Section
  8]{Linnell(1998)}). That means a densely defined unbounded operator
  $D$ on $l^2\Gamma$ belongs to $\universalU\Gamma$ if and only if all
  its spectral projections belong to $\NeumannN\Gamma$.
  It is a classical fact that the ring $\universalU\Gamma$ is a
  (non-commutative) localization of $\NeumannN\Gamma$, which here
  means it is obtained from $\NeumannN\Gamma$ by inverting all
  non-zero divisors of $\NeumannN\Gamma$.

  Fix a subfield $K\subset \complexs$ which is closed under complex
  conjugation. We have to consider $K\Gamma\subset
  \NeumannN \Gamma\subset \universalU \Gamma$. Define $D_K\Gamma$ as the
  \emph{division closure}
  of $K\Gamma$ in $\universalU \Gamma$. By definition, this is the
  smallest subring of
  $\universalU \Gamma$ which contains $K\Gamma$ and which has the
  property that,
  whenever $x\in D_K\Gamma$ is invertible in $\universalU \Gamma$, then
  $x^{-1}\in  D_K\Gamma$.
\end{definition}

\begin{remark}
  Since the operators which belong to $\universalU\Gamma$ are only
  densely defined, one has to be careful when defining the sum or
  product of two such operators. This is done by first defining these
  operators on the obvious (common) domain, but then taking there
  closure, i.e.~to extend the domain of definition as far as
  possible. One has to check that this indeed gives a reasonable
  ring. This is a classical result which uses the $\Gamma$-dimension.
\end{remark}

\begin{example}
  Again, we turn to the example $\Gamma=\integers$. We have seen that, 
  via Fourier transform, $\NeumannN\Gamma$ becomes $L^\infty(S^1)$
  acting on $L^2(S^1)$ by pointwise multiplication. The ring
  $\universalU\integers$ becomes the ring of \emph{all} measurable
  functions on $S^1$, still acting by pointwise multiplication. These
  operators are in general unbounded and not defined on all of
  $L^2(S^1)$, because the product of an $L^2$-function with an arbitrary 
  measurable function belongs not necessarily to $L^2(S^1)$.
  It is not hard to show that every measurable function $f$ on $S^1$ is the
  quotient of two bounded functions $g,h\in L^\infty(S^1)$, $f=g/h$,
  where the set of zeros of $h$ has measure zero. This reflects the
  fact that $\universalU\Gamma$ is a localization of $\NeumannN\Gamma$.

  If $K\subset\complexs$, then $K\integers$ are the Laurent
  polynomials $K[z,z^{-1}]$ identified with functions on $S^1$ (by
  substituting $z\in S^1$ for the variable). In the same way,
  $D_K\Gamma$ is the field of rational functions $K(z)$, identifies
  with functions on $S^1$ by substituting $z\in S^1$ for the variable.
\end{example}

The (very strong) connection between the Atiyah conjecture and ring theoretic
properties of $D_K\Gamma$ is given by the following theorem.

\begin{theorem}\label{theo:skew}
  Let $\Gamma$ be a torsion-free group, and let $K$ be a subfield of
  $\complexs$ which is closed under complex conjugation.

 $K\Gamma$ fulfills the strong
  Atiyah conjecture in the sense of Conjecture \ref{conj:general_algebraic_Atiyah}
  if and only if the division closure $D_k\Gamma$ of $K\Gamma$ in
  $\universalU \Gamma$ is a skew field.
\end{theorem}

In other words, we have a canonical candidate $D_K\Gamma$ for a skew
field, into which $K\Gamma$ embeds, and this ring is a skew field if
and only if $K\Gamma$ satisfies the Atiyah conjecture.
\begin{proof}[Proof of Theorem \ref{theo:skew}.]
  If $D_K\Gamma$ is a skew field, then each matrix $A\in M(d\times
  d,K\Gamma)$ acts on $(D_K\Gamma)^d$, and its kernel is a finite
  dimensional vector space over the field $D_K\Gamma$. In particular, 
  its $D_K\Gamma$-dimension of course is an integer. Then, one can establish
  that the $\Gamma$-dimension of $\ker(A\colon
  (l^2\Gamma)^d\to(l^2\Gamma)^d)$ coincides with this
  dimension. Details are given in \cite[Lemma 3]{Schick(2000)}. 

  For the converse, given an element $0\ne a\in D_K\Gamma$ one can,
  using a matrix trick for division closures due to Cohn, produce a
  $d\times d$-matrix $A$ over $K\Gamma$ such that the
  $\Gamma$-dimension of $\ker(A)$ is evidently strictly smaller than
  $1$, and which is non-zero if and only if $a$ is invertible in
  $\universalU\Gamma$ (slogan: ``a non-trivial
  kernel of $a$ gives rise to a non-trivial kernel of $A$''). Because
  $\dim_\Gamma(\ker(A))\in\integers$ by assumption,
  $\dim_\Gamma(\ker(A))=0$, i.e.~$a$ is invertible in
  $\universalU\Gamma$. Since $D_K\Gamma$ is division closed, $a$ is
  invertible in $D_K\Gamma$, as well.
\end{proof}

This property allows to prove the Atiyah conjecture for the first
interesting class of groups (containing non-abelian groups), namely
the class of elementary amenable groups.

\begin{theorem}\label{theo:amext} Fix a subfield $K\subset \complexs$
  which is closed under complex conjugation.
  
 Let $1\to H\to G\to A\to 1$ be an exact sequence of
  groups. Assume that $G$ is torsion free and $A$
  is elementary amenable. For every finite subgroup
  $E\subgroup A$ let $H_E$ be the inverse
  image of $E$ in $G$. Assume for all
  finite subgroups
  $E\subgroup G$ that $K H_E$ fulfills
  the Atiyah
  conjecture \ref{conj:general_algebraic_Atiyah}. Then $KG$
  fulfills also
  the Atiyah conjecture.
\end{theorem}
\begin{proof}
  The proof is given in \cite{Linnell(1993)} for
  $K=\rationals$. Essentially the same proof works for arbitrary $K$,
  compare \cite[Proposition 3.1]{Schick(2000)}.
\end{proof}

\begin{corollary}\label{tfamext} Fix a subfield $K=\overline
  K\subset\complexs$.

 Suppose
  $H$ is torsion-free and $KH$ fulfills the Atiyah conjecture. If $G$ is an
  extension of $H$ with elementary amenable torsion-free quotient then
  $KG$ fulfills the Atiyah conjecture.

  In particular (with $H=1$) if $G$ is a torsion-free elementary
  amenable group
  then $KG$ satisfies the Atiyah conjecture.
\end{corollary}
\begin{proof}
  By assumption, the only finite subgroup of $G/H$ is the trivial
  group and the Atiyah conjecture is true for its inverse image $H$.
\end{proof}

\subsubsection{Atiyah conjecture and approximation}

Here, we describe the last method of proof for the Atiyah
conjecture. It is based on the approximation results of Section
\ref{sec:appr-l2-betti}.

\begin{theorem}
  Assume $\Gamma$ is a torsion-free discrete group with a nested
  sequence of normal subgroups $\Gamma\supergroup \Gamma_1\supergroup
  \Gamma_2\supergroup\cdots$ such that
  $\bigcap_{k\in\naturals}\Gamma/\Gamma_k=\{1\}$ and such that
  $\Gamma/\Gamma_k\in\RAgroups$ for each $k\in\naturals$. Moreover,
  assume that all the quotient groups $\Gamma/\Gamma_k$ are torsion-free
  and satisfy the Atiyah conjecture \ref{conj:Atiyah}.

  For example, all the groups $\Gamma/\Gamma_k$ might be torsion-free
  elementary amenable.

  Then $\Gamma$ also satisfies the Atiyah conjecture
  \ref{conj:Atiyah}.
\end{theorem}
\begin{proof}
  Given any $\Gamma$-covering $\overline X\to X$ (with a finite
  CW-complex $X$), we have to prove that $b^p_{(2)}(\overline
  X,\Gamma)\in\integers$ for each $p$.

  Now, the sequence of normal subgroups $\Gamma_k$ provides us with a
  sequences of normal coverings $X_k:=\overline X/\Gamma_k$ of $X$,
  with covering group $\Gamma/\Gamma_k$. Since
  $\Gamma/\Gamma_k\in\RAgroups$ for each $k\in\naturals$, by Theorem
  \ref{theo:general_approximation}
  \begin{equation*}
    b^p_{(2)}(\overline X,\Gamma) =\lim_{k\to\infty}
    b^p_{(2)}(X_k,\Gamma/\Gamma_k).
  \end{equation*}
  By assumption, each term on the right hand side is an integer,
  since the Atiyah conjecture holds for $\Gamma/\Gamma_k$. Since
  $\integers$ is discrete in $\reals$, the same will be true for its
  limit, and this is exactly what we have to prove.
\end{proof}

As observed above, this translates to a statement about the integral
group ring of $\Gamma$. To extend this result from the integral group
ring to more general
coefficient rings, which is interesting because of the algebraic
consequences, we have to generalize the approximation results of
Section \ref{sec:appr-l2-betti} to algebraic approximating results
for more general coefficient rings. In fact, we have the following
result of \cite{Dodziuk-Linnell-Mathai-Schick-Yates(2001)}.

\begin{theorem}\label{theo:algebraic_approximation}
  Assume $\Gamma$ is a discrete group with a nested
  sequence of normal subgroups $\Gamma\supergroup \Gamma_1\supergroup
  \Gamma_2\supergroup\cdots$ such that
  $\bigcap_{k\in\naturals}\Gamma/\Gamma_k=\{1\}$ and such that
  $\Gamma/\Gamma_k\in\RAgroups$ for each $k\in\naturals$.

  For example, all the groups $\Gamma/\Gamma_k$ might be
  elementary amenable, e.g.~finite.

  Let $\overline\rationals$ be the field of algebraic numbers (over
  $\rationals$) in $\complexs$. Assume $A\in M(d\times
  d,\overline\rationals\Gamma)$.

  The projection $\Gamma\to\Gamma/\Gamma_k$ extends canonically to the 
  group rings and to matrix rings over the group rings. Let $A_k\in
  M(d\times d,\overline\rationals\Gamma/\Gamma_k)$ be the image of $A$ 
  under this induced homomorphism.

  Then $A$ acts on $(l^2\Gamma)^d$ and $A_k$ acts on
  $(l^2\Gamma/\Gamma_k)^d$.

  The following approximation result for the kernels of these
  operators holds:
  \begin{equation*}
    \dim_\Gamma(\ker(A)) = \lim_{k\to\infty} \dim_{\Gamma/\Gamma_k} (\ker(A_k)).
  \end{equation*}

In particular, if all the groups $\Gamma/\Gamma_k$ are
torsion-free and $\overline\rationals[\Gamma/\Gamma_k]$ satisfies the
algebraic Atiyah conjecture \ref{conj:general_algebraic_Atiyah}, then, 
using the same argument as above, $\overline\rationals\Gamma$ also
satisfies the Atiyah conjecture \ref{conj:general_algebraic_Atiyah}. 
\end{theorem}

\begin{remark}
  The approximation result \ref{theo:general_approximation} is a
  special case, since if $\Delta$ is a matrix representative for a
  combinatorial Laplacian for
  $\overline X$, then $\Delta_k$, constructed as in Theorem
  \ref{theo:algebraic_approximation}, is a matrix representative for
  the corresponding combinatorial Laplacian of $X_k$.

  The proof uses the fact that one is working with algebraic
  coefficients. So far, no generalization to $\complexs\Gamma$ 
  has been obtained.

  However, it is a well known fact that, if
  $\overline\rationals\Gamma$ has no non-trivial zero divisors, then
  the same is true for $\complexs\Gamma$ (compare
  e.g.~\cite{Dodziuk-Linnell-Mathai-Schick-Yates(2001)}). Therefore,
  from the point of view of the zero divisor conjecture
  \ref{conj:zero_divisor}, there is no need to generalize Theorem
  \ref{theo:algebraic_approximation}. 
\end{remark}

\begin{remark}
  It is not hard to see that $\extendedC$ is contained in
  $\RAgroups$. Therefore, Theorem \ref{theo:algebraic_approximation}
  provides the last step for the proof of Theorem
  \ref{theo:best_result_about_AC}. 
\end{remark}

\subsubsection{Atiyah conjecture for braid groups}
\label{sec:atiy-conj-braid}

\begin{definition}
  Fix $n\in\naturals$.  A \emph{braid with $n$-strings} is an
  embedding 
  \begin{equation*}
\phi\colon\{1,\dots,n\}\times [0,1]\to\complexs\times
[0,1]
\end{equation*}
such that $\phi(p,0)=(p,0)$ and
$\phi(p,1)\in \{1,\dots,n\}\times\{1\}$ for $p=1,\dots,n$. Two braid
are considered
equal if they are isotopic where the isotopy fixes the top and the
bottom.

Isotopy classes of braids from a group by stacking two braids
together, the so called \emph{Artin braid group $B_n$}. Note that the
$p$-th string is not necessarily stacked on the $p$-th string, since
the $p$-th string might lead from $(p,0)$ to $(\sigma(p),1)$ for some
permutation $\sigma$ of $\{1,\dots,n\}$. We have 
to account for this when we define the ``stacked'' map
$\{1,\dots,n\}\times [0,1]\to \complexs\times[0,1]$.

The braid group $B_n$ contains a normal subgroup $P_n$, the \emph{pure 
  braid group}, where we require $\phi(p,1)=(p,1)$ for each
$p\in\{1,\dots,n\}$. The quotient $B_n/P_n$ is the symmetric group
$S_n$ of permutations of $\{1,\dots,n\}$, where the image permutation
is given as above. 
\end{definition}

In Example \ref{ex:groups_in_extendedC}, we assert that all the braid
groups belong to $\extendedC$. Indeed, every braid group
$B_n$ has a nested sequence of normal subgroups $B_n\supergroup
P_{1,n}\supergroup P_{2,n}\supergroup \cdots$ with $B_n/P_{k,n}$
torsion-free elementary amenable for each $k$ and such that
$\bigcap_{k\in\integers} P_{k,n}=\{1\}$. The
corresponding result for the pure braid groups is proved in
\cite[Theorem 2.6]{Falk-Randell(1988)}. However, to extend such a
result from a subgroup of finite index to a bigger group is highly
non-trivial and in general not possible. Indeed, for the full braid
group it was conjectured for a while that it has no non-trivial
torsion-free quotients at all, opposite to what we need. Using certain
totally disconnected completions of the groups involved, and
cohomology of these completions (Galois cohomology, which takes the
topology of the completions into account) the above result is proved
in \cite{Linnell-Schick(2000)}. Actually, it is proved there that
every torsion-free finite extension of $P_n$ has a sequence of
subgroups as above, and therefore belongs to $\extendedC$. This paper also contains
generalizations, where the pure braid groups are replaced by other
kinds of groups, still with the result that the property to belong to
$\extendedC$ passes to finite extensions (as long as they are
torsion-free). In \cite{Linnell-Schick(2000c)}, it is shown how this
applies to fundamental groups of certain complements of links in
$\reals^3$ (a link is an embedding of the disjoint union of finitely
many circles).

The proof of the Baum-Connes conjecture for the full braid group
\cite{Schick(2000a)} mentioned in Section
\ref{sec:status-BC-conjecture} is
based on the same results.

\subsection{Atiyah conjecture for groups with torsion}

The Atiyah conjecture has also been obtained for many groups $\Gamma$ with
torsion, as long as $Fin^{-1}(\Gamma)$ is a discrete subset of
$\reals$. We are not discussing them here because of lack of space and 
time, and because there is no zero divisor conjecture for groups with
torsion. The ring $D_K\Gamma$ also exists for groups with torsion. It
can not be a skew field but, under the assumption that
$Fin^{-1}(\Gamma)$ is discrete, it often turns out to be a semi-simple 
Artinian ring. More details can be found e.g.~in the original sources
\cite{Linnell(1998),Linnell(1993),Schick(2000)}.

\bibliographystyle{plain}
\bibliography{script}

\begin{thebibliography}{10}

\bibitem{Atiyah(1976)}
M.~F. Atiyah.
\newblock Elliptic operators, discrete groups and von {N}eumann algebras.
\newblock In {\em Colloque ``Analyse et Topologie'' en l'Honneur de Henri
  Cartan (Orsay, 1974)}, pages 43--72. Ast\'erisque, No. 32--33. Soc. Math.
  France, Paris, 1976.

\bibitem{Atiyah-Collected2}
Michael Atiyah.
\newblock {\em Collected works. {V}ol. 2}.
\newblock The Clarendon Press Oxford University Press, New York, 1988.
\newblock $K$-theory.

\bibitem{Ballmann-Bruening(2000)}
W.~Ballmann and J.~Br\"uning.
\newblock On the spectral theory of manifolds with cusps.
\newblock SFB256-Preprint 689, Bonn 2000.

\bibitem{Baum-Connes-Higson(1994)}
Paul Baum, Alain Connes, and Nigel Higson.
\newblock Classifying space for proper actions and ${K}$-theory of group
  ${C}\sp \ast$-algebras.
\newblock In {\em $C\sp \ast$-algebras: 1943--1993 (San Antonio, TX, 1993)},
  pages 240--291. Amer. Math. Soc., Providence, RI, 1994.

\bibitem{Baum-Schick(2001)}
Paul Baum and Thomas Schick.
\newblock Equivariant {K}-homology as via equivariant $({M,E},\phi)$-theory.
\newblock in preparation (2001).

\bibitem{Beguni-Bettaieb-Valette(1999)}
C{\'e}dric B{\'e}guin, Hela Bettaieb, and Alain Valette.
\newblock ${K}$-theory for ${C}\sp \ast$-algebras of one-relator groups.
\newblock {\em $K$-Theory}, 16(3):277--298, 1999.

\bibitem{Berline-Getzler-Vergne(1992)}
Nicole Berline, Ezra Getzler, and Mich{\`e}le Vergne.
\newblock {\em Heat kernels and {D}irac operators}.
\newblock Springer-Verlag, Berlin, 1992.

\bibitem{Blackadar(1998)}
Bruce Blackadar.
\newblock {\em ${K}$-theory for operator algebras}.
\newblock Cambridge University Press, Cambridge, second edition, 1998.

\bibitem{Borel(1985)}
A.~Borel.
\newblock The ${L}\sp 2$-cohomology of negatively curved {R}iemannian symmetric
  spaces.
\newblock {\em Ann. Acad. Sci. Fenn. Ser. A I Math.}, 10:95--105, 1985.

\bibitem{Botvinnik-Gilkey-Stolz(1995)}
Boris Botvinnik, Peter Gilkey, and Stephan Stolz.
\newblock The {G}romov-{L}awson-{R}osenberg conjecture for groups with periodic
  cohomology.
\newblock {\em J. Differential Geom.}, 46(3):374--405, 1997.

\bibitem{Broecker-Dieck(1995)}
Theodor Br{\"o}cker and Tammo tom Dieck.
\newblock {\em Representations of compact {L}ie groups}.
\newblock Springer-Verlag, New York, 1995.
\newblock Translated from the German manuscript, Corrected reprint of the 1985
  translation.

\bibitem{Cao-Xavier(2001)}
J.~Cao and F.~Xavier.
\newblock {K}{\"a}hler parabolicity and the {E}uler number of compact manifolds
  of non-positive sectional curvature.
\newblock {\em Math. Annalen}, pages 483--491, 2001.

\bibitem{Chern(1968)}
S.~S. Chern.
\newblock {\em Minimal submanifolds in a {R}iemannian manifold}.
\newblock Univ. of Kansas Lawrence, Kan., 1968.

\bibitem{Davis-Lueck(1998)}
James~F. Davis and Wolfgang L{\"u}ck.
\newblock Spaces over a category and assembly maps in isomorphism conjectures
  in ${K}$- and ${L}$-theory.
\newblock {\em $K$-Theory}, 15(3):201--252, 1998.

\bibitem{Dicks-Schick(2001)}
Warren Dicks and Thomas Schick.
\newblock The spectral measure of certain elements of the group ring of a
  wreath product.
\newblock Preprint-Series SFB 478, M\"unster (2001), to appear in Geometriae
  Dedicata, available via
  \href{http://www.math.uni-muenster.de/u/schickt/publ/}{http://www.math.uni-m%
uenster.de/u/schickt/publ/}.

\bibitem{Dodziuk(1977)}
Jozef Dodziuk.
\newblock de {R}ham-{H}odge theory for ${L}\sp{2}$-cohomology of infinite
  coverings.
\newblock {\em Topology}, 16(2):157--165, 1977.

\bibitem{Dodziuk(1979)}
Jozef Dodziuk.
\newblock ${L}\sp{2}$\ harmonic forms on rotationally symmetric {R}iemannian
  manifolds.
\newblock {\em Proc. Amer. Math. Soc.}, 77(3):395--400, 1979.

\bibitem{Dodziuk-Linnell-Mathai-Schick-Yates(2001)}
Jozef Dodziuk, Peter Linnell, Varghese Mathai, Thomas Schick, and Stuart Yates.
\newblock Approximating ${L}^2$-invariants, and the {A}tiyah conjecture.
\newblock preprint (2001).

\bibitem{Donnelly-Xavier(1984)}
Harold Donnelly and Frederico Xavier.
\newblock On the differential form spectrum of negatively curved {R}iemannian
  manifolds.
\newblock {\em Amer. J. Math.}, 106(1):169--185, 1984.

\bibitem{Eckmann(2000)}
Beno Eckmann.
\newblock Introduction to $l\sb 2$-methods in topology: reduced $l\sb
  2$-homology, harmonic chains, $l\sb 2$-{B}etti numbers.
\newblock {\em Israel J. Math.}, 117:183--219, 2000.
\newblock Notes prepared by Guido Mislin.

\bibitem{Effros(1989)}
Edward~G. Effros.
\newblock Why the circle is connected: an introduction to quantized topology.
\newblock {\em Math. Intelligencer}, 11(1):27--34, 1989.

\bibitem{Falk-Randell(1988)}
Michael Falk and Richard Randell.
\newblock Pure braid groups and products of free groups.
\newblock In {\em Braids (Santa Cruz, CA, 1986)}, pages 217--228. Amer. Math.
  Soc., Providence, RI, 1988.

\bibitem{Federer(1969)}
H.~Federer.
\newblock {\em Geometric measure theory}, volume 153 of {\em Grundlehren der
  math. Wissenschaften}.
\newblock Springer, 1969.

\bibitem{Grigorchuk-Linnell-Schick-Zuk(2000)}
Rostislav~I. Grigorchuk, Peter Linnell, Thomas Schick, and Andrzej {\.Z}uk.
\newblock On a question of {A}tiyah.
\newblock {\em C. R. Acad. Sci. Paris S\'er. I Math.}, 331(9):663--668, 2000.

\bibitem{Grigurchuk-Zuk(1999)}
Rostislav~I. Grigorchuk and Andrzej {\.Z}uk.
\newblock On the asymptotic spectrum of random walks on infinite families of
  graphs.
\newblock In {\em Random walks and discrete potential theory (Cortona, 1997)},
  pages 188--204. Cambridge Univ. Press, Cambridge, 1999.

\bibitem{Gromov(1991)}
M.~Gromov.
\newblock {K}\"ahler hyperbolicity and ${L}\sb 2$-{H}odge theory.
\newblock {\em J. Differential Geom.}, 33(1):263--292, 1991.

\bibitem{Gromov-Lawson(1983)}
M.~Gromov and H.B. Lawson.
\newblock Positive scalar curvature and the dirac operator on complete
  {R}iemannian manifolds.
\newblock {\em Publ. Math. IHES}, 58:83--196, 1983.

\bibitem{Hambleton-Pedersen(2001)}
Ian Hambleton and Erik Pedersen.
\newblock Identifying assembly maps in {K}- and {L}-theory.
\newblock preprint (2001), available via
  \href{http://www.math.binghamton.edu/erik/}{http://www.math.binghamton.edu/e%
rik/}.

\bibitem{Higson-Kasparov(1997)}
Nigel Higson and Gennadi Kasparov.
\newblock Operator ${K}$-theory for groups which act properly and isometrically
  on {H}ilbert space.
\newblock {\em Electron. Res. Announc. Amer. Math. Soc.}, 3:131--142
  (electronic), 1997.

\bibitem{Higson-Lafforgue-Skandalis(2001)}
Nigel Higson, Vincent Lafforgue, and George Skandalis.
\newblock Counterexamples to the {B}aum-{C}onnes conjecture.
\newblock preprint, Penn State University, 2001, available via
  \href{http://www.math.psu.edu/higson/research.html}{http://www.math.psu.edu/%
higson/research.html}.

\bibitem{Higson-Roe(2000)}
Nigel Higson and John Roe.
\newblock Amenable group actions and the {N}ovikov conjecture.
\newblock {\em J. Reine Angew. Math.}, 519:143--153, 2000.

\bibitem{Higson-Roe(2001)}
Nigel Higson and John Roe.
\newblock {\em {Analytic {K}-homology.}}
\newblock {Oxford Mathematical Monographs}. Oxford University Press, Oxford,
  2001.

\bibitem{Hitchin(1974b)}
N.~Hitchin.
\newblock Harmonic spinors.
\newblock {\em Advances in Mathem.}, 14:1--55, 1974.

\bibitem{Joachim-Schick(2000)}
Michael Joachim and Thomas Schick.
\newblock Positive and negative results concerning the
  {G}romov-{L}awson-{R}osenberg conjecture.
\newblock In {\em Geometry and topology: Aarhus (1998)}, pages 213--226. Amer.
  Math. Soc., Providence, RI, 2000.

\bibitem{Jost-Xin(2000)}
J.~Jost and Y.~L. Xin.
\newblock Vanishing theorems for ${L}\sp 2$-cohomology groups.
\newblock {\em J. Reine Angew. Math.}, 525:95--112, 2000.

\bibitem{Jost-Zuo(2000)}
J{\"u}rgen Jost and Kang Zuo.
\newblock Vanishing theorems for ${L}\sp 2$-cohomology on infinite coverings of
  compact {K}\"ahler manifolds and applications in algebraic geometry.
\newblock {\em Comm. Anal. Geom.}, 8(1):1--30, 2000.

\bibitem{Julg-Kasparov(1995)}
Pierre Julg and Gennadi Kasparov.
\newblock Operator ${K}$-theory for the group ${\rm {s}{u}}(n,1)$.
\newblock {\em J. Reine Angew. Math.}, 463:99--152, 1995.

\bibitem{Kasparov(1995)}
G.~G. Kasparov.
\newblock ${K}$-theory, group ${C}\sp *$-algebras, and higher signatures
  (conspectus).
\newblock In {\em Novikov conjectures, index theorems and rigidity, Vol.\ 1
  (Oberwolfach, 1993)}, pages 101--146. Cambridge Univ. Press, Cambridge, 1995.

\bibitem{Kasparov-Skandalis(1991)}
G.~G. Kasparov and G.~Skandalis.
\newblock Groups acting on buildings, operator ${K}$-theory, and {N}ovikov's
  conjecture.
\newblock {\em $K$-Theory}, 4(4):303--337, 1991.

\bibitem{Kasparov-Skandalis(1994)}
Guennadi Kasparov and Georges Skandalis.
\newblock Groupes ``boliques'' et conjecture de {N}ovikov.
\newblock {\em C. R. Acad. Sci. Paris S\'er. I Math.}, 319(8):815--820, 1994.

\bibitem{Kazdan-Warner(1975)}
Jerry~L. Kazdan and F.~W. Warner.
\newblock Prescribing curvatures.
\newblock In {\em Differential geometry (Proc. Sympos. Pure Math., Vol. XXVII,
  Stanford Univ., Stanford, Calif., 1973), Part 2}, pages 309--319. Amer. Math.
  Soc., Providence, R.I., 1975.

\bibitem{Kwasik-Schultz(1990)}
S.~Kwasik and R.~Schultz.
\newblock Positive scalar curvature and periodic fundamental groups.
\newblock {\em Comm. Math. Helv.}, 65:271--286, 1990.

\bibitem{Lafforgue(1999)}
Vincent Lafforgue.
\newblock Compl\'ements \`a la d\'emonstration de la conjecture de
  {B}aum-{C}onnes pour certains groupes poss\'edant la propri\'et\'e ({T}).
\newblock {\em C. R. Acad. Sci. Paris S\'er. I Math.}, 328(3):203--208, 1999.

\bibitem{Lawson-Michelsohn(1989)}
H.~Blaine Lawson, Jr. and Marie-Louise Michelsohn.
\newblock {\em Spin geometry}.
\newblock Princeton University Press, Princeton, NJ, 1989.

\bibitem{Lichnerowicz(1963)}
{Lichnerowicz, A.}
\newblock {Spineurs harmoniques}.
\newblock {\em C.R. Acad. Sci. Paris, s\'eries 1}, 257:7--9, 1963.
\newblock Zentralblatt 136.18401.

\bibitem{Linnell-Schick(2000)}
P.~Linnell and T.~Schick.
\newblock Finite group extensions and the {A}tiyah conjecture.
\newblock preliminary version.

\bibitem{Linnell-Schick(2000c)}
P.~Linnell and T.~Schick.
\newblock {G}alois cohomology of completed link groups.
\newblock Preprint, M{\"u}nster, 2000.

\bibitem{Linnell(1993)}
Peter~A. Linnell.
\newblock Division rings and group von {N}eumann algebras.
\newblock {\em Forum Math.}, 5(6):561--576, 1993.

\bibitem{Linnell(1998)}
Peter~A. Linnell.
\newblock Analytic versions of the zero divisor conjecture.
\newblock In {\em Geometry and cohomology in group theory (Durham, 1994)},
  pages 209--248. Cambridge Univ. Press, Cambridge, 1998.

\bibitem{Lueck(1997)}
W.~L{\"u}ck.
\newblock ${L}^2$-invariants of regular coverings of compact manifolds and
  ${CW}$-complexes.
\newblock to appear in ``Handbook of Geometry'', Elsevier, 1998.

\bibitem{Lueck(2001)}
W.~L{\"u}ck.
\newblock ${L}^2$-invariants and ${K}$-theory.
\newblock to be published by Springer, 2001.

\bibitem{Lueck(2001a)}
W.~L{\"u}ck.
\newblock The relation between the {B}aum-{C}onnes conjecture and the trace
  conjecture.
\newblock Preprintreihe SFB 478 --- Geometrische Strukture in der Mathematik,
  Heft 151, 2001.

\bibitem{Lueck-Schick(2001)}
W.~L\"uck and T.~Schick.
\newblock Approximating ${L}^2$-signatures bey their compact analogues.
\newblock Preprintreihe SFB 478 --- Geometrische Strukture in der Mathematik,
  Heft 190, 2001, available via
  href{http://arxiv.org/abs/math.GT/0110328}{http://arxiv.org/abs/math.GT/0110%
328}.

\bibitem{Lueck(1998a)}
Wolfgang L{\"u}ck.
\newblock Dimension theory of arbitrary modules over finite von {N}eumann
  algebras and ${L}\sp 2$-{B}etti numbers. {I}. {F}oundations.
\newblock {\em J. Reine Angew. Math.}, 495:135--162, 1998.

\bibitem{Lueck(1998b)}
Wolfgang L{\"u}ck.
\newblock Dimension theory of arbitrary modules over finite von {N}eumann
  algebras and ${L}\sp 2$-{B}etti numbers. {I}{I}. {A}pplications to
  {G}rothendieck groups, ${L}\sp 2$-{E}uler characteristics and {B}urnside
  groups.
\newblock {\em J. Reine Angew. Math.}, 496:213--236, 1998.

\bibitem{Mathai(1999)}
Varghese Mathai.
\newblock ${K}$-theory of twisted group ${C}\sp *$-algebras and positive scalar
  curvature.
\newblock In {\em Tel Aviv Topology Conference: Rothenberg Festschrift (1998)},
  pages 203--225. Amer. Math. Soc., Providence, RI, 1999.

\bibitem{Mineyev-Yu(2001)}
Igor Mineyev and Guoliang Yu.
\newblock The {B}aum-{C}onnes conjecture for hyperbolic groups.
\newblock preprint 2001, available via
  \href{http://www.math.usouthal.edu/~mineyev/math/}{http://www.math.usouthal.%
edu/~mineyev/math/}.

\bibitem{Mislin(2001)}
Guido Mislin.
\newblock Equivariant {K}-homology of the classifying space for proper actions.
\newblock Lecture notes, in preparation (2001).

\bibitem{Morgan(1988)}
F.~Morgan.
\newblock {\em Geometric measure theory, a beginner's guide}.
\newblock Academic Press, 1988.

\bibitem{Olbricht(2000)}
M.~Olbricht.
\newblock ${L}^2$-invariants of locally symmetric spaces.
\newblock preprint, G\"ottingen, 2000.

\bibitem{Oyono(1998)}
Herv{\'e} Oyono-Oyono.
\newblock La conjecture de {B}aum-{C}onnes pour les groupes agissant sur les
  arbres.
\newblock {\em C. R. Acad. Sci. Paris S\'er. I Math.}, 326(7):799--804, 1998.

\bibitem{Palmer(1970)}
T.~W. Palmer.
\newblock Real ${C}\sp*$-algebras.
\newblock {\em Pacific J. Math.}, 35:195--204, 1970.

\bibitem{Rosenberg(1983)}
J.~Rosenberg.
\newblock ${C}^*$-algebras, positive scalar curvature, and the {N}ovikov
  conjecture.
\newblock {\em Publ. Math. IHES}, 58:197--212, 1983.

\bibitem{Rosenberg(1987)}
J.~Rosenberg.
\newblock ${C}^*$-algebras, positive scalar curvature, and the {N}ovikov
  conjecture iii.
\newblock {\em Topology}, 25:319--336, 1986.

\bibitem{Rosenberg-Stolz(1995)}
J.~Rosenberg and S.~Stolz.
\newblock The ``stable'' version of the {G}romov-{L}awson conjecture.
\newblock In {\em Proc. of the Cech Centennial Homotopy Theory Conference},
  volume 181 of {\em Contemporary Mathematics}, pages 405--418, 1995.

\bibitem{Rosenberg(1995)}
Jonathan Rosenberg.
\newblock Analytic {N}ovikov for topologists.
\newblock In {\em Novikov conjectures, index theorems and rigidity, Vol.\ 1
  (Oberwolfach, 1993)}, pages 338--372. Cambridge Univ. Press, Cambridge, 1995.

\bibitem{Roy(1999)}
Ranja Roy.
\newblock The trace conjecture---a counterexample.
\newblock {\em $K$-Theory}, 17(3):209--213, 1999.

\bibitem{Schick(2000a)}
Thomas Schick.
\newblock Finite group extensions and the {B}aum-{C}onnes conjecture.
\newblock preprint M\"unster 2000, available via
  \href{http://www.math.uni-muenster.de/u/schickt/publ}{http://www.math.uni-mu%
enster.de/u/schickt/publ}.

\bibitem{Schick(1996)}
Thomas Schick.
\newblock {\em Analysis on $\partial$-manifolds of bounded geometry, {H}odge-de
  {R}ham isomorphism and ${L}^2$-index theorem}.
\newblock Shaker Verlag, Aachen, 1996.
\newblock Dissertation, Johannes Gutenberg-Universit{\"a}t Mainz.

\bibitem{Schick(2000)}
Thomas Schick.
\newblock Integrality of ${L}\sp 2$-{B}etti numbers.
\newblock {\em Math. Ann.}, 317(4):727--750, 2000.

\bibitem{Schick(1998a)}
Thomas Schick.
\newblock ${L}\sp 2$-determinant class and approximation of ${L}\sp 2$-{B}etti
  numbers.
\newblock {\em Trans. Amer. Math. Soc.}, 353(8):3247--3265 (electronic), 2001.

\bibitem{Schick(2001)}
Thomas Schick.
\newblock ${L}\sp 2$-index theorem for elliptic differential boundary problems.
\newblock {\em Pacific J. Math.}, 197(2):423--439, 2001.

\bibitem{Schick(2001a)}
Thomas Schick.
\newblock The trace on the ${K}$-theory of group ${C}\sp *$-algebras.
\newblock {\em Duke Math. J.}, 107(1):1--14, 2001.

\bibitem{Schick-Stolz(2000)}
Thomas Schick and Stephan Stolz.
\newblock A counterexample to the {G}romov-{L}awson conjecture with
  torsion-free fundamental group.
\newblock Preliminary version.

\bibitem{Schoen-Yau(1979a)}
R.~Schoen and S.-T. Yau.
\newblock Existence of incompressible minimal surfaces and the topology of
  three dimensional manifolds with non-negative scalar curvature.
\newblock {\em Annals of Math.}, 110:127--142, 1979.
\newblock MR 81k:5802.

\bibitem{Schoen-Yau(1979)}
R.~Schoen and S.T. Yau.
\newblock On the structure of manifolds with positive scalar curvature.
\newblock {\em manuscr. mathem.}, 28:159--183, 1979.

\bibitem{Schroeder(1993)}
Herbert Schr{\"o}der.
\newblock {\em ${K}$-theory for real ${C}\sp *$-algebras and applications}.
\newblock Longman Scientific \& Technical, Harlow, 1993.

\bibitem{Shubin(1987)}
M.~A. Shubin.
\newblock {\em Pseudodifferential operators and spectral theory}.
\newblock Springer-Verlag, Berlin, 1987.
\newblock Translated from the Russian by Stig I. Andersson.

\bibitem{Smale(1993)}
Nathan Smale.
\newblock Generic regularity of homologically area minimizing hypersurfaces in
  eight-dimensional manifolds.
\newblock {\em Comm. Anal. Geom.}, 1(2):217--228, 1993.

\bibitem{Stolz(1992)}
Stephan Stolz.
\newblock Simply connected manifolds of positive scalar curvature.
\newblock {\em Annals of Math.}, 136:511--540, 1992.

\bibitem{Stolz(1994)}
Stephan Stolz.
\newblock Splitting certain ${M}\, {\rm {s}pin}$-module spectra.
\newblock {\em Topology}, 33(1):159--180, 1994.

\bibitem{Stolz(1998)}
Stephan Stolz.
\newblock Concordance classes of positive scalar curvature metrics.
\newblock preprint, Notre Dame, available via
  \href{http://www.nd.edu/~stolz/concordance.ps}{http://www.nd.edu/~stolz/conc%
ordance.ps}, 1998.

\bibitem{Suslin-Wodzicki(1992)}
Andrei~A. Suslin and Mariusz Wodzicki.
\newblock Excision in algebraic ${K}$-theory.
\newblock {\em Ann. of Math. (2)}, 136(1):51--122, 1992.

\bibitem{Swan(1962)}
Richard~G. Swan.
\newblock Vector bundles and projective modules.
\newblock {\em Trans. Amer. Math. Soc.}, 105:264--277, 1962.

\bibitem{Dieck(1972)}
Tammo tom Dieck.
\newblock Orbittypen und \"aquivariante {H}omologie. {I}.
\newblock {\em Arch. Math. (Basel)}, 23:307--317, 1972.

\bibitem{Valette(2001)}
A.~Valette.
\newblock Introduction to the {B}aum-{C}onnes conjecture.
\newblock preprint 2001, to appear as ETHZ lecture note, published by
  Birkh\"auser.

\bibitem{Valette(2000)}
A.~Valette.
\newblock On the {B}aum-{C}onnes assembly map for discrete groups.
\newblock unpublished preprint, perhaps to appear as appendix to ``Introduction
  to the Baum-Connes conjecture'', to appear as ETHZ lecture note, published by
  Birkh\"auser.

\bibitem{Wegge-Olsen(1993)}
N.~E. Wegge-Olsen.
\newblock {\em ${K}$-theory and ${C}\sp *$-algebras}.
\newblock The Clarendon Press Oxford University Press, New York, 1993.
\newblock A friendly approach.

\end{thebibliography}

\end{document}